\documentclass[10pt]{article}

\usepackage{amsmath}
\usepackage{amsmath, amsthm, amsfonts, amssymb, color}
\usepackage{mathrsfs}
\usepackage{latexsym,bm}
\usepackage{graphicx}
\usepackage[toc,page]{appendix}

\newtheorem{theorem}{Theorem}[section]
\newtheorem{proposition}[theorem]{Proposition}

\newtheorem{lemma}[theorem]{Lemma}

\newtheorem{prop}[theorem]{Proposition}
\newtheorem{remark}[theorem]{Remark}
\newtheorem{example}[theorem]{Example}

\def\FF{\mathcal F}

\def\eps{\varepsilon}

\definecolor{wco}{rgb}{0.5,0.2,0.3}

\newcommand{\p}{\mathrm{P}}
\newcommand{\q}{\mathrm{Q}}
\newcommand{\pp}{\mathbb{P}}
\newcommand{\qq}{\mathbb{Q}}
\newcommand{\mc}{\mathcal{M}_{c}(E)}
\newcommand{\mf}{\mathcal{M}_{F}(E)}
\newcommand{\mof}{\mathcal{M}^+_{F} (E)}
\newcommand{\mh}{\mathcal{M}^{h}_{F}(E)}

\newcommand{\mw}{\mathcal{M}^{w}_{F}(E)}
\newcommand{\B}{\mathcal{B}}
\newcommand{\R}{\mathbb{R}}
\newcommand{\F}{\mathcal{F}}
\newcommand{\K}{\kappa}
\newcommand{\ia}{\mathbf{a}}
\newcommand{\ib}{\mathbf{b}}
\numberwithin{equation}{section} 

\begin{document}

\allowdisplaybreaks

\title{\bf Skeleton decomposition and law of large numbers for supercritical superprocesses}
\author{{\bf Zhen-Qing Chen}\thanks{The research partially supported by  Simons Foundation grant 520542 and
NNSFC  Grant 11731009.},
\quad {\bf Yan-Xia Ren}
\thanks{The research of this author is supported by NNSFC (Grant No.  11731009 and 11671017).}
 \quad  \hbox{and}\quad
{\bf Ting Yang}\footnote{Corresponding author. The research of this author is supported by NNSF of China (Grant No. 11501029).}}
\date{}
\maketitle

\begin{abstract}
The goal of this paper has two-folds. First, we establish
skeleton and spine decompositions for superprocesses
whose underlying processes are general symmetric Hunt processes.
Second, we use
these decompositions to obtain
weak and strong law of large numbers for supercritical superprocesses where the spatial motion is a symmetric Hunt process on a locally compact metric space $E$ and the branching mechanism takes the form
$$\psi_{\beta}(x,\lambda)=-\beta(x)\lambda+\alpha(x)\lambda^{2}+\int_{(0,{\infty})}(e^{-\lambda y}-1+\lambda y)\pi(x,dy)$$
with $\beta\in\mathcal{B}_{b}(E)$, $\alpha\in \mathcal{B}^{+}_{b}(E)$ and $\pi$ being a kernel from $E$ to $(0,{\infty})$ satisfying $\sup_{x\in E}\int_{(0,{\infty})} (y\wedge y^{2}) \pi(x,dy)<{\infty}$. The limit theorems are established under the assumption that an associated Schr\"{o}dinger operator has a spectral gap.   Our results cover
many interesting examples of superprocesses,
including super Ornstein-Uhlenbeck process and super stable-like process.
The strong law of large numbers for supercritical superprocesses are obtained under the assumption
that the strong law of large numbers for an associated supercritical
branching Markov process holds along a
discrete sequence of times,
extending an earlier result of Eckhoff, Kyprianou and Winkel \cite{EKW}
for superdiffusions to a large class of superprocesses.
The key for such a result is due to the skeleton decomposition of
 superprocess, which represents a superprocess as an immigration process along a supercritical branching Markov process.
\end{abstract}

\medskip


\medskip


\section{Introduction}

 Recently there is a lot of work
(see \cite{CRW,EKW,E,EK,ET,EW,KR,LRS,W} and the references therein)
on
limit theorems for superprocesses using the principal eigenvalue and ground state of the
corresponding Schr\"{o}dinger equations. For superdiffusions, weak (convergence in law or in probability) and strong (almost sure convergence) laws of large numbers have been established successively in
\cite{EKW,E,EK,ET,EW,LRS,W}. That the underlying process is a diffusion plays an essential role
in these papers.
We refer to \cite{EKW} for a survey on the recent developments of
limit theorems for superdiffusions. Unlike the case for superdiffusions, there is much less work on limit theorems for superprocesses when the spatial motion is a general Hunt process.
As far as we know, \cite{CRW} is the first paper to establish
strong law of large numbers (SLLN in brief)
for superprocesses. The spatial motions in \cite{CRW} are symmetric Hunt processes
which can have discontinuous sample paths. Amongst other assumptions, a spectral gap condition was used in \cite{CRW} to obtain a Poincar\'{e} inequality,
which is
the main ingredient in the proof of almost sure convergence along lattice times. Later \cite{KR} used a Fourier analysis approach to establish
a SLLN for super-stable processes with spatially independent quadratic branching mechanisms. Very recently,
SLLN is obtained in \cite{CRSZ}
for a class of superprocesses where the spatial motion can be a non-symmetric Hunt process. In \cite{CRSZ}, the mean semigroup of the superprocess is assumed to be compact
in $L^{2}$ space
and the branching mechanism is assumed to satisfy a second moment condition. Under these conditions, the mean semigroup of the superprocess automatically has a spectral gap. In the existing literature for almost sure convergence, either the spatial motion has to be a diffusion, or the branching mechanism has to obey restrictive conditions (in \cite{CRW} the coefficients of the branching mechanism have to satisfy
a Kato class condition,
while in \cite{KR, W}
the branching mechanisms are spatially independent, in \cite{CRSZ,CRW} the branching mechanisms have to satisfy a second moment condition, and in \cite{EKW} it has to satisfy a $p\,$th-moment condition
for some $p\in (1,2]$).

The present paper is devoted to establish weak and strong laws of large numbers for superprocesses where the underlying spatial motions are symmetric Hunt processes and the branching mechanisms are general. A key ingredient of
our approach is a
\textit{skeleton decomposition} for superprocesses, which provides a pathwise description of
the superprocess
in terms of immigration along a branching Markov process, called \textit{skeleton}.
The skeleton decomposition for superdiffusions was developed in
\cite{BKM, KPR}.
Very recently this decomposition was used in \cite{EKW} to study limit theorems for superdiffusions. In the present paper we
extend this decomposition to superprocesses whose spatial motion can be discontinuous,
and use it to make a connection between the asymptotic behavior of a branching Markov process and that of a superprocess. Our proof of
SLLN follows two main steps,
first to obtain SLLN along lattice times and then extend it to all times.
Our approach to SLLN along lattice times is different from that of \cite{EKW}.
Motivated by \cite{CRY}, where we established laws of large numbers for branching symmetric Hunt processes under a Kesten-Stigum $L\log L$ type condition, we use only $L^{1}$-convergence technique to establish the SLLN along lattice times, whereas the approach in \cite{EKW} involves a
$p$-th moment calculation for some $p\in (1,2]$.
As a result, the branching mechanism in the present paper has to obey
 an
 $L\log L$ type condition, while the branching mechanism in \cite{EKW} has to obey a
$p$-th moment condition with $p\in (1, 2]$.
So even within the superdiffusion case, our milder moment assumptions include superprocesses where the assumptions in \cite{EKW} may fail (see, for example, Example \ref{E:8.3} in Section \ref{sec8} below).
 For the transition from discrete to continuous time, \cite{EKW} adapted the idea of \cite{LRS}, where this transition is obtained through approximation of the indicator functions by resolvent functions. This idea works for superdiffusions
 because the spatial motion has continuous sample paths.
 But it
is not applicable for superprocesses when the underlying motion is a general symmetric Hunt process.
We overcome this difficulty by imposing a continuity condition on the $h$-transformed process.
Examples which satisfy these assumptions are given in Section \ref{sec8}.
In particular, the spatial motions in Example \ref{E:8.4} and Example \ref{E:8.5}
have discontinuous paths.

The remainder of this paper is organized as follows. We start Section \ref{sec2} with a review on definitions and basic properties of symmetric Hunt processes and superprocesses.
A key ingredient in constructing skeleton as well as spine decomposition of superprocesses is the Kuznetsov measures
or excursion measures (also known as $\mathbb{N}$-measures in literature) for the superprocesses.
 However, the existence of such a measure was taken for granted
in \cite{BKM, EKW, KPR}. In Subsection \ref{sup-Br} and Appendix \ref{sup-diffusion} of this paper,
we give sufficient conditions for the existence of such measures for general superprocesses and superdiffusions, respectively,
and thus also filling the missing pieces in \cite{BKM, EKW, KPR}. We then present our working hypothesis and main results
in Subsection \ref{sec2.3}.
We give a spine decomposition for the superprocess in Section \ref{sec3}, and then use it to prove
the $L^{1}$-convergence of an associated martingale in Section \ref{sec4}. In Section \ref{sec5} we give a detailed description of the skeleton space for superprocesses and show that the martingale limits for the superprocess and skeleton coincide. The proofs of weak and strong laws of large numbers
are given, respectively,
in Section \ref{sec6} and Section \ref{sec7}.
We present several examples to illustrate our results in Section \ref{sec8} and give a detailed verification of those examples in
Appendix \ref{verification} of this paper.

Throughout this paper, we use
``$:=$"
as a way of definition. For a function $f$,
$\|f\|_{\infty}:=\sup |f(x)|$.
For $a,b\in \mathbb{R}$, $a\wedge
b:=\min\{a,b\}$ and $a\vee b:=\max\{a,b\}$.
$\log^{+}x:=\log (x\vee 1)$
and
\begin{equation} \nonumber
\log^{*}x:=\left\{ \begin{aligned}
         x/e, &\quad x\le e, \\
  \log x,&\quad x> e.
                          \end{aligned} \right.
\end{equation}
For two positive functions $f$ and $g$, we use $f\stackrel{c}{\asymp}g$ to denote that there is a positive constant $c$ such that $c^{-1}f\le g\le c f$ on their common domain of definition. We also write $\asymp$ for
$\stackrel{c}{\asymp}$ if $c$ is unimportant.  We use $B(x,r)$ to denote the ball in $\R^{d}$ centered at $x$ with radius $r$.

\section{Preliminary}\label{sec2}

\subsection{Spatial process}\label{sec2.1}

Suppose $E$ is a locally compact separable metric space.
Let $E_{\partial} := E \cup \{ \partial \}$  be its one point compactification.
Denote by $\mathcal{B}(E)$ the Borel $\sigma$-field on $E$.
 The notation
$B\Subset E$ means that its closure $\bar B$ is compact in $E$.

 We use $\mathcal{B}_{b}(E)$ (respectively,
$\mathcal{B}^{+}(E)$) to denote
the space of bounded (respectively,
nonnegative) measurable functions on $(E,\mathcal{B}(E))$.
The space of continuous (and compactly supported) functions on $E$ will be denoted as $C(E)$ (and $C_{c}(E)$ resp.). Any functions $f$ on $E$ will be automatically extended to $E_{\partial}$ by setting $f(\partial)=0$.
 Suppose that $m$ is a $\sigma$-finite nonnegative Radon measure on $E$
 with full support. When $\mu$ is a measure on $\mathcal{B}(E)$ and $f$ is a measurable function, let $\langle f,\mu\rangle:=\int_{E}f(x)\mu(dx)$ whenever the right hand side makes sense. In particular, if $\mu$ has a density $\rho$ with respect to the measure
$m$, we write $\langle f,\rho\rangle$ for $\langle f,\mu\rangle$.
 If $g(t,x)$ is a measurable function on $[0,{\infty})\times E$, we say $g$ is locally bounded
 if $\sup_{t\in[0,T]}\sup_{x\in E}g(t,x)<{\infty}$ for every $T\in (0,{\infty})$.

 Let $\xi=(\Omega,\mathcal{H},\mathcal{H}_{t},\theta_{t},\xi_{t},\Pi_{x},\zeta)$ be an $m$-symmetric Hunt process on $E$.
 Here $\{\mathcal{H}_{t}:\ t\ge 0\}$ is the minimal admissible filtration, $\{\theta_{t}:\ t\ge
 0\}$ the time-shift operator of $\xi$ satisfying
 $\xi_{t}\circ\theta_{s}=\xi_{t+s}$ for $s,t\ge 0$, and $\zeta:=\inf\{t>0:\
 \xi_{t}=\partial\}$ the life time of $\xi$.
  Suppose for each $t>0$, $\xi$ has a
  transition density function   $p(t,x,y)$ with respect to the measure $m$,
	where  $p(t,x,y)$ is positive, continuous and symmetric in $(x,y)$.
  Let $\{P_{t}:\ t\ge 0\}$ be
 the Markovian transition semigroup of $\xi$, i.e.,
 $$P_{t}f(x):=\Pi_{x}\left[f(\xi_{t})\right]=\int_{E} p(t,x,y) f(y) m(dy)$$
 for any nonnegative measurable function $f$. The symmetric
 Dirichlet form on $L^{2}(E,m)$ generated by $\xi$ will be denoted as
 $(\mathcal{E},\mathcal{F})$:
 $$
 \mathcal{F}= \Big \{u\in L^{2}(E,m):\ \lim_{t\to 0}\frac{1}{t}\int_{E}\left(u(x)-P_{t}u(x)\right)u(x)m(dx)<{\infty} \Big\},
 $$
 $$\mathcal{E}(u,v)=\lim_{t\to 0}\frac{1}{t}\int_{E}\left(u(x)-P_{t}u(x)\right)v(x)m(dx),\quad u,v\in\mathcal{F}.$$
 It is known (cf. \cite{CF})
 that $(\mathcal{E},\mathcal{F})$ is quasi-regular and hence is quasi-homeomorphic to a regular Dirichlet form on a locally compact separable
 metric space.
 For any set $B\in\mathcal{B}(E)$, let $\tau_{B}:=\inf\{t\ge 0: \xi_{t}\not\in B\}$.

 \subsection{Superprocesses}\label{sec2.2}

 Suppose $\mf$ is the space of finite measures on $E$ equipped with the topology of weak convergence. The set of finite and compactly supported measures on $E$ is denoted by $\mc$.
 The main process of interest in this paper is an $\mf$-valued Markov process $X=\{X_{t}:t\ge 0\}$ with evolution depending on two quantities $P_{t}$ and $\psi_{\beta}$. Here $P_{t}$ is the semigroup of the spatial process $\xi$ and $\psi_{\beta}$ is the so-called
  \textit{branching mechanism},
 which takes the form
 \begin{equation}
 \psi_{\beta}(x,{\lambda})=-\beta(x){\lambda}+\alpha(x){\lambda}^{2}+\int_{(0,{\infty})}\left(e^{-{\lambda} y}-1+{\lambda} y\right)\pi(x,dy)\quad\mbox{ for }x\in E \mbox{ and }{\lambda}\ge 0,\label{2.1}
 \end{equation}
 with $\beta\in\mathcal{B}_{b}(E)$, $\alpha\in\mathcal{B}^{+}_{b}(E)$ and for each $x\in E$, $\pi(x,dy)$ being a measure concentrated on $(0,{\infty})$ such that $x\mapsto \int_{(0,{\infty})}(y\wedge y^{2})\pi(x,dy)$ is bounded from above.
 The distribution of $X$ is denoted by $\p_{\mu}$ if it is started in $\mu\in\mf$. $X$ is called a \textit{$(P_{t},\psi_{\beta})$-superprocess} (or \textit{$(\xi,\psi_{\beta})$-superprocess}) if $X$ is
 an
 $\mf$-valued process such that for all $\mu\in \mf$, $f\in\mathcal{B}^{+}_{b}(E)$ and $t\ge 0$,
 \begin{equation}
 \p_{\mu}\left[e^{-\langle f,X_{t}\rangle}\right]=e^{-\langle u_{f}(t, \cdot),\mu\rangle},\label{eq0}
 \end{equation}
 where $u_{f}(t,x):=-\log \p_{\delta_{x}}\left(e^{-\langle f,X_{t}\rangle}\right)$ is
 the unique nonnegative locally bounded solution
 to the following integral equation:
 \begin{equation}\label{eq1}
 u_{f}(t,x)=P_{t}f(x)-\int_{0}^{t}P_{s}\left(\psi_{\beta}(\cdot,u_{f}(t-s, \cdot))\right)(x)ds \quad\mbox{ for all }x\in E \mbox{ and }t\ge 0.
 \end{equation}
 The existence of such a process $X$ is established by \cite{D93}. Moreover, such a superprocess $X$ has a Hunt realization in $\mf$ (see, for example, \cite[Theorem 5.12]{Li}) such that $t\mapsto \langle f,X_{t}\rangle$ is almost surely right continuous for all bounded continuous functions $f$. We shall always work with this version.

 We define the Feynman-Kac semigroup $P^{\beta}_{t}$ by
 \begin{equation}
 P^{\beta}_{t}f(x):=\Pi_{x}\left[
 \exp\left(\int_{0}^{t}\beta(\xi_{s})ds\right)
 f(\xi_{t})\right]\quad\mbox{ for  }
 f\in\mathcal{B}^{+}_{b}(E).\nonumber
 \end{equation}
Then it follows by \cite[Lemma A.1.5]{D93} that the equation \eqref{eq1} is equivalent to the following integral equation:
 \begin{equation}\label{eq2}
 u_{f}(t,x)=P^{\beta}_{t}f(x)-\int_{0}^{t}P^{\beta}_{s}\left(\psi_{0}(\cdot,u_{f}(t-s, \cdot))\right)ds\quad\mbox{ for all }x\in E \mbox{ and }t\ge 0,
 \end{equation}
 where
 $$
 \psi_{0}(x,{\lambda}):=\alpha(x){\lambda}^{2}+\int_{(0,{\infty})}\left(e^{-{\lambda} y}-1+{\lambda} y\right)\pi(x,dy).
 $$
It is known (cf. \cite{D93}) that for every $\mu\in\mf$ and $f\in\mathcal{B}^{+}_{b}(E)$, the first moment of $\langle f,X_{t}\rangle$ exists and can be expressed as
 \begin{equation}\label{mean}
 \p_{\mu}\left(\langle f,X_{t}\rangle\right)=\langle P^{\beta}_{t}f,\mu\rangle.
 \end{equation}
Moreover,
the second moment of $\langle f,X_{t}\rangle$, if exists, can be expressed as
\begin{equation}\label{var}
\mbox{Var}_{\mu}\left(\langle f,X_{t}\rangle\right)=\int_{0}^{t}\langle P^{\beta}_{s}\left(\left(2\alpha+\int_{(0,{\infty})}y^{2}\pi(\cdot,dy)\right)\left(P^{\beta}_{t-s}f\right)^{2}\right),\mu\rangle ds.
\end{equation}
 If $\alpha(x)+\pi(x,(0,{\infty}))=0$ for $m$-a.e. $x\in E$,
then \eqref{var} implies that $\langle f,X_{t}\rangle=\langle P^{\beta}_{t}f,\mu\rangle$ for all $t\ge 0$
 $\p_{\mu}$-almost surely for all $\mu\in\mf$ and all bounded continuous functions $f$.
This trivial case is excluded in this paper, and we always assume that $m\left(\{x\in E:\ \alpha(x)+\pi(x,(0,{\infty}))>0\}\right)>0$.

Modelling
superprocess as a system of exit measures from  time-space open sets
has been
systematically developed in \cite{D93,D01,DK}. In particular, branching property and Markov properties of
such system
are
established there.
We take $\mathcal{B}^{+}_{b}([0,t]\times E)$ to be the space of nonnegative bounded measurable functions on $[0,t]\times E$. For every $\widetilde{f}\in\mathcal{B}^{+}_{b}([0,t]\times E)$, we extend $\widetilde{f}$ to $[0,t]\times \{\partial\}$ by setting $\widetilde{f}(s, \partial)=0$.
It follows from \cite[Theorem I.1.1 and Theorem I.1.2]{D93} that
for any open set $D\subset E$ and $t\ge 0$, there exists a random measure $\widetilde{X}^{D}_{t}$ on $[0,{\infty})\times E$ such that for every $\mu\in\mf$ and $\widetilde{f}\in\mathcal{B}^{+}_{b}([0,t]\times E)$,
 \begin{equation}
 \p_{\mu}\left[e^{-\langle\widetilde{f},\widetilde{X}^{D}_{t}\rangle}\right]=e^{-\langle\widetilde{u}^{D}_{\widetilde{f}}(t, \cdot),\mu\rangle},\label{2.4}
 \end{equation}
 where $\widetilde{u}^{D}_{\widetilde{f}}(t,x)$ is the unique nonnegative locally bounded solution to the following integral equation:
 \begin{equation}
 \widetilde{u}^{D}_{\widetilde{f}}(t, x)=\Pi_{x}\left[\widetilde{f}(t\wedge \tau_{D}, \xi_{t\wedge \tau_{D}})\right]-
 \Pi_{x}\left[\int_{0}^{t\wedge \tau_{D}}\psi_{\beta}\left(\xi_{s},\widetilde{u}^{D}_{\widetilde{f}}(t-s, \xi_{s})\right)ds\right]\quad\mbox{ for all }x\in E \mbox{ and }t\ge 0.\label{2.5}
 \end{equation}
Here $\tau_{D}$ is the fist time of $\xi$ leaving from $D$.
For an arbitrary function $f\in \mathcal{B}^{+}_{b}(D)$, let $\widetilde{f}(s,x):=f(x)$ if $x\in D$ and otherwise $\widetilde{f}(s, x)=0$. We can define a random measure $X^{D}_{t}$ on $D$ by setting $\langle f,X^{D}_{t}\rangle=\langle \widetilde{f},\widetilde{X}^{D}_{t}\rangle$ for all $f\in\mathcal{B}^{+}_{b}(D)$. This definition implies that $X^{D}_{t}$ is the
projection of $\widetilde{X}^{D}_{t}$ on $ D$.
For $f\in \mathcal{B}^+_b (D)$, we write $u^{D}_{f}(t,x)$ for $\widetilde{u}^{D}_{\widetilde{f}}(t,x)$.
It follows that $u^{D}_{f}(t,x)$ is the unique nonnegative locally bounded solution to the equation:
\begin{equation}
u^{D}_{f}(t,x)=\Pi_{x}\left[f(\xi_{t});t<\tau_{D}\right]-
 \Pi_{x}\left[\int_{0}^{t\wedge \tau_{D}}\psi_{\beta}\left(\xi_{s},u^{D}_{f}(t-s, \xi_{s})\right)ds\right]\quad\mbox{ for all }x\in D \mbox{ and }t\ge 0.\label{2.6}
\end{equation}
For an arbitrary $f\in\mathcal{B}^{+}(D)$, there exists a sequence of functions $f_{k}\in\mathcal{B}^{+}_{b}(D)$ such that $f_{k}\uparrow f$ pointwise. By \eqref{2.4}, $u^{D}_{f_{k}}(t,x)$ is
increasing in $k$ and we denote this limit by $u^{D}_{f}(t,x)\in [0,{\infty})$. With this notion, the
monotone convergence theorem implies that \eqref{2.4} is valid for all $f\in\mathcal{B}^{+}(D)$. The same argument shows that \eqref{eq2} holds for all $f\in\mathcal{B}^{+}(E)$.
If we define the subprocess $\xi^{D}$ of $\xi$ by
\begin{equation} \nonumber
\xi^{D}_{t}=\left\{ \begin{aligned}
         \xi_{t}, &\quad\mbox{ if } t<\tau_{D}, \\
  \partial,&\quad\mbox{ if } t\ge \tau_{D},
                          \end{aligned} \right.
                          \end{equation}
then \eqref{2.6} is equivalent to
\begin{equation}
u^{D}_{f}(t,x)=\Pi_{x}\left[f(\xi^{D}_{t})\right]-
 \Pi_{x}\left[\int_{0}^{t}\psi_{\beta}\left(\xi^{D}_{s},u^{D}_{f}(t-s, \xi^{D}_{s})\right)ds\right]\quad\mbox{ for all }x\in D \mbox{ and }t\ge 0.\nonumber
\end{equation}
As a process in time, $X^{D}:=\{X^{D}_{t}:t\ge 0\}$ is a $(\xi^{D},\psi_{\beta})$-superprocess. One may think of $X^{D}_{t}$ as describing the mass
of $X$ at time $t$ that historically avoids leaving $D$.
The following proposition follows directly from the properties of exit measures (cf. \cite[Theorem 1.2 and Theorem 1.3]{D93}). We omit its proof here.
\begin{proposition}[monotonicity in $D$]\label{prop1}
Suppose $D_{1}$ and $D_{2}$ are two open sets in $E$ with $D_{1}\subset D_{2}$. Then for every $\mu\in\mf$ and $t\ge 0$,
$$\p_{\mu}\left(X_{t}^{D_{1}}(B)\le X^{D_{2}}_{t}(B)\mbox{ for all }B\in\mathcal{B}(D_{1})\right) =1.
$$
\end{proposition}
Proposition \ref{prop1} implies that for any $f\in\mathcal{B}^{+}_{b}(D_{1})$, $u^{D_{1}}_{f}(t,x)\le u^{D_{2}}_{f}(t,x)$ for all $x\in D_{1}$ and $t\ge 0$. This monotonicity is also obtained in \cite[Section A.2]{EKW} through arguments on the integral equations \eqref{eq1} and \eqref{eq2}.
Although the underlying spatial motion is a diffusion in \cite{EKW}, their approach works generally.

\subsection{Kuznetsov measures for superprocesses}\label{sup-Br}

Suppose $X:=((X_{t})_{t\ge 0};\p_{\mu},\mu\in\mf)$ is a $(\xi,\psi_{\beta})$-superprocess, where the spatial motion $\xi$ is
an $m$-symmetric Borel right process
on a Luzin topological space $E$ with transition semigroup $(P_{t})_{t\ge 0}$ and the branching mechanism $\psi_{\beta}$ is given in \eqref{2.1}.
By \cite[Theorem 5.12]{Li}, the $(\xi,\psi_{\beta})$-superprocess has a right realization in $\mf$.
Let $\mathcal{W}^{+}_{0}$ denote the space of right continuous paths from $[0,{\infty})$ to $\mf$ having null
measure as a trap. Without loss of generality, we assume $X_{t}$ is the coordinate process in $\mathcal{W}^{+}_{0}$ and that $(\mathcal{F},(\mathcal{F}_{t})_{t\ge 0})$ is the minimal augmented
$\sigma$-fields on $\mathcal{W}^{+}_{0}$ generated by the coordinate process.

In this paper we will use two decompositions of superprocesses called spine decomposition and skeleton decomposition. To introduce these two decompositions we need to introduce Kuznetsov measures or excursion measures for
 superprocesses.

Let $\{Q_{t}(\mu,\cdot):=\p_{\mu}\left(X_{t}\in\cdot\right):t\ge 0,\ \mu\in\mf\}$ be the transition
semigroup of $X_{t}$ and $V_{t}$ be an operator on $\mathcal{B}^{+}_{b}(E)$ such that $V_{t}f(x):=u_{f}(t,x)$ for $t\ge0$ and $f\in\mathcal{B}^{+}_{b}(E)$.
Here $u_{f}(t,x)$ is the unique nonnegative locally bounded solution to \eqref{eq1}.
Then by \eqref{eq0}, we have
\begin{equation}
\int_{\mf}e^{-\langle f,\nu\rangle}Q_{t}(\mu,d\nu)=\exp\left(-\langle V_{t}f,\mu\rangle\right)\quad\mbox{ for }\mu\in\mf\mbox{ and }t\ge 0.
\nonumber
\end{equation}
It implies that $Q_{t}(\mu_{1}+\mu_{2},\cdot)=Q_{t}(\mu_{1},\cdot)*Q_{t}(\mu_{2},\cdot)$ for any $\mu_{1},\mu_{2}\in\mf$, and hence $Q_{t}(\mu,\cdot)$ is an infinitely divisible probability measure on $\mf$. By the semigroup property of $Q_{t}$, $V_{t}$ satisfies that
\begin{equation}
V_{s}V_{t}=V_{t+s}\quad\mbox{ for all }s,t\ge 0.\label{2.11}
\end{equation}
Moreover, by the infinite divisibility of $Q_{t}$, each operator $V_{t}$ has the representation
\begin{equation}
V_{t}f(x)=\int_{E}f(y){\Lambda}_{t}(x,dy)
+\int_{\mof} \left(1-e^{-\langle f,\nu\rangle}\right)L_{t}(x,d\nu)\label{2.12}
\end{equation}
for $t>0$ and $f\in\mathcal{B}^{+}_{b}(E)$
where ${\Lambda}_{t}(x,dy)$ is a bounded kernel on $E$, $\mof:=\mf\setminus\{0\}$ and $(1\wedge \langle 1,\nu\rangle)L_{t}(x,d\nu)$ is a bounded kernel from $E$ to $\mof$.
The operators $(V_{t})_{t\ge 0}$ satisfying \eqref{2.11} and \eqref{2.12} is called a \textit{cumulant semigroup}.

Let $Q_{t}^{0}$ be the restriction of $Q_{t}$ to $\mof$, and
$$
E_{0}:=\{x\in E:{\Lambda}_{t}(x,E)=0\mbox{ for all }t >0\}.
$$
Then  $x\in E_{0}$ if and only if
\begin{equation}\nonumber
V_{t}f(x)=\int_{\mof} \left(1-e^{-\langle f,\nu\rangle}\right)L_{t}(x,d\nu)\quad\mbox{ for all } t>0\mbox{ and }f\in\mathcal{B}^{+}_{b}(E).\label{4}
\end{equation}
  It follows by \cite[Proposition 2.8 and Theorem A.40]{Li} that for every $x\in E_{0}$, the family of measures $\{L_{t}(x,\cdot):t>0\}$ on $\mof$ constitutes an entrance law for the semigroup $\{Q_{t}^{0}:t\ge 0\}$ on $\mof$. Hence it
corresponds a unique $\sigma$-finite measure $\mathbb{N}_{x}$ on $(\mathcal{W}^{+}_{0},\mathcal{F})$ such that $\mathbb{N}_{x}(\{0\})=0$, and that
  for any $0<t_{1}<t_{2}<\cdots<t_{n}<{\infty}$,
  \begin{equation}
  \mathbb{N}_{x}\left(X_{t_{1}}\in d\nu_{1},X_{t_{2}}\in d\nu_{2},\cdots, X_{t_{n}}\in d\nu_{n}\right)
  =L_{t_{1}}(x,d\nu_{1})Q_{t_{2}-t_{1}}^{0}(\nu_{1},d\nu_{2})\cdots Q_{t_{n}-t_{n-1}}^{0}(\nu_{n-1},d\nu_{n}).\nonumber
  \end{equation}
It follows that
for all $t>0$ and $f\in\mathcal{B}^{+}_{b}(E)$,
\begin{equation}
\mathbb{N}_{x}\left(1-e^{-\langle f,X_{t}\rangle}\right)
=\int_{\mof} \left(1-e^{-\langle f,\nu\rangle}\right)L_{t}(x,d\nu)=u_{f}(t,x).\label{N}
\end{equation}
This measure $\mathbb{N}_{x}$ is called the \textit{Kuznetsov measure}
or {\it excursion measure}
(also known as $\mathbb{N}$-measure in \cite{DK})
of the $(P_{t},\psi_{\beta})$-superprocess corresponding to the entrance law $\{L_{t}(x,\cdot):t>0\}$,
and $E_0$ is
the collection of $x\in E$ such that
there exists the Kuznetsov measure $\mathbb{N}_{x}$ on $(\mathcal{W}^{+}_{0},\mathcal{F})$ corresponding to the $(P_{t},\psi_{\beta})$-superprocess.

For a constant ${\lambda}>0$, we use $V_{t}{\lambda}$ to denote $V_{t}f$ when the function $f(x)\equiv {\lambda}$. It follows by \eqref{2.12} that for every $x\in E$ and $t>0$,
\begin{equation}\nonumber
V_{t}{\lambda}(x)=
\lambda
{\Lambda}_{t}(x,E) +\int_{\mof}\left(1-e^{-{\lambda}\langle 1,\nu\rangle}\right)L_{t}(x,d\nu).
\end{equation}
The left hand side tends to $-\log\p_{\delta_{x}}\left(X_{t}=0\right)$ as ${\lambda}\to{\infty}$. Therefore, if $x\in E$ satisfies
\begin{equation}
\p_{\delta_{x}}\left(X_{t}=0\right)>0\quad\mbox{ for all }t>0,\label{condi0}
\end{equation}
then ${\Lambda}_{t}(x,E)=0$ for all $t>0$ and, consequently, $x\in E_{0}$.
Thus we have
\begin{equation}
\nonumber
\{x\in E: \p_{\delta_{x}}\left(X_{t}=0\right)>0\mbox{ for all }t>0\}\subset  E_{0}.
\end{equation}
In the study of spine decomposition, a key question is whether $\{x\in E:\alpha(x)>0\}\subset E_{0}$.
 Later in this subsection we will give sufficient conditions for this formula to be true.

In the study of  skeleton  decomposition
of general
superprocesses, we need the following condition:

\medskip

\noindent\textbf{Condition 1.}
There is a bounded positive function $w$ on $E$ such that
\begin{equation}
\p_{\delta_{x}}\left(e^{-\langle w,X_{t}\rangle}\right)=e^{-w(x)}\quad\mbox{ for all }x\in E\mbox{ and }t\ge 0.\label{a3}
\end{equation}

If the spatial motion $\xi$ is a
diffusion,
then this $(\xi,\psi_\beta)$-superprocess is called a superdiffusion. We refer to Appendix \ref{sup-diffusion} for an explicit definition. When $X$ is a superdiffusion, we
only need the following condition which is weaker than Condition 1.

\medskip

\noindent\textbf{Condition 1'.}
There is a locally bounded positive function $w$ on $E$
satisfying \eqref{a3}.
\medskip

Suppose Condition 1 holds.
Since $w$ is bounded on $E$,
equation \eqref{eq0} implies that
$u_{w}(t,x)=w(x)$, which is independent of $t\geq 0$.
Thus by \eqref{eq2} $w$ is the unique nonnegative bounded solution
to the following integral equation:
\begin{equation}
w(x)+\int_{0}^{t}P^{\beta}_{s}\left(\psi_{0}\left(\cdot,w(\cdot)\right)\right)(x)ds=P^{\beta}_{t}w(x)\quad\mbox{ for all }x\in E.\label{1.6}
\end{equation}
It follows from the Markov property of $X$ that for every $\mu\in\mf$ and $s,t\ge 0$,
$$\p_{\mu}\left(e^{-\langle w,X_{t+s}\rangle }\,|\,\mathcal{F}_{t}\right)=\p_{X_{t}}\left(e^{-\langle w,X_{s}\rangle }\right)=e^{-\langle u_{w}(s,\cdot),X_{t}\rangle}=e^{-\langle w,X_{t}\rangle}.$$
Hence
$\{e^{-\langle w,X_{t}\rangle}:\ t\ge 0\}$ is bounded positive $\p_{\mu}$-martingale with respect to the filtration $(\mathcal{F}_{t})_{t\ge 0}$.
 We now study the effects of change of measures via this martingale on the superprocess $X$.

\begin{proposition}\label{prop2.4}
Suppose Condition 1 is satisfied.
For every $\mu\in\mf$, define $\p^{*}_{\mu}$ by
$$d\p^{*}_{\mu}
:=e^{-\langle w,X_{t}\rangle+\langle w,\mu\rangle}d\p_{\mu}
\quad\mbox{ on }\mathcal{F}_{t}
\quad\mbox{ for }t\ge 0.$$
Then for every $\mu\in \mf$, $f\in\mathcal{B}^{+}_{b}(E)$ and $t\ge 0$,
\begin{equation}\label{Laplace-*}
\p^{*}_{\mu}\left(e^{-\langle f,X_{t}\rangle}\right)=e^{-\langle u^{*}_{f}(t,\cdot),\mu\rangle},
\end{equation}
 and $u^{*}_{f}(t,x)=u_{w+f}(t,x)-w(x)$ is the unique nonnegative locally bounded solution to the following integral equation
\begin{equation}
u^{*}_{f}(t,x)+\int_{0}^{t}P^{\beta^{*}}_{s}\left(\psi^{*}_{0}(\cdot,u^{*}_{f}(\cdot,t-s))\right)(x)ds=P^{\beta^{*}}_{t}f(x)
\quad\mbox{ for }x\in E\mbox{ and }t\ge 0,\label{1.8}
\end{equation}
where $P^{\beta^{*}}_{t}f(x):=\Pi_{x}\left[\exp\left(\int_{0}^{t}\beta^{*}(\xi_{s})ds\right) f(\xi_{t})\right]$ for $f\in\mathcal{B}^{+}_{b}(E)$,
    $$\beta^{*}(x):=\beta(x)-2\alpha(x)w(x)-\int_{(0,{\infty})}\left(1-e^{-w(x)y}\right)y\pi(x,dy),$$
    $$\psi^{*}_{0}(x,{\lambda}):=\alpha(x){\lambda}^{2}+\int_{(0,{\infty})}\left(e^{-{\lambda} y}-1+{\lambda} y\right)\pi^{*}(x,dy)$$
    and $\pi^{*}(x,dy):=e^{-w(x)y}\pi(x,dy)$.
\end{proposition}

\proof
By the definition of $\p^{*}$,
we have for every $\mu\in\mf$, $f\in \mathcal{B}^{+}_{b}(E)$ and $t\ge 0$,
\begin{equation}
\p^{*}_{\mu}\left(e^{-\langle f,X_{t}\rangle}\right)=e^{\langle w,\mu\rangle}\p_{\mu}\left(e^{-\langle w+f,X_{t}\rangle}\right)=e^{-\langle u_{w+f}(t,\cdot)-w(\cdot),\mu\rangle}=e^{-\langle u^{*}_{f}(t,\cdot),\mu\rangle}.\nonumber
\end{equation}
Recall that $u_{w+f}(t,x)$ satisfies the following equation
\begin{equation}
u_{w+f}(t,x)+ \int_{0}^{t}P^{\beta}_{s}\left(\psi_{0}\left(\cdot,u_{w+f}(\cdot,t-s)\right)\right)(x)ds= P^{\beta}_{t}(w+f)(x),\label{1.9}
\end{equation}
for $x\in E$ and $t\ge 0$.
Using \eqref{1.6} and \eqref{1.9} it is straightforward to check that $u^{*}_{f}(t,x)$ is a nonnegative locally bounded solution to the integral equation
\begin{equation}
u^{*}_{f}(t,x)+ \int_{0}^{t}P^{\beta}_{s}\left[\psi_{0}(\cdot,u^{*}_{f}(t-s,\cdot)+w(\cdot))-\psi_{0}(\cdot,w(\cdot))\right](x)ds
=P^{\beta}_{t}f(x),\label{1.10}
\end{equation}
for $x\in E$ and $t\ge 0$.
It then follows from \cite[Lemma A.1(ii)]{EKW} (by setting $A=B=E$, $T=t$, $g_{1}(x,s)=\beta^{*}(x)$, $g_{2}(x,s)=\beta(x)-\beta^{*}(x)$, $f_{1}(x)=f(x)$ and $f_{2}(x,s)=\psi_{0}(x,w(x))-\psi_{0}(x,u_{w+f}(s,x))$) that $u^{*}_{f}(t,x)$ satisfying \eqref{1.10} is also a solution to \eqref{1.8}. The uniqueness of the solution to \eqref{1.8} follows from \cite[Appendix A.2]{EKW} where the uniqueness of the solution for a more general class of integral equations is obtained.\qed

\medskip

\begin{remark}\label{suferdiffusion-*}\rm
If $X$ is a superdiffusion which satisfies Condition 1', we can also define a probability measure $\p^{*}_{\mu}$ such that the function $u^{*}_{f}(t,x)$ given by \eqref{Laplace-*} is the unique nonnegative locally bounded solution to \eqref{1.8}.
We defer its details to Appendix \ref{sup-diffusion}.
\end{remark}

 \medskip

Let $\psi^{*}_{\beta^{*}}(x,{\lambda}):=-\beta^{*}(x){\lambda}+\psi^{*}_{0}(x,{\lambda})$ for $x\in E$ and ${\lambda}\ge 0$. It is easy to check that $\psi^{*}_{\beta^{*}}(x,{\lambda})=\psi_{\beta}(x,{\lambda}+w(x))-\psi_{\beta}(x,w(x))$. Proposition \ref{prop2.4}
implies that under Condition 1,
$((X_{t})_{t\ge 0};\p^{*}_{\mu},\mu\in\mf)$
is a $(\xi,\psi^{*}_{\beta^{*}})$-superprocess.
We use $Q^{*}_{t}$ and $V^{*}_{t}$ to denote the transition semigroup and the cumulant semigroup of this $(\xi,\psi^{*}_{\beta^{*}})$-superprocess, respectively.
Denote by $E^{*}_{0}$ the collection of $x\in E$ such that
there exists the Kuznetsov measure $\mathbb{N}^{*}_{x}$ on $(\mathcal{W}^{+}_{0},\mathcal{F})$ corresponding to the $(\xi,\psi^{*}_{\beta^{*}})$-superprocess.
It follows that for every $x\in E^{*}_{0}$, every $t>0$ and $f\in\mathcal{B}^{+}_{b}(E)$,
\begin{equation}
\mathbb{N}^{*}_{x}\left(1-e^{-\langle f,X_{t}\rangle}\right)
=u^{*}_{f}(t,x).\nonumber
\end{equation}
Applying similar argument as in the beginning of this subsection, we have
\begin{equation}
\nonumber
\{x\in E:\p^{*}_{\delta_{x}}\left(X_{t}=0\right)>0\mbox{ for all }t>0\} \subset E^{*}_{0}.
\end{equation}
In the study of skeleton decomposition, a key question is whether $\{x\in E:\alpha(x)>0\}\subset E^*_{0}$.
Sufficient conditions for this will be given at the end of this subsection.

 The following assumption will be used later to establish the spine and skeleton decompositions.

 \medskip
\noindent\textbf{Assumption 1.}
 Suppose either one of the following conditions holds.
\begin{description}
\item{\rm (i)} $X$ is a $(\xi,\psi_{\beta})$-superprocess satisfying Condition 1 and
 $$
E_{+}:=\{x\in E:\ \alpha(x)>0\}\subset E_0\cap E^*_0.
$$

\item{\rm (ii)}  $X$ is a $(\xi,\psi_{\beta})$-superdiffusion satisfying Condition 1'.

\end{description}

 \medskip

\begin{remark}\label{R:2.4}  \rm
We note that $E_{+}\subset E_{0}$ is a necessary condition for the spine decomposition
while $E_{+}\subset E^{*}_{0}$ is necessary for the skeleton decomposition
for the  $(\xi,\psi_{\beta})$-superprocess $X$.
If $X$ is a superdiffusion that satisfies Condition 1',
then $E_{+}\subset E_0\cap E^*_0$ always holds.
See Appendix \ref{sup-diffusion} for its proof.
\end{remark}

\medskip

For general superprocesses, we now give some sufficient conditions for Assumption 1(i).
Consider a special superprocess denoted by $\widetilde{X}$, where
the spatial motion is a conservative Borel right process on $E$ and the branching mechanism
$\psi$ is given by
\begin{equation}
\psi({\lambda}):= -b{\lambda}+a{\lambda}^{2}+\int_{(0,{\infty})}\left(e^{-{\lambda} y}-1+{\lambda} y\right)\eta(dy)\quad\mbox{ for }x\in E,\ {\lambda}\ge 0,\label{psi2}
\end{equation}
where $b\in\mathbb{R}$, $a\in\mathbb{R}^{+}$ and $\eta$ is a measure supported in $(0,{\infty})$ satisfying $\int_{(0,{\infty})}(x\wedge x^{2})\eta(dx)<{\infty}$.
The
total mass process $\{Y_{t}:=\langle 1,\widetilde{X}_{t}\rangle;\ t\ge 0\}$ is a one-dimensional \textit{continuous-state branching process} (CB process
in abbreviation) with branching mechanism $\psi$.
It is well-known that for the CB process $Y$,
\begin{equation}
\p_{\mu}\left[ e^{-{\lambda} Y_{t}}\right]=e^{-v_{t}({\lambda})\mu(E)}\quad\mbox{ for all }\mu\in\mf,\ {\lambda}\ge 0\mbox{ and }t\ge 0,\nonumber
\end{equation}
where $v_{t}({\lambda})$ is the unique nonegative solution of
$$v_{t}({\lambda})={\lambda}-\int_{0}^{t}\psi(v_{s}({\lambda}))ds\quad\mbox{ for all }t,{\lambda}\ge 0.$$
The family of maps $\{{\lambda}\mapsto v_{t}({\lambda}); t\ge 0\}$ is called the cumulant semigroup of the CB process. It follows by \cite[Theorem 3.7]{Li} that $\p_{\delta_{x}}\left(Y_{t}=0\right)>0$ for every $t>0$ and every $x\in E$ if and only if $\psi({\infty})={\infty}$ and $\int_{N}^{{\infty}}\psi({\lambda})^{-1}d{\lambda}<{\infty}$ for some $N>0$.
For the $(\xi,\psi_{\beta})$-superprocess,
the following result gives a sufficient condition under which \eqref{condi0} holds for every $x\in E$.

\begin{proposition}\label{prop2.1}
 Suppose that there is a spatially independent branching mechanism $\psi$ in the form of \eqref{psi2} such that $\psi({\infty})={\infty}$, $\int^{{\infty}}\psi({\lambda})^{-1}d{\lambda}<{\infty}$, and
\begin{equation}
\psi_{\beta}(x,{\lambda})\ge \psi({\lambda})\quad\mbox{ for all }x\in E\mbox{ and } {\lambda}\ge 0.\label{condi1}
\end{equation}
  Then $\p_{\delta_{x}}\left(X_{t}=0\right)>0$ for every $x\in E$ and every $t>0$, and hence $E_{0}=E$.
\end{proposition}
\proof
Let $V_{t}$ and $\hat{V}_{t}$ denote the cumulant semigroups of a $(\xi,\psi_{\beta})$- and a $(\xi,\psi)$-superprocesses, respectively.
Given \eqref{condi1},
we have by \cite[Corollary 5.18]{Li} (cf.\cite[Lemma 4.5]{EKW}) that
for every ${\lambda}\ge 0$ and $t\ge 0$, $V_{t}{\lambda}(x)\le \hat{V}_{t}{\lambda}(x)$ pointwise.
Let $v_{t}({\lambda})$ be the cumulant semigroup of the CB-process with branching mechanism $\psi$. Then for every ${\lambda}\ge 0$ and $t\ge 0$, we have $\hat{V}_{t}{\lambda}(x)\le v_{t}({\lambda})$ for all $x\in E$, with the equality holds if $\xi$ is a conservative process on $E$. Under the conditions on $\psi$, we know by \cite[Theorem 3.7]{Li} that
$\lim_{{\lambda}\to{\infty}}v_{t}({\lambda})=-\log\p_{\delta_{x}}\left(Y_{t}=0\right)<{\infty}$ for all $t>0$. Thus $-\log\p_{\delta_{x}}\left(X_{t}=0\right)=\lim_{{\lambda}\to{\infty}}V_{t}{\lambda}(x)<{\infty}$ for every $x\in E$ and $t>0$, and hence we complete the proof.
\qed

\begin{proposition}\label{prop2.5}
Under Condition 1 and the condition of Proposition \ref{prop2.1}, we have $E^{*}_{0}=E$.
 Consequently, Assumption 1(i) holds.
\end{proposition}

\proof
Since the event $\{X_{t}=0\}\in\mathcal{F}_{t}$, it follows that
$$
\p^{*}_{\mu}\left(X_{t}=0\right)=e^{\langle w,\mu\rangle}\p_{\mu}\left(e^{-\langle w,X_{t}\rangle};X_{t}=0\right)=e^{\langle w,\mu\rangle }\p_{\mu}(X_{t}=0)
$$
 for every $\mu\in\mf$. Thus
$\{x\in E:\p_{\delta_{x}}\left(X_{t}=0\right)>0\mbox{ for all }t>0\}= \{x\in E:\p^{*}_{\delta_{x}}\left(X_{t}=0\right)>0\mbox{ for all }t>0\}\subset E^{*}_{0}$. Note that by Proposition \ref{prop2.1}, $\{x\in E:\p_{\delta_{x}}\left(X_{t}=0\right)>0\mbox{ for all }t>0\}=E$.
Hence we get $E^{*}_{0}=E$.
\qed

\begin{remark}\label{rm1}
\rm
For the CB process $Y$ with branching mechanism $\psi$, it is known that if $\psi({\infty})={\infty}$ and $\int^{{\infty}}\psi({\lambda})^{-1}d{\lambda}<{\infty}$, then for every $x\in E$, $\p_{\delta_{x}}\left(Y_{t}=0\mbox{ for some }t\ge 0\right)=e^{-z_{\psi}}$ where $z_{\psi}:=\sup\{{\lambda}\ge 0:\psi({\lambda})\le 0\}\in [0,{\infty})$.
For the $(\xi,\psi_{\beta})$-superprocess $X$, let
$w(x):=-\log\p_{\delta_{x}}\left(\mathcal{E}\right)$, where the event $\mathcal{E}:=\{X_{t}=0\mbox{ for some }t\ge 0\}$.
Then under Condition \eqref{condi1},
we have $$w(x)=\lim_{t\to{\infty}}\lim_{{\lambda}\to{\infty}}V_{t}{\lambda}\le \lim_{t\to{\infty}}\lim_{{\lambda}\to{\infty}}v_{t}({\lambda})=z_{\psi}<{\infty},\quad\forall x\in E.$$  Moreover, by the bounded convergence theorem, $w$ satisfies \eqref{a3}.
In particular, if the  bounded function $w$ in Condition 1
is given by $w(x)=-\log\p_{\delta_{x}}(\mathcal{E})$, then the process $((X_{t})_{t\ge 0};\p^{*}_{\mu},\mu\in\mf)$
can be obtained from $X$ by conditioning on $\mathcal{E}$, i.e.,
$\p^{*}_{\mu}(X_{t}\in\cdot)
=\p_{\mu}(X_{t}\in\cdot\,|\,\mathcal{E})$, cf. \cite{KPR}.
\end{remark}

 \subsection{Assumptions and main results}\label{sec2.3}

 Since $\xi$ has a transition density function $p(t,x,y)$ with respect to the measure $m$,
 it follows
 that for each $t>0$, $P^{\beta}_{t}$ admits
 an integral kernel
 with respect to $m$. We denote this kernel by
 $p^{\beta }(t,x,y)$. It is positive, symmetric and continuous in $(x,y)$ for each $t>0$ and it satisfies that
 \begin{equation}\label{density-FK}
 e^{-t\|\beta\|_{\infty}}p(t,x,y)\le p^{\beta}(t,x,y)\le  e^{t\|\beta\|_{\infty}}p(t,x,y)\quad\mbox{ for all }(t,x,y)\in(0,{\infty})\times E\times E.
 \end{equation}
  This semigroup $P^{\beta }_{t}$ associates with
 a quadratic form $(\mathcal{E}^{(\beta ) },\mathcal{F})$,  where
$$
\mathcal{E}^{(\beta ) }(u,u)  :=  \mathcal{E}(u,u)-\int_{E}u(x)^{2}\beta(x)m(dx),
 \quad   u\in \mathcal{F}.
$$
Since $\beta$ is a bounded function, $\mathcal{E}^{(\beta ) }(u,u)\ge -\|\beta\|_{\infty}\int_{E}u(x)^{2}m(dx)$ for all $u\in\mathcal{F}$.
Thus by \cite{ABM}, $\{P^{\beta}_{t}:t\ge 0\}$ is a strongly continuous semigroup on $L^{2}(E,m)$.
We define
\begin{equation}\label{l1}
\lambda_{1}:=\inf\left\{\mathcal{E}^{(\beta ) }(u,u):u\in
\mathcal{F}\hbox{ with }  \int_{E}u(x)^{2}m(dx)=1\right\}.
\end{equation}
Obviously $\lambda_{1}\ge
-\|\beta^{+}\|_{\infty}
$ by definition.

\noindent\textbf{Assumption 2.}
$\lambda_{1}<0$ and
there is a positive continuous function $h \in \mathcal{F}$ with
$\int_{E}h(x)^{2}m(dx)=1$
so that
$ \mathcal{E}^{(\beta ) }(h, h) = \lambda_1$.

\medskip

Observe that if $u$ is a minimizer for \eqref{l1}, then so is $| u|$.
Hence one can always take a non-negative minimizer. Assumption 2
requires that there is a minimizer  for \eqref{l1} that can be chosen
to be positive everywhere. Clearly the following property holds for $h$:
\begin{equation}
\mathcal{E}^\beta (h,v)=\lambda_{1}\langle  h,v\rangle
\quad \hbox{for every }  v\in\mathcal{F}.\label{1.3}
\end{equation}
Let
$\sigma(\mathcal{E}^{(\beta ) })$ denote the spectrum of the self-adjoint operator associated with $\mathcal{E}^{(\beta ) }$.
Then $\lambda_{1}$ is a simple eigenvalue in $\sigma(\mathcal{E}^{(\beta ) })$ with eigenfunction $h$. It holds that $h=e^{\lambda_{1}t}P^{\beta}_{t}h$ on $E$.

\begin{remark}\rm
From the definition of $\mathcal{E}^{(\beta ) }$,
Assumption 2
implies that $m\left(\{x\in E:\ \beta(x)>0\}\right)>0$.
\end{remark}

Let $\lambda_{2}$ be the second bottom of $\sigma(\mathcal{E}^{(\beta ) })$, that is,
$$
\lambda_{2}:=\inf \Big\{ \mathcal{E}^{(\beta ) }(u,u):\ u\in\mathcal{F}, \ \int_{E}u(x)h(x)m(dx)=0,\ \int_{E}u(x)^{2}m(dx)=1
\Big\}.
$$

\medskip

\noindent\textbf{Assumption 3.}
There is a positive spectral gap in
$\sigma(\mathcal{E}^{(\beta ) })$:
$\lambda_{h}:=\lambda_{2}-\lambda_{1} >0$.
\medskip

Define the $h$-transformed semigroup $\{P^h_t; t\geq 0\}$ from
$\{ P^{\beta }_{t}; t\geq 0\}$ by
\begin{equation}\label{e:1.8}
P^{h}_{t}f(x)=
\frac{e^{\lambda_{1}t}}{h(x)}P^{\beta }_{t}(hf)(x)
\quad\mbox{ for }x\in E\mbox{ and }f\in\mathcal{B}^{+}_{b}(E).
\end{equation}
Then it is easy to see that $\{P^h_t: t\geq 0\}$ is an $\widetilde m$-symmetric semigroup, where $\widetilde m:=h^2 m$, and
$1$ is an eigenfunction of $P^{h}_{t}$
with eigenvalue $1$. Furthermore the spectrum of the infinitesimal generator of
$\{P^h_t: t\geq 0\}$ in $L^2(E; \widetilde m)$ is the spectrum of the infinitesimal generator of $\{ P^{\beta }_{t}: t\geq 0\}$ in $L^2(E; m)$ shifted by
$\lambda_1$.
Hence under Assumption 3,
we have the following Poincar\'{e} inequality:
\begin{equation} \label{p1}
\|P^{h}_{t}\varphi\|_{L^{2}(E,\widetilde{m})}\le
e^{-\lambda_h t}\|\varphi\|_{L^{2}(E,\widetilde{m})}
\end{equation}
for all $\varphi\in L^{2}(E,\widetilde{m})$ with
$\int_{E}\varphi(x)\widetilde{m}(dx)=0$.

We use $\mh$ to denote the space of all finite measures $\mu$ on $E$ with $\langle h,\mu\rangle<{\infty}$. By \eqref{mean} and the fact that $e^{\lambda_{1}t}P^{\beta}_{t}h=h$,
we can verify that under Assumption 2,
$$W^{h}_{t}(X):=e^{\lambda_{1}t}\langle h,X_{t}\rangle\quad\mbox{ for }t\ge 0$$
is a nonnegative $\p_{\mu}$-martingale for all $\mu\in\mh$. We define $W^{h}_{\infty}(X):=\lim_{t\to{\infty}}W^{h}_{t}(X)$.

\begin{theorem}\label{them2}
Suppose Assumptions 1-3 hold. If
\begin{equation}\label{llogl1}
\langle \alpha\log^{+}h,h^{2}\rangle+\langle \int_{(0,{\infty})}r\log^{*}(rh(\cdot))\pi(\cdot,dr),h^{2}\rangle<{\infty}  ,
\end{equation}
then the non-negative
martingale $W^{h}_{t}(X)$ converges to $W^h_\infty (X)$
as $t\to \infty$
$\p_{\mu}$-a.s. and in $L^{1}(\p_{\mu})$ for every $\mu\in\mh$.
In particular,
 $W^{h}_{\infty}(X)$ is non-degenerate in the sense that
$\p_{\mu}\left(W^{h}_{\infty}(X)>0\right)>0$ for any nontrivial $\mu\in \mh$.
\end{theorem}

\begin{remark}\label{rm for them2}\rm
We notice that only the spine decomposition of superprocesses is applied in the proof of this theorem. Therefore the above result
 in fact holds with the same proof
 under a weaker condition, say, Assumptions 2-3 and that $E_{+}\subset E_{0}$.
\end{remark}

As we mentioned earlier,
Assumption 1 is a fundamental assumption for the existence of the skeleton space.
The skeleton space for supercritical superprocesses offers a pathwise description of the superprocess in terms of an $\mf$-valued Markov process $X^{*}$ and a supercritical branching Markov process $Z$ dressed with an immigration process $I$ on a rich  probability space with probability measures $\pp_{\mu}$. A more detailed description and discussion of skeleton space is deferred to Section \ref{sec5}. In \cite{KPR}, the skeleton space was constructed for superdiffusions using the function $w(x)=-\log\p_{\delta_{x}}\left(\exists t\ge0: \langle 1,X_{t}\rangle=0\right)$. The key property of $w$ used in the skeleton construction is that $w$ gives rise to the multiplicative $\p_{\mu}$-martingale $\{e^{-\langle w,X_{t}\rangle}:t\ge 0\}$.
Later  the skeleton space is established in \cite{EKW} for superdiffusions assuming only
the existence of such a martingale function $w$. The main reason why skeleton space is constructed only for superdiffusions in \cite{KPR} and \cite{EKW}
is due to the proof of \cite[Lemma 6.1]{KPR}, where a comparison principle for elliptic differential operators is used. This comparison principle allows one to conclude that if $B_{1}$ and $B_{2}$ are domains with $\mathrm{supp}(f)\subseteq B_{1}\subseteq B_{2}$, then $u^{B_{1}}_{f}(t,x)\le u^{B_{2}}_{f}(t,x)$ pointwise.
 Proposition \ref{prop1} implies that this monotonicity in fact holds more generally as
 in the set-up  of  this paper, where the underlying spatial motion can be discontinuous.
 Therefore we  can establish  the existence of skeleton space for superprocesses defined
 in Section \ref{sec5}.

In Theorem \ref{them1} below we show that under
Assumptions 1-3,
$$
W^{h/w}_{t}(Z):=e^{\lambda_{1}t}\langle\, \frac{h}{w}\, ,Z_{t}\rangle\quad\mbox{ for }t\ge 0
$$
equals $W^h_t(X)$ and so
is a nonnegative $\pp_{\mu}$-martingale for all $\mu\in\mh$. Let
$ W^{h/w}_{\infty}(Z)
:=\lim_{t\to{\infty}}W^{h/w}_{t}(Z)$.

\begin{theorem}\label{them1}
Suppose Assumptions 1-3 hold. For every $\mu\in\mh$ and $t\ge 0$, $W^{h}_{t}(X)=W^{h/w}_{t}(Z)$ $\mathbb{P}_{\mu}$-a.s.
and in particular
$W^{h}_{\infty}(X)=W^{h/w}_{\infty}(Z)$ $\mathbb{P}_{\mu}$-a.s.
\end{theorem}

The proof of Theorem \ref{them1} will be given in Section \ref{sec5}.
We will establish laws
of large numbers under the following moment conditions.

\noindent\textbf{Assumption 4.}
\begin{description}
\item{\rm (i)} {\rm ($L \log L$ condition)}
 $\langle \log^{+}h,h^{2}\rangle+\langle \int_{(0,{\infty})}r\log^{*}(rh(\cdot))\pi(\cdot,dr),h^{2}\rangle<{\infty}.$

\item{\rm (ii)} $\langle (\int_{(0,{\infty})}r^{2}e^{-w(\cdot)r}\pi(\cdot,dr))^{2},1\wedge h^{4}\rangle<{\infty}.$
\end{description}

\begin{theorem}[Weak law of large numbers]\label{them3}
Suppose Assumptions 1-3 and
Assumption 4(i)
hold. Then for all $\mu\in\mh$ and all $f\in\mathcal{B}^{+}(E)$ with $f/h$ bounded,
\begin{equation}\nonumber
\lim_{t\to{\infty}}e^{\lambda_{1}t}\langle f,X_{t}\rangle=\langle f,h\rangle W^{h}_{\infty}(X)\quad\mbox{ in }L^{1}(\p_{\mu}).
\end{equation}
\end{theorem}

The proof of Theorem \ref{them3} will be given in Section \ref{sec6}.
The next assumption assumes that the strong law of large numbers
holds for the supercritical branching Markov process $Z$
along an increasing  sequence of  discrete times.

\noindent\textbf{Assumption 5.} For all $\mu\in\mc$, $\sigma>0$ and $\phi\in\mathcal{B}^{+}_{b}(E)$,
$$\lim_{n\to{\infty}}e^{\lambda_{1}n\sigma}\langle \frac{h}{w}\phi,Z_{n\sigma}\rangle =\langle \phi,h^{2}\rangle W^{h/w}_{\infty}(Z)\quad\pp_{\mu}\mbox{-a.s.}$$

We will prove at the end of Section \ref{sec6} that under Assumptions 1-3 and Assumption 4(i), Assumption 5 is equivalent to a weaker condition as follows.

\noindent\textbf{Assumption 5'.} For all $\mu\in\mc$, $\sigma>0$ and $\phi\in\mathcal{B}^{+}_{b}(E)$,
there exists $m\in \mathbb{N}$ such that
$$\lim_{n\to{\infty}}e^{\lambda_{1}n\sigma}\langle \frac{h}{w}P^{h}_{m\sigma}\phi,Z_{n\sigma}\rangle =\langle \phi,h^{2}\rangle W^{h/w}_{\infty}(Z)\quad\pp_{\mu}\mbox{-a.s.}
$$

\medskip

\noindent\textbf{Assumption 6.} The semigroup $\{P^{h}_{t}:t\ge 0\}$ satisfies that
$$\lim_{t\to 0}\|P^{h}_{t}f-f\|_{\infty}=0\quad\mbox{ for all }f\in C_{0}(E).$$
Here $C_{0}(E):=\{f\in C(E)\mbox{ with }\lim_{x\to\partial}f(x)=0\}$.

The following is the main result of this paper, which extends the main results
of \cite{EKW,LRS} to superprocesses having possibly discontinuous spatial motions.
Its proof will be given in Section \ref{sec7}.

\begin{theorem}[Strong law of large numbers]\label{them4}
Suppose Assumptions 1-6 hold. Then there exists $\Omega_{0}\subset \Omega$ of $\p_{\mu}$-full probability for every $\mu\in\mh$ such that on $\Omega_{0}$, for every $m$-almost everywhere continuous nonnegative measurable function $f$ with $f/h$ bounded, we have
\begin{equation}
\lim_{t\to{\infty}}e^{\lambda_{1}t}\langle f,X_{t}\rangle=\langle f,h\rangle W^{h}_{\infty}(X).\label{2.23}
\end{equation}
The convergence in \eqref{2.23} also holds in $L^{1}(\p_{\mu})$.
\end{theorem}

\medskip

\section{Change of measures}\label{sec3}

Recall that $h\in \FF$ is the minimizer in Assumption 2.
Since $h\in \mathcal{F}$, by Fukushima's decomposition, we have for
q.e. $x\in E$, $\Pi_{x}$-a.s.
$$
h(\xi_{t})-h(\xi_{0})=M^{h}_{t}+N^{h}_{t}\quad \mbox{ for }   t\ge 0,
$$
where $M^{h}_{t}$ is a martingale additive functional of $\xi$ having
finite energy and $N^{h}_{t}$ is a continuous additive functional of
$\xi$ having zero energy. It follows from \eqref{1.3} and \cite[Theorem 5.4.2]{FOT} that $N^{h}_{t}$ is of
bounded variation, and
$$N^{h}_{t}=-\lambda_{1}\int_{0}^{t}h(\xi_{s})ds-\int_{0}^{t}h(\xi_{s})\beta(\xi_{s})ds\quad\mbox{ for } t\ge 0.$$
Following \cite[Section 2]{CFTYZ}
(see also \cite[Section 2]{CRY}),
we define
a local martingale on the random time interval
$[0,\zeta_{p})$ by
\begin{equation} \label{(1)}
M_{t}:=\int_{0}^{t}\frac{1}{h(\xi_{s-})}dM^{h}_{s},
\quad  t\in [0,\zeta_{p}),
\end{equation}
where $\zeta_{p}$ is the predictable part of the life time $\zeta$
of $\xi$. Then the solution $R_{t}$ of the stochastic differential
equation
\begin{equation}
R_{t}=1+\int_{0}^{t}R_{s-}dM_{s},  \quad
 t\in [0,\zeta_{p}),\label{(2)}
\end{equation}
is a positive local martingale on $[0,\zeta_{p})$ and hence a
supermartingale. As a result, the formula
$$d\Pi^{h}_{x}=R_{t}d\Pi_{x}\quad \mbox{on }\mathcal{H}_{t}\cap\{t<\zeta\}\quad\mbox{ for }x\in E$$
uniquely determines a family of subprobability measures
$\{\Pi^{h}_{x}:x\in E\}$ on $(\Omega,\mathcal{H})$. We denote $\xi$
under $\{\Pi^{h}_{x}:x\in E\}$ by $\xi^{h}$, that is
$$\Pi^{h}_{x}\left[f(\xi^{h}_{t})\right]=\Pi_{x}\left[R_{t}f(\xi_{t}):t<\zeta\right]
\quad\mbox{ for all }t\ge 0\mbox{ and }f\in\mathcal{B}^{+}_{b}(E).$$
It follows from
\cite[Theorem 2.6]{CFTYZ} that the process $\xi^{h}$ is an irreducible recurrent
$\widetilde m$-symmetric right Markov process, where $\widetilde{m}(dy)=h(y)^{2}m(dy)$.
 Note that by \eqref{(1)},
\eqref{(2)} and Dol\'{e}an-Dade's formula,
\begin{equation}\label{(3)}
R_{t}=\exp\big( M_{t}-\frac{1}{2}\langle M^{c}\rangle_{t}\big) \prod_{0<s\le
t}\frac{h(\xi_{s})}{h(\xi_{s-})}\exp\left( 1-\frac{h(\xi_{s})}{h(\xi_{s-})}\right),
\quad t\in [0, \zeta_p),
\end{equation}
where $M^{c}$ is the continuous martingale part of $M$. Applying
Ito's formula to $\log h(\xi_{t})$, we obtain that for q.e. $x\in E$,
$\Pi_{x}$-a.s. on  $[0,\zeta)$,
\begin{equation}
\log h(\xi_{t})-\log h(\xi_{0})=M_{t}-\frac{1}{2}\langle
M^{c}\rangle_{t}+\sum_{s\le t}
\left(\log\frac{h(\xi_{s})}{h(\xi_{s-})}-\frac{h(\xi_{s})-h(\xi_{s-})}{h(\xi_{s-})}\right)-\lambda_{1}t-\int_{0}^{t}\beta(\xi_{s})ds.\label{(4)}
\end{equation}
By \eqref{(3)} and \eqref{(4)}, we get
$$
R_{t}=\exp\left( \lambda_{1}t+\int_{0}^{t}\beta(\xi_{s})ds\right)
\frac{h(\xi_{t})}{h(\xi_{0})}\quad\mbox{on }[0,\zeta).
$$
Therefore for any $f\in\mathcal{B}^{+}_{b}(E)$,
\begin{eqnarray}
\Pi^{h}_{x}\left(f(\xi^{h}_{t})\right)
=\frac{e^{\lambda_{1}t}}{h(x)}\Pi_{x}\left(e^{\int_{0}^{t}\beta(\xi_{s})ds}h(\xi_{t})f(\xi_{t})\right)
 = \frac{e^{\lambda_{1}t}}{h(x)}P^{\beta }_{t}(hf)(x)
 =P^{h}_{t}f(x) .
 \nonumber
\end{eqnarray}
This implies that
the transition semigroup of $\xi^h$ is exactly the semigroup
$\{P^{h}_{t}: t\geq 0\} $
obtained from $P^{\beta }_{t}$ through Doob's $h$-transform. Let $(\mathcal{E}^{h},\mathcal{F}^{h})$ be the symmetric Dirichlet
form in $L^{2}(E; \widetilde{m})$ generated by $\xi^{h}$.
Then $f\in \mathcal{F}^{h}$ if and only if $fh\in\mathcal{F} $, and
$$
\mathcal{E}^{h}(f,f)=\mathcal{E}^{(\beta ) }(fh,fh)-\lambda_{1}\int_{E}f(x)^{2}h(x)^{2}m(dx).
$$
In other words, $\Phi^{h}:f\mapsto fh$ is an isometry from $(\mathcal{E}^{h},\mathcal{F}^{h})$ onto $(\mathcal{E}^{\beta +\lambda_{1}m},\mathcal{F} )$ and from $L^{2}(E,\widetilde{m})$ onto $L^{2}(E,m)$.
Let
$\sigma(\mathcal{E}^{h})$ denote the spectrum of
$- \widetilde {\mathcal L}$, where $\widetilde {\mathcal L}$ is the self-adjoint operator associated with
the Dirichlet form $(\mathcal{E}^{h}, \mathcal{F}^h)$ in $L^{2}(E; \widetilde{m})$.
We know from  \cite[Theorem 2.6]{CFTYZ} that the constant function $1$ belongs to
$\mathcal{F}^{h}$, and $\mathcal{E}^{h}(1,1)=0$.
Hence $0\in \sigma(\mathcal{E}^{h})$ is a simple
eigenvalue and $1$ is the corresponding eigenfunction.
In particular,
$$
\lambda^{h}_{1}:=\inf\left\{\mathcal{E}^{h}(u,u):\ u\in\mathcal{F}^{h} \hbox{ with }  \int_{E}u(x)^{2}\widetilde{m}(dx)=1\right\}=0.
$$
Let $\lambda^{h}_{2}$ be the second bottom of $\sigma(\mathcal{E}^{h})$, i.e.
$$
\lambda^{h}_{2}:=\inf\left\{\mathcal{E}^{h}(u,u):\ u\in\mathcal{F}^{h} \hbox{ with }  \int_{E}u(x)\widetilde{m}(dx)=0
\hbox{ and }  \int_{E}u(x)^{2}\widetilde{m}(dx)=1\right\}.
$$
In view of the isometry $\Phi^h$,  we have $\lambda^{h}_{2}= \lambda_2 -\lambda_1$.
So Assumption 3
is equivalent to assuming $\lambda^h_2>0$.
The $h$-transformed process $\xi^{h}$ has a transition density function $p^{h}(t,x,y)$ with respect to $\widetilde{m}$, which
is positive, symmetric and continuous in $(x,y)$ for each $t>0$ and
is
related to $p^{\beta}(t,x,y)$ by
\begin{equation}\label{density-trans}
p^{h}(t,x,y)=e^{\lambda_{1}t}\frac{p^{\beta}(t,x,y)}{h(x)h(y)}\quad\mbox{ for }x,y\in E\mbox{ and }t\ge 0.
\end{equation}
Define
\begin{equation}
\widetilde{a}_{t}(x):=p^{h}(t,x,x)
\quad \mbox{ for }  t>0 \mbox{ and } x\in E.\nonumber
\end{equation}
Using the Poincar\'{e} inequality \eqref{p1}, we can prove that (cf. \cite{CRY}) for every $g\in L^{2}(E,\widetilde{m})$,
\begin{equation}
|P^{h}_{t}g(x)-\langle g,h^{2}\rangle |
\le e^{-\lambda_h (t-s)}\widetilde{a}_{2s}(x)^{1/2}\|g\|_{L^{2}(E,\widetilde{m})}.\label{p3}
\end{equation}
Moreover for $t>s>0$ and $x,y\in E$,
\begin{equation}
\left|p^{h}(t,x,y)-1\right|\le e^{-\lambda_{h}(t-s)}\widetilde{a}_{s}(x)^{1/2}\widetilde{a}_{s}(y)^{1/2}.\label{p2}
\end{equation}

For $\mu\in\mh$ with $\mu\not= 0$, we use $\Pi^{h}_{h\mu}$ to denote the probability measure where $\xi^{h}$ is the recurrent motion with
starting point
randomised according to $\frac{h(x)\mu(dx)}{\langle h,\mu\rangle}$. In other words,
$$\Pi^{h}_{h\mu}(\cdot)=\frac{1}{\langle h,\mu\rangle }\int_{E}\Pi^{h}_{x}(\cdot)h(x)\mu(dx).$$
Since $W^{h}_{t}(X)$ is a nonnegative $\p_{\mu}$-martingale, we can define a new probability measure $\q_{\mu}$ by
$$
d\q_{\mu}:=\frac{W^{h}_{t}(X)}{\langle h,\mu\rangle}\,d\p_{\mu}\quad\mbox{ on }\sigma(X_{s};s\in [0,t]).
$$

The next result can be proved in the same way as that for \cite[Lemma 2.17]{EKW}.

 \begin{lemma}\label{lem2.1} Suppose Assumption 2 holds.
 For all $\mu\in\mh$ with $\mu\not= 0$,  $f,g\in\B^{+}_{b}(E)$ and $t\ge 0$, we have
\begin{equation}
\q_{\mu}\left[e^{-\langle f,X_{t}\rangle}\frac{\langle gh,X_{t}\rangle}{\langle h,X_{t}\rangle}\right]=\p_{\mu}\left[e^{-\langle f,X_{t}\rangle}\right]
\Pi^{h}_{h\mu}\left[g(\xi_{t})\exp\left(-\int_{0}^{t}\frac{\partial \psi_{0}}{\partial \lambda}(\xi_{s},u_{f}(t-s,\xi_{s}))ds\right)\right],
 \nonumber
\end{equation}
where $u_{f}(t,x)$ is the unique nonnegative locally bounded solution to the integral equation \eqref{eq1}.
\end{lemma}

\medskip

It follows from Lemma \ref{lem2.1} that
\begin{eqnarray}
&&\q_{\mu}\left[e^{-\langle f,X_{t}\rangle}\right]\nonumber\\
&=&\p_{\mu}\left[e^{-\langle f,X_{t}\rangle}\right]
\Pi^{h}_{h\mu}\left[\exp\left(-\int_{0}^{t}\frac{\partial \psi_{0}}{\partial \lambda}(\xi_{s},u_{f}(t-s,\xi_{s}))\right)\right]\nonumber\\
&=&\p_{\mu}\left[e^{-\langle f,X_{t}\rangle}\right]\Pi^{h}_{h\mu}\left[\exp\left(-\int_{0}^{t}2\alpha(\xi_{s})
u_{f}(t-s,\xi_{s})ds-\int_{0}^{t}\int_{(0,{\infty})}\left(1-e^{-u_{f}(t-s,\xi_{s})y}\right)y\pi(\xi_{s},dy)ds\right)\right].\nonumber
\end{eqnarray}
This formula offers a probabilistic view of the superprocess $X_{t}$ under the new measure $\q_{\mu}$ that is stated in the following proposition.

\begin{proposition}\label{prop3}
Suppose Assumption 2 holds and that $E_{+}\subset E_{0}$.
For every $\mu\in\mh$ with $\mu\not= 0$, there exists a probability space with probability measure $\qq_{h\mu}$ that carries the following processes
\begin{description}
\item{\rm(i)}
$((\xi_{t})_{t\ge 0};\qq_{h\mu})$ is equal in distribution to $((\xi_{t})_{t\ge 0};\Pi^{h}_{h\mu})$.
We call $((\xi_{t})_{t\ge 0};\qq_{h\mu})$ the \textit{spine}.

\item{\rm(ii)}
$(n;\qq_{h\mu})$ is a random measure such that given
$(\xi, \qq_{h\mu})$,
 $n$ is a Poisson random measure which issues $\mf$-valued processes $X^{n,t}:=\{X^{n,t}_{s}; s\ge 0\}$
 at space-time point $(\xi_{t},t)$ with rate
    $$ d\mathbb{N}_{\xi_{t}}\times 2\alpha(\xi_{t})dt.$$
Here for every $x\in E_{+}=\{x\in E:\ \alpha(x)>0\}$, $\mathbb{N}_{x}$ is the Kuznetsov measure on $\mathcal{W}^{+}_{0}$ corresponding to the $(P_{t},\psi_{\beta})$-superprocess, while for $x\not\in E_{+}$ $\mathbb{N}_{x}$
is the null
measure on $\mathcal{W}^{+}_{0}$.
Note that, given $\xi$, the immigration happens only at space-time point $(\xi_{t},t)$ with $\alpha(\xi_{t})>0$.
    Let $D^{n}$ denote the almost surely countable set of immigration times, and $D^{n}_{t}:=D^{n}\cap [0,t]$.
    Given $\xi$, the processes $\{X^{n,t}:t\in D^{n}\}$ are mutually independent. We refer to $\{X^{n,t}:t\in D^{n}\}$ as the \textit{continuous immigration}.

 \item{\rm(iii)}
$(m;\qq_{h\mu})$ is a random measure such that given $\xi$,
 $m$ is a Poisson random measure which issues $\mf$-valued processes $X^{m,t}:=(X^{m,t}_{s})_{s\ge 0}$ at space-time point $(\xi_{t},t)$ with initial mass $y$ at rate
 $$y\pi(\xi_{t},dy)\times d\p_{y\delta_{\xi_{t}}}\times dt.$$
 Here $\p_{y\delta_{x}}$ denotes the law of the $(P_{t},\psi_{\beta})$-superprocess starting in $y\delta_{x}$.
 Let $D^{m}$ denote the almost surely countable set of immigration times, and $D^{m}_{t}:=D^{m}\cap [0,t]$.
 Given $\xi$, the processes $\{X^{m,t}:t\in D^{m}\}$ are mutually independent and independent of $n$ and $\{X^{n,t}:t\in D^{n}\}$. We refer to $\{X^{m,t}:t\in D^{m}\}$ as the \textit{discontinuous immigration}.

  \item{\rm(iv)}
 $((X_{t})_{t\ge 0};\qq_{h\mu})$ is equal in distribution to $((X_{t})_{t\ge 0};\p_{\mu})$. Moreover $((X_{t})_{t\ge 0};\qq_{h\mu})$ is independent of $\xi$, $n$, $m$ and all immigration processes.

\end{description}

Let $X^{n}_{t}:=\sum_{s\in D^{n}_{t}}X^{n,s}_{t-s}$ and $X^{m}_{t}:=\sum_{s\in D^{m}_{t}}X^{m,s}_{t-s}$. Define $$\Gamma_{t}:=X_{t}+X^{n}_{t}+X^{m}_{t} \mbox{ for }t\ge 0.$$
Then $(\Gamma:=(\Gamma_{t})_{t\ge 0};\qq_{h\mu})$ is equal in distribution to $(X:=(X_{t})_{t\ge 0};\q_{\mu})$.
\end{proposition}

\medskip

The proof of Proposition \ref{prop3} is similar to that of \cite[Theorem 5.2]{KLMR}
so it is omitted here.
For $s\ge 0$, define
$$
I^{m}_{s}:=\langle 1,X^{m,s}_{0}\rangle \quad\mbox{ if }s\in D^{m}\quad \mbox{ and }\quad I^{m}_{s}:=0\quad\mbox{elsewise}.
$$
Then, given $\xi$,  $\{I^{m}_{s},s\ge 0\}$  is a Poisson point process with characteristic measure $y\pi(\xi_{s},dy)$.
Let $\mathcal{G}$ be the $\sigma$-field generated by $\xi$, the random measures $n$ and $m$, and the process $\{I^{m}_{s},s\ge 0\}$.

\begin{lemma}
For $\mu\in\mh$ with
$\mu\not=0$,
$f\in\B^{+}_{b}(E)$ and $t\ge 0$,
\begin{equation}
\qq_{h\mu}\left[\langle f,\Gamma_{t}\rangle |
\mathcal{G}
\right]=
\langle P^{\beta}_{t}f,\mu\rangle +\sum_{s\in D^{n}_{t}}P^{\beta}_{t-s}f(\xi_{s})+
\sum_{s\in D^{m}_{t}}I^{m}_{s}P^{\beta}_{t-s}f(\xi_{s})\quad\qq_{h\mu}\mbox{-a.s.}\label{sd2}
\end{equation}
\end{lemma}

\proof
By \eqref{N}, we have for every $x\in E_{0}$, $f\in\B^{+}_{b}(E)$ and $t>0$,
$$\mathbb{N}_{x}\left(\langle f,X_{t}\rangle\right)=\p_{\delta_{x}}\left(\langle f,
X_{s}\rangle\right)=P^{\beta}_{t}f(x).$$
Thus by the definition of $\Gamma_{t}$, under $\qq_{h\mu}$,
\begin{eqnarray}
\qq_{h\mu}\left[\langle f,\Gamma_{t}\rangle|\mathcal{G}\right]
&=&\qq_{h\mu}\left(\langle f,X_{t}\rangle\right)+\sum_{s\in D^{n}_{t}}\qq_{h\mu}\left[\langle f,X^{n,s}_{t-s}\rangle|\mathcal{G}\right]+\sum_{s\in D^{m}_{t}}\qq_{h\mu}\left[\langle f,X^{m,s}_{t-s}\rangle|\mathcal{G}\right]\nonumber\\
&=&\p_{\mu}\left(\langle f,X_{t}\rangle\right)+\sum_{s\in D^{n}_{t}}\mathbb{N}_{\xi_{s}}\left(\langle f,X_{t-s}\rangle\right)+\sum_{s\in D^{m}_{t}}\p_{I^{m}_{s}\delta_{\xi_{s}}}\left(\langle f,X_{t-s}\rangle\right)\nonumber\\
&=&\langle P^{\beta}_{t}f,\mu\rangle+\sum_{s\in D^{n}_{t}}P^{\beta}_{t-s}f(\xi_{s})+\sum_{s\in D^{m}_{t}}I^{m}_{s}P^{\beta}_{t-s}f(\xi_{s}).\nonumber
\end{eqnarray}\qed

\medskip

By the monotone convergence theorem,
\eqref{sd2} holds for any $f\in\mathcal{B}^{+}(E)$.
We call formula \eqref{sd2} the \textit{spine decomposition} for $((X_{t})_{t\ge 0};\q_{\mu})$.

\medskip

\section{Martingale convergence of $W^{h}_{t}(X)$}\label{sec4}

In this section we prove Theorem \ref{them2}
 under Assumptions 1-3.

\begin{lemma}\label{lem2.2}
$g(x):=h(x)^{-1}\p_{\delta_{x}}\left(W^{h}_{\infty}(X)\right)$
is a constant function on $E$.
\end{lemma}

\proof We first claim that for any $\mu\in\mh$,
\begin{equation}\label{2.2.1}
\p_{\mu}\left(W^{h}_{\infty}(X)\right)=\int_{E}\p_{\delta_{x}}\left(W^{h}_{\infty}(X)\right)\mu(dx).
\end{equation}
The proof of this claim is the same as that of \cite[Lemma 5.4]{RSY}, we omit its details here.
By the Markov property of $X$,
we have for any $t\ge 0$ and $x\in E$,
\begin{eqnarray}
g(x)&=&\frac{1}{h(x)}\p_{\delta_{x}}\left[\lim_{s\to{\infty}}e^{\lambda_{1}(t+s)}\langle h,X_{t+s}\rangle\right]\nonumber\\
&=&\frac{e^{\lambda_{1}t}}{h(x)}\p_{\delta_{x}}\left[\p_{X_{t}}\left(\lim_{s\to{\infty}}W^{h}_{s}(X)\right)\right]\nonumber\\
&=&\frac{e^{\lambda_{1}t}}{h(x)}\p_{\delta_{x}}\left[\p_{X_{t}}\left(W^{h}_{\infty}(X)\right)\right]
\nonumber.
\end{eqnarray}
This together with \eqref{2.2.1} implies that
\begin{equation}
g(x)=\frac{e^{\lambda_{1}t}}{h(x)}\p_{\delta_{x}}\left(\langle gh,X_{t}\rangle\right)
=\frac{e^{\lambda_{1}t}}{h(x)}P^{\beta}_{t}(gh)(x)=P^{h}_{t}g(x).\nonumber
\end{equation}
This means that $g$ is an invariant function for the irreducible recurrent process $(\xi^{h};\Pi^{h})$.
Hence $g$ is a constant function on $E$.\qed

\medskip

\noindent\textit{Proof of Theorem \ref{them2}:}
Without loss of generality, we assume $\mu\in\mh$ and $\mu\not= 0$.
Since $W^{h}_{t}(X)$ is a nonnegative martingale,
to show it is a closed martingale,
it suffices to prove
 \begin{equation}\label{4.7}
\p_{\mu}\left(W^{h}_{\infty}(X)\right)=\langle h,\mu\rangle.
\end{equation}
 First we claim that \eqref{4.7} is true for
 $ \mu_{B}(dy):=
 1_{B}(y)h(y)m(dy)$ with $B\Subset E$ and $m(B)\not= 0$.
 It is  straightforward  to see from the change of measure methodology (see, for example,
 \cite[Theorem 5.3.3]{Durrett}
 ) that the proof for this claim is completed as soon as we can show that
\begin{equation}
\q_{\mu_{B}}\left(\limsup_{t\to{\infty}}W^{h}_{t}(X)<{\infty}\right)=1.\label{4.8}
\end{equation}
Since $((X_{t})_{t\ge 0};\q_{\mu_{B}})$ is equal in distribution to $((\Gamma_{t})_{t\ge 0};\qq_{h\mu_{B}})$, \eqref{4.8} is equivalent to that
\begin{equation}
\qq_{h\mu_{B}}\left(\limsup_{t\to{\infty}}W^{h}_{t}(\Gamma)<{\infty}\right)=1.\label{4.9}
\end{equation}
Recall that $\mathcal{G}$ is the $\sigma$-field generated by $\xi$, $n$, $m$ and $\{I^{m}_{s},s\ge 0\}$.
By the spine decomposition formula \eqref{sd2}, for any $t>0$,
\begin{eqnarray}
\qq_{h\mu_{B}}\left(W^{h}_{t}(\Gamma)|\mathcal{G}\right)
&=&e^{\lambda_{1}t}\langle P^{\beta}_{t}h,\mu_{B}\rangle +e^{\lambda_{1}t}\sum_{s\in D^{n}_{t}}P^{\beta}_{t-s}h(\xi_{s})
+e^{\lambda_{1}t}\sum_{s\in D^{m}_{t}}I^{m}_{s}P^{\beta}_{t-s}h(\xi_{s})\nonumber\\
&=&\langle h,\mu_{B}\rangle +\sum_{s\in D^{n}_{t}}e^{\lambda_{1}s}h(\xi_{s})+\sum_{s\in D^{m}_{t}}e^{\lambda_{1}s}I^{m}_{s}h(\xi_{s})\nonumber\\
&\le & \langle h,\mu_{B}\rangle +\sum_{s\in D^{n}}e^{\lambda_{1}s}h(\xi_{s})+\sum_{s\in D^{m}}e^{\lambda_{1}s}I^{m}_{s}h(\xi_{s}).\label{5.8}
\end{eqnarray}
Applying almost the same argument as in the proof of \cite[Theorem 5.1]{RSY}, we can show that the last two terms in \eqref{5.8} are finite almost surely,
and hence $\limsup_{t\to {\infty}}\qq_{h\mu_{B}}\left(W^{h}_{t}(\Gamma)|\mathcal{G}\right)<{\infty}$
$\qq_{h\mu_{B}}$-a.s. By Fatou's lemma, $\qq_{h\mu_{B}}\left(\liminf_{t\to{\infty}}W^{h}_{t}(\Gamma)<{\infty}\right)=1$. Note that $W^{h}_{t}(\Gamma)^{-1}$ is a nonnegative $\qq_{h\mu_{B}}$-supermartingale. Hence $$\qq_{h\mu_{B}}\left(\limsup_{t\to{\infty}}W^{h}_{t}(\Gamma)<{\infty}\right)=1.$$
This proves \eqref{4.9} and consequently,
\begin{equation}
\p_{\mu_{B}}\left(W^{h}_{\infty}(X)\right)=\langle h,\mu_{B}\rangle.\label{4.15}
\end{equation}
Note that $\p_{\mu_{B}}\left(W^{h}_{\infty}(X)\right)=\langle \p_{\delta_{x}}\left(W^{h}_{\infty}(X)\right),\mu_{B}\rangle$ by \eqref{2.2.1}
 and then
  \begin{equation}\label{4.15'}\langle \p_{\delta_{x}}\left(W^{h}_{\infty}(X)\right),\mu_{B}\rangle=\langle h,\mu_{B}\rangle.\end{equation}
 By Fatou's lemma, for every $x\in E$,
  $\p_{\delta_{x}}\left(W^{h}_{\infty}(X)\right)\le h(x)$.
 We get by
 \eqref{4.15'}
  that $\p_{\delta_{x}}\left(W^{h}_{\infty}(X)\right)= h(x)$ $m$-a.e. on $B$. Since $B$ is arbitrary, $\p_{\delta_{x}}\left(W^{h}_{\infty}(X)\right)= h(x)$ $m$-a.e. on $E$ and hence everywhere on $E$ by Lemma \ref{lem2.2}. Therefore by \eqref{2.2.1}, $\p_{\mu}\left(W^{h}_{\infty}(X)\right)=\langle h,\mu\rangle$ holds for all $\mu\in\mh$. This completes the proof.
\qed

\medskip

\section{Skeleton space}\label{sec5}

In this section, we  work under Assumptions 1-3.
We first consider  general superprocess case, for which as a part of Assumption 1,
Condition 1 holds and $w$ is the bounded positive function in \eqref{a3}.
 By \eqref{eq0} and \eqref{eq1}, $u_w (t, x)=w(x)$
is the unique nonnegative bounded solution to the following equation:
\begin{equation}\label{e:5.1}
w(x) =P_t w(x) - \int_0^t P_s (\psi_\beta (\cdot, w (\cdot)) (x) ds
\quad \hbox{for every } x\in E \hbox{ and } t\geq 0.
\end{equation}
Note that
\begin{equation}\label{e:5.2}
\frac{\psi_\beta (x, \lambda)}{\lambda} = -\beta(x)+ a(x)\lambda  + \frac1{\lambda} \int_{(0, \infty)} \left(e^{-\lambda y} -1 + \lambda y \right) \pi (x, dy),
\end{equation}
which is a  bounded function function on $E\times (0, M]$ for every $M>0$.
Let
$$
M_t:=w(\xi_t)-\int_0^t \psi_\beta (\xi_s, w(\xi_s)) ds,
$$
which is bounded in $t\in [0, T]$ for every $T>0$.
 Let $\{\mathcal{H}_{t}: t\geq 0\}$ be the minimal augmented $\sigma$-fields generated by the process $\xi=\{\xi_t, t\geq 0; \Pi_x, x\in E\}$.
By the Markov property of $\xi_t$ and \eqref{e:5.1}, we have for every $t> s\geq 0$ and $x\in E$,
\begin{eqnarray*}
\Pi_x \left[ M_t | \mathcal{H}_{s}\right]
&=& \Pi_{\xi_s}\left[ w(\xi_{t-s}\right]
 - \int_0^s \psi_\beta (\xi_r, w(\xi_r)) dr
-\Pi_{\xi_s} \left[  \int_0^{t-s} \psi_\beta (\xi_r, w(\xi_r)) dr \right] \\
&=& P_{t-s}w(\xi_s) - \Pi_{\xi_s} \left[ \int_0^{t-s}   \psi_\beta (\cdot , w(\cdot )) dr \right]
 -  \int_0^s \psi_\beta (\xi_r, w(\xi_r)) dr \\
&=& w(\xi_s) -   \int_0^s \psi_\beta (\xi_r, w(\xi_r)) dr = M_s.
\end{eqnarray*}
In other words, $\{M_t; t\geq 0\}$ is a martingale additive functional of $\xi$.
Observe that it follows from \eqref{e:5.2} that
$$
\gamma(x):=\frac{\psi_{\beta}(x ,w(x))}{w(x)}
$$
 is a bounded function.
Since $w(\xi_t)= M_t+  \int_0^t \psi_\beta (\xi_s, w(\xi_s)) ds$,
 we have by Ito's formula
\begin{eqnarray*}
d\left( w(\xi_t) e^{-\int_0^t \gamma (\xi_s) ds} \right)
&=& -w(\xi_t) \gamma (\xi_t)  e^{-\int_0^t \gamma (\xi_s) ds}  dt
+ e^{-\int_0^t \gamma (\xi_s) ds}  (dM_t + \psi_\beta (\xi_t, w(\xi_t)) dt \\
&=& e^{-\int_0^t \gamma (\xi_s) ds}    dM_t.
\end{eqnarray*}
Thus we have shown that $t\mapsto w(\xi_t) e^{-\int_0^t \gamma (\xi_s) ds}$ is a martingale.

Next we consider the
superdiffusion case, for which Assumption 1 is just Condition 1';
that is, there is a   locally bounded positive function $w$ satisfying \eqref{a3}.
 By truncating $w$ and then using the monotone convergence theorem, one can prove that
\begin{equation}\nonumber
w(x)=\Pi_{x}\left[w(\xi_{t\wedge \tau_{B}})\right]-\Pi_{x}\left[\int_{0}^{t\wedge \tau_{B}}\psi_{\beta}\left(\xi_{s},w(\xi_{s})\right)ds\right]
\end{equation}
for any bounded open set $B\Subset E$. Since $x\mapsto w(x)$ is bounded on $B$, one can generalize the above argument and show that
$t\mapsto w(\xi_{t\wedge\tau_{B}})\exp\left(-\int_{0}^{t\wedge \tau_{B}}\gamma(\xi_{s})ds\right)$ is a nonnegative $\Pi_{x}$-martingale.
This implies that $\{w(\xi_{t})\exp\left(-\int_{0}^{t}\gamma(\xi_{s})ds\right):t\ge 0\}$ is a nonnegative $\Pi_{x}$-local martingale and hence a $\Pi_{x}$-supermartingale.

Hence in both cases we can define a family of (sub)probability measures
$\{\Pi^{w}_{x}:x\in E\}$ on $\mathcal{H}$ by
$$
d\Pi^{w}_{x}
:=\frac{w(\xi_{t})}{w(x)}\exp\left(-\int_{0}^{t}\frac{\psi_{\beta}(\xi_{s},w(\xi_{s}))}{w(\xi_{s})}\,ds\right) d\Pi_{x}
\quad\mbox{ on }\mathcal{H}_{t}
\quad\mbox{ for every } t\ge 0.
$$
Clearly, the transformed process $\xi^{w}:=((\xi_{t})_{t\ge 0};\Pi^{w}_{x},x\in E)$ is a Borel right process.
Through the arguments after Remark \ref{rm for them2},
we can also construct the skeleton space for a $(P_{t},\psi_{\beta})$-superprocesses as follows.

\begin{proposition}\label{prop2}
 Suppose Assumptions 1-3 hold.
Let $\mathcal{M}^{loc}_{A}(E)$ denote the set of locally finite integer-valued measures on $(E,\mathcal{B}(E))$. For every $\mu\in \mf$ and every $\nu\in\mathcal{M}^{loc}_{A}(E)$, there exists a probability space with probability measure $\pp_{\mu,\nu}$ that carries the following processes:
\begin{description}
\item{\rm(i)}
$\left(Z:=(Z_{t})_{t\ge 0};\pp_{\mu,\nu}\right)$ is a supercritical branching Markov process with spatial motion $\xi^{w}$, branching rate function $q(x)$ and offspring distribution function $\{p_{k}(x):k\ge 2\}$ uniquely defined by
    \begin{equation}\label{G}
    G(x,s):=q(x)\sum_{k=2}^{{\infty}}p_{k}(x)(s^{k}-s)
    :=
    \frac{1}{w(x)}
    \left(
    \psi_{0}(x,(1-s)w(x))
    -(1-s)\psi_{0}(x,w(x))
    \right),
    \end{equation}
   and
   $\pp_{\mu,\nu}(Z_{0}=\nu)=1$. We use the classical Ulam-Harris notations to refer to the particles in the genealogical tree $\mathcal{T}$ of $Z$. For a particle $u\in \mathcal{T}$, we use $b_{u}$ and $d_{u}$ for its birth and death times, and $z_{u}(t)$ for its spatial location at time $t\in [b_{u},d_{u}]$.

\item{\rm(ii)}
$\left(X^{*}:=(X^{*}_{t})_{t\ge 0};\pp_{\mu,\nu}\right)$ is a
$(P_{t},\psi^{*}_{\beta^{*}})$-superprocess
with $\pp_{\mu,\nu}(X^{*}_{0}=\mu)=1$ such that for all $\mu\in\mf$, $f\in\mathcal{B}^{+}_{b}(E)$ and $t\ge 0$,
    \begin{equation}
    \pp_{\mu,\nu}\left[e^{-\langle f,X^{*}_{t}\rangle}\right]=e^{-\langle u^{*}_{f}(t,\cdot),\mu\rangle},\label{3.3}
    \end{equation}
 where $u^{*}_{f}(t,x)$ is the unique nonnegative locally bounded solution to
    equation \eqref{1.8},
    and
    $\psi^{*}_{\beta^{*}}(x,\lambda):=-\beta^{*}(x)\lambda+\psi^{*}_{0}(x,\lambda)$
    is given in Proposition \ref{prop2.4}.
    By \eqref{1.8} and \eqref{3.3}, the mean of $\langle f,X^{*}_{t}\rangle$ can
    be expressed by
    \begin{equation}
    \pp_{\mu,\nu}\left(\langle f,X^{*}_{t}\rangle\right)
    =
    \langle P^{\beta^{*}}_{t}f,\mu\rangle\quad\mbox{ for }f\in\mathcal{B}^{+}_{b}(E).\nonumber
    \end{equation}
    The distribution of $X^{*}$ under $\pp_{\mu,\nu}$ does not depend on $\nu$. Moreover, under $\pp_{\mu,\nu}$, $X^{*}$ is independent of $Z$.

\item{\rm(iii)}
$\left(I:=(I_{t})_{t\ge 0};\pp_{\mu,\nu}\right)$ is a $\mf$-valued process with $\pp_{\mu,\nu}(I_{0}=0)=1$, which
is given by $I=I^{\ia}+I^{\ib}$ where $I^{\ia}$ and $I^{\ib}$ are described as follows:
    \begin{description}
    \item{\rm(a)}
    $(\ia;\pp_{\mu,\nu})$ is a random measure such that given $Z$, $\ia$ is a Poisson random measure that issues for every $u\in\mathcal{T}$, $\mf$-valued processes $X^{\ia,u,r}:=(X^{\ia,u,r}_{t})_{t\ge 0}$ along the space-time trajectory $\{(z_{u}(r),r):r\in (b_{u},d_{u}]\}$ with rate
$$
		\left(2\alpha(z_{u}(r))d\mathbb{N}^{*}_{z_{u}(r)}+\int_{(0,{\infty})}y\pi^{*}(z_{u}(r),dy)\times d\p^{*}_{y\delta_{z_{u}(r)}}\right)\times dr,
$$
    where $\p^{*}_{\mu}$ is the law of the $(P_{t},\psi^{*}_{\beta^{*}})$-superprocess starting in $\mu$, for every $x\in E_{+}= \{x\in E:\ \alpha(x)>0\}$, $\mathbb{N}^{*}_{x}$ denotes the Kuznetsov measure on $\mathcal{W}^{+}_{0}$ corresponding to the $(P_{t},\psi^{*}_{\beta^{*}})$-superprocess, and for $x\not\in E_{+}$,  $\mathbb{N}^{*}_{x}$
is the null
measure on $\mathcal{W}^{+}_{0}$.
    The processes $\{X^{\ia,u,r}:u\in\mathcal{T},\ r\in (b_{u},d_{u}]\}$ are independent of $X^{*}$, and, given $Z$, are mutually independent.
    $$I^{\ia}_{t}:=\sum_{u\in\mathcal{T}}\sum_{r\in (b_{u},d_{u}\wedge t]}X^{\ia,u,r}_{t-r}\quad\mbox{ for }t\ge 0.$$
		
\item{\rm(b)}
    $(\ib;\pp_{\mu,\nu})$ is a random measure such that given $Z$, $\ib$ issues for every $u\in\mathcal{T}$, at space-time point $(z_{u}(d_{u}),d_{u})$,
        an
        $\mf$-valued process $X^{\ib,u}:=(X^{\ib,u}_{t})_{t\ge 0}$ with law $\p^{*}_{Y_{u}\delta_{z_{u}(d_{u})}}$ such that given that $u$ gives birth to $k$ particles at its death time, the independent $\mathbb{R}^{+}$-valued random variable $Y_{u}$ is distributed according to the measure
        $$\left.\frac{1}{q(x)w(x)p_{k}(x)}\left(\alpha(x)w(x)^{2}\delta_{0}(dy)1_{\{k=2\}}+
        w(x)^{k}\frac{y^{k}}{k!}\pi^{*}(x,dy)\right)\right|_{x=z_{u}(d_{u})}.$$
        The processes $\{X^{\ib,u}:u\in\mathcal{T}\}$ are independent of $X^{*}$, and, given $Z$, are mutually independent and independent of $\ia$.
        $$I^{\ib}_{t}:=\sum_{u\in\mathcal{T}}1_{\{d_{u}\le t\}}X^{\ib,u}_{t-d_{u}}\quad\mbox{ for }t\ge 0.$$
    \end{description}
          \end{description}
For every $\mu\in\mf$, let $\pp_{\mu}$ denote the measure $\pp_{\mu,\nu}$ with $\nu$ replaced by a Poisson random measure with intensity $w(x)\mu(dx)$.
Then $\left(
\widehat X :=
X^{*}+I;\pp_{\mu}\right)$ is Markovian and
has the same distribution as
$(X;\p_{\mu})$.
Moreover, under $\pp_{\mu}$, given $\widehat  X_{t}$,
the measure $Z_{t}$ is a Poisson random measure with intensity
$w(x)\widehat X_{t}(dx)$.
\end{proposition}

With the results from Subsection \ref{sup-Br} and Appendix \ref{sup-diffusion}
on the existence of ${\mathbb N}^*_x$ under Assumption 1,
the proof of Proposition \ref{prop2} is very similar to
that of \cite[Corollary 2]{KPR} and we omit
its long computations here.
We call the probability space
in Proposition \ref{prop2}
 the \textit{skeleton space}. The process $Z$ is called \textit{skeleton process} and $I$ is called \textit{immigration process}.
We call $(X^{*}+I;\pp_{\mu} )$ the skeletion decomposition of $X$. Since
$({\widehat X};\pp_{\mu})$
is equal in distribution to the $(P_{t},\psi_{\beta})$-superprocess $(X;\p_{\mu})$, we may work on this skeleton space whenever it is convenient.
For notational simplification, we will abuse the notation and denote $\widehat X$ by $X$.
 Since the distributions of $X^{*}$ (resp. $I$) under $\pp_{\mu,\nu}$ do not depend on $\nu$ (resp. $\mu$), we sometimes write $\pp_{\mu,\cdot}$ (resp. $\pp_{\cdot,\nu}$) for $\pp_{\mu,\nu}$.

For $t\ge 0$, we write $Z_{t}=\sum_{i=1}^{N_{t}}\delta_{z_{i}(t)}$ where $N_{t}$ denotes the number of skeleton particles at time $t$ and $\{z_{i}(t):i=1,\cdots,N_{t}\}$ their spatial locations. Let $m(x):=\sum_{k=2}^{{\infty}}kp_{k}(x)$. We have by \eqref{G} that $q(x)(m(x)-1)=\frac{\partial G}{\partial s}(x,s)|_{s=1}=\psi_{0}(x,w(x))/w(x)
$. Then by the many-to-one formula for branching Markov process (cf. \cite[Lemma 3.3]{Shiozawa}), for $f\in\mathcal{B}^{+}_{b}(E)$, $x\in E$ and $t\ge 0$,
\begin{eqnarray}
\pp_{\cdot,\delta_{x}}\left(\langle f,Z_{t}\rangle\right)
&=&\Pi^{w}_{x}\left[\exp\left(\int_{0}^{t}q(\xi_{s})(m(\xi_{s})-1)ds\right)f(\xi_{t})\right]\nonumber\\
&=&\frac{1}{w(x)}\Pi_{x}\left[\exp\left(\int_{0}^{t}\frac{-\psi_{\beta}(\xi_{s},w(\xi_{s}))+\psi_{0}(\xi_{s},w(\xi_{s}))}{w(\xi_{s})}
\,ds\right)w(\xi_{t})f(\xi_{t})\right]\nonumber\\
&=&\frac{1}{w(x)}P^{\beta}_{t}(wf)(x).\label{3.4}
\end{eqnarray}
Note that under $\pp_{\mu}$ with $\mu\in\mf$, $Z_{0}$ is a Poisson random measure with intensity $w(x)\mu(dx)$. We get from \eqref{3.4} that
\begin{equation}
\pp_{\mu}\left(\langle f,Z_{t}\rangle\right)=\pp_{\mu}\left(\sum_{i=1}^{N_{0}}\pp_{\cdot,\delta_{z_{i}(0)}}\left(\langle f,Z_{t}\rangle\right)\right)=\pp_{\mu}\left(\langle \frac{1}{w}P^{\beta}_{t}(wf),Z_{0}\rangle\right)=\langle P^{\beta}_{t}(wf),\mu\rangle.\label{meanz}
\end{equation}
For $t\ge 0$, let $\mathcal{F}_{t}$ denote the $\sigma$-filed generated by $Z$, $X^{*}$ and $I$ up to time $t$. Denote by $I^{*,t}_{s}$ the immigration at time $t+s$ that occurred along the skeleton before time $t$. For $i\in\{1,2,\cdots,N_{t}\}$, denote by $I^{i,t}_{s}$ the immigration at time $t+s$ that occurred along the subtree of the skeleton rooted at the $i$th particle with location $z_{i}(t)$. We have
$$X_{s+t}=X^{*}_{s+t}+I^{*,t}_{s}+\sum_{i=1}^{N_{t}}I^{i,t}_{s}\quad\mbox{ for all }s,t\ge 0.$$
We know by \cite{EKW} that, given $\mathcal{F}_{t}$, $(X^{*}_{s+t}+I^{*,t}_{s})_{s\ge 0}$ is equal in distribution to
$((X^{*}_{s})_{s\ge 0};\pp_{X_{t}})$
and $I^{i,t}:=(I^{i,t}_{s})_{s\ge 0}$ is equal in distribution to $(I;\pp_{\cdot,\delta_{z_{i}(t)}})$. Moreover, the processes $\{I^{i,t}:\ i=1,2,\cdots, N_{t}\}$ are mutually independent. Using these properties, we have for $\mu\in\mf$, $f\in\mathcal{B}^{+}_{b}(E)$ and $t,s\ge 0$, under $\pp_{\mu}$
\begin{eqnarray}
\pp_{\mu}\left(\langle f,X_{t+s}\rangle|\mathcal{F}_{t}\right)&=&\pp_{\mu}\left(\langle f,X^{*}_{t+s}+I^{*,t}_{s}\rangle|\mathcal{F}_{t}\right)+\pp_{\mu}\left(
\left.  \sum_{i=1}^{N_{t}}\langle f,I^{i,t}_{s}\rangle
\right|\mathcal{F}_{t}\right)\nonumber\\
&=&\pp_{X_{t}}\left(\langle f,X^{*}_{s}\rangle\right)+\sum_{i=1}^{N_{t}}\pp_{\cdot,\delta_{z_{i}(t)}}\left(\langle f,I_{s}\rangle\right).\label{2.20}
\end{eqnarray}
Note that under $\pp_{\mu}$, $N_{0}$ is a Poisson random variable with mean $\langle w,\mu\rangle$.
We have for each $x\in E$,
\begin{eqnarray}
\pp_{\delta_{x}}\left(\langle f,I_{t}\rangle\right)
&=&
\pp_{\delta_{x}}\left(\sum_{i=1}^{N_{0}}\pp_{\cdot,\delta_{z_{i}(0)}}\left(\langle f,I_{t}\rangle\right)\right)\nonumber\\
&=&\pp_{\delta_{x}}\left(N_{0}\pp_{\cdot,\delta_{x}}\left(\langle f,I_{t}\rangle\right)\right)=w(x)\pp_{\cdot,\delta_{x}}\left(\langle f,I_{t}\rangle\right).\label{2.19}
\end{eqnarray}
Thus $\pp_{\cdot,\delta_{x}}\left(\langle f,I_{t}\rangle\right)=w(x)^{-1}\pp_{\delta_{x}}\left(\langle f,I_{t}\rangle\right)=w(x)^{-1}\pp_{\delta_{x}}\left(\langle f,X_{t}-X^{*}_{t}\rangle\right)=w(x)^{-1}(P^{\beta}_{t}f(x)-P^{\beta^{*}}_{t}f(x))$.
Using this and \eqref{2.20}, we have for every $\mu\in\mf$, $f\in\mathcal{B}^{+}_{b}(E)$ and $s,t\ge 0$,
\begin{equation}\label{sd1}
\pp_{\mu}\left[\langle f,X_{t+s}\rangle|\mathcal{F}_{t}\right]=\langle P^{\beta^{*}}_{s}f,X_{t}\rangle+\langle
\frac{P^{\beta}_{s}f}{w},Z_{t}\rangle-\langle\frac{P^{\beta^{*}}_{s}f}{w},Z_{t}\rangle\quad\pp_{\mu}\mbox{-a.s.}
\end{equation}
By the monotone convergence theorem, the boundedness of $f$ in \eqref{sd1} is unnecessary and \eqref{sd1} holds for $f\in\mathcal{B}^{+}(E)$
which satisfies $\langle P^{\beta}_{t}f,\mu\rangle <{\infty}$ for any $t\ge 0$,

\begin{lemma}\label{lem1.1}
Suppose Assumptions 1-3 hold.
For every $f\in\mathcal{B}^{+}(E)$
so that $f/h$ is bounded,
define
$$\theta^{*}_{f}(t,x):=\frac{e^{\lambda_{1}t}}{h(x)}P^{\beta^{*}}_{t}f(x)\quad\mbox{ for }x\in E\mbox{ and }t\ge 0.$$
Then $\theta^{*}_{f}(t,x)$ is bounded on $[0,{\infty})\times E$, and for every $x\in E$,
$$\theta^{*}_{f}(t,x)\to 0\mbox{  as }t\to{\infty}.$$
\end{lemma}

\proof Let  $g(x):=\beta(x)-\beta^{*}(x)=2\alpha(x)w(x)+\int_{(0,{\infty})}\left(1-e^{-w(x)y}\right)y\pi(x,dy)\in\mathcal{B}^{+}(E)$
and $c_{1}:=\sup_{x\in E}f(x)/h(x)\in [0,{\infty})$. By the definition, for every $x\in E$ and $t\ge 0$
\begin{eqnarray}
\theta^{*}_{f}(t,x)&=&\frac{e^{\lambda_{1}t}}{h(x)}\Pi_{x}\left[\exp\left(\int_{0}^{t}\beta^{*}(\xi_{s})ds\right)f(\xi_{t})\right]
\nonumber\\
&=&\frac{e^{\lambda_{1}t}}{h(x)}\Pi_{x}\left[\exp\left(\int_{0}^{t}\beta(\xi_{s})-g(\xi_{s})ds\right)h(\xi_{t})
\frac{f(\xi_{t})}{h(\xi_{t})}\right]\nonumber\\
&=&\Pi_{x}^{h}\left[\exp\left(-\int_{0}^{t}g(\xi_{s})ds\right)\frac{f(\xi_{t})}{h(\xi_{t})}\right]\nonumber\\
&\le&c_{1}\Pi_{x}^{h}\left[\exp\left(-\int_{0}^{t}g(\xi_{s})ds\right)\right].\nonumber
\end{eqnarray}
Immediately $\theta^{*}_{f}(t,x)\le c_{1}$ for all $(t,x)\in [0,{\infty})\times E$. Let $l(x):=\Pi_{x}^{h}\left[\exp\left(-\int_{0}^{{\infty}}g(\xi_{s})ds\right)\right]$. To prove the second claim of this lemma, it suffices to prove that $l(x)\equiv 0$ on $E$.
Let $A:=\{x\in E: g(x)>0\}=\{x\in E:\ \alpha(x)+\pi(x,(0,{\infty}))>0\}$. We have $\widetilde{m}(A)>0$ since $m(A)>0$.
For an arbitraty $M>0$ and $x\in E$, define $P_{M}(x):=\Pi^{h}_{x}\left(\int_{0}^{{\infty}}g(\xi_{s})ds>M\right)$.
By the Markov property of $\xi$,
for $t\ge 0$ and $x\in E$,
\begin{eqnarray}
P^{h}_{t}P_{M}(x)&=&\Pi^{h}_{x}\left[\Pi^{h}_{\xi_{t}}\left(\int_{0}^{{\infty}}g(\xi_{s})ds>M\right)\right]\nonumber\\
&=&\Pi^{h}_{x}\left(\int_{t}^{{\infty}}g(\xi_{s})ds>M\right)\nonumber\\
&\le&P_{M}(x).\nonumber
\end{eqnarray}
Moreover, $\lim_{t\to 0}P^{h}_{t}P_{M}(x)=P_{M}(x)$ by the monotone convergence theorem. Hence $P_{M}(x)$ is an excessive function for the irreducible recurrent process $(\xi^{h};\Pi^{h})$. It follows from \cite[Lemma 3.5 and Lemma A.2.17]{CF} that $P_{M}$ is a constant function on $E$.
We claim that there exists $M_{1}>0$ and $c_{2}\in (0,1]$ such that
\begin{equation}\label{5.0}
P_{M_{1}}(x)\ge c_{2}>0\quad\mbox{ for all }x\in E.
\end{equation}
If \eqref{5.0} is not true, that is, $P_{M}(x)\equiv 0$ for all $M>0$, then $\Pi^{h}_{x}\left(\int_{0}^{{\infty}}g(\xi_{s})ds=0\right)=1$. This contradicts the fact that $\Pi^{h}_{x}\left(\int_{0}^{{\infty}}g(\xi_{s})ds\right)=\int_{0}^{{\infty}}\int_{A}p^{h}(s,x,y)g(y)\widetilde{m}(dy)ds>0$ since $\widetilde{m}(A)>0$ and $p^{h}(s,x,\cdot)$ is positive on $A$.
In view of \eqref{5.0}, we have
$$l(x)=\Pi^{h}_{x}\left[\exp\left(-\int_{0}^{{\infty}}g(\xi_{s})ds\right)\right]\le c_{2}e^{-M_{1}}+(1-c_{2})=:c_{3}<1.$$
Thus by the Markov property of $\xi$,
for any $t>0$,
\begin{eqnarray}
l(x)&=&\Pi^{h}_{x}\left[\exp\left(-\int_{0}^{t}g(\xi_{s})ds\right)\exp\left(-\int_{t}^{{\infty}}g(\xi_{s})ds\right)\right]\nonumber\\
&=&\Pi^{h}_{x}\left[\exp\left(-\int_{0}^{t}g(\xi_{s})ds\right)l(\xi_{t})\right]\nonumber\\
&\le&c_{3}\Pi^{h}_{x}\left[\exp\left(-\int_{0}^{t}g(\xi_{s})ds\right)\right].\label{5.24}
\end{eqnarray}
Letting $t\to {\infty}$ in \eqref{5.24}, we get $l(x)\le c^{2}_{3}$ for all $x\in E$. Applying the above argument recursively for $n$ times, we get $l(x)\le c^{n}_{3}$ for all $x\in E$. Hence we conclude $l(x)\equiv 0$ by letting $n\to{\infty}$.\qed

\medskip

\noindent\textit{Proof of Theorem \ref{them1}:} We only need to prove the first claim. For any $t,s\ge 0$, by the skeleton decomposition \eqref{sd1}, we have under $\pp_{\mu}$
\begin{eqnarray}
W^{h}_{t}(X)&=&\pp_{\mu}\left(W^{h}_{t+s}(X)|\mathcal{F}_{t}\right)\nonumber\\
&=&e^{\lambda_{1}(t+s)}\pp_{\mu}\left(\langle h,X_{t+s}\rangle |\F_{t}\right)\nonumber\\
&=&e^{\lambda_{1}(t+s)}\left(\langle P^{\beta^{*}}_{s}h,X_{t}\rangle+\langle\frac{P^{\beta}_{s}h}{w},Z_{t}\rangle-\langle\frac{P^{\beta^{*}}_{s}h}{w},Z_{t}\rangle\right)\nonumber\\
&=&e^{\lambda_{1}t}\langle\theta^{*}_{h}(s,\cdot)h,X_{t}\rangle+e^{\lambda_{1}t}\langle\frac{h}{w},Z_{t}\rangle
-e^{\lambda_{1}t}\langle\theta^{*}_{h}(s,\cdot)\frac{h}{w},Z_{t}\rangle\nonumber\\
&=&e^{\lambda_{1}t}\langle\theta^{*}_{h}(s,\cdot)h,X_{t}\rangle+W^{h/w}_{t}(Z)
-e^{\lambda_{1}t}\langle\theta^{*}_{h}(s,\cdot)\frac{h}{w},Z_{t}\rangle.\label{them1.eq1}
\end{eqnarray}
Letting $s\to{\infty}$, the first term and the third term in \eqref{them1.eq1} converge to $0$ $\pp_{\mu}$-almost surely by Lemma \ref{lem1.1} and the dominated convergence theorem and hence we get $W^{h}_{t}(X)=W^{h/w}_{t}(Z)$ $\pp_{\mu}$-a.s.\qed

\medskip

\section{Weak law of large numbers}\label{sec6}

Throughout this section,  we assume  Assumptions 1-3 and Assumption 4(i) hold.

\begin{lemma}\label{lem3.1}
 For any
$\phi\in\mathcal{B}^{+}_{b}(E)$ and $t>0$,
\begin{equation}
\int_{E}\p_{\delta_{x}}\left(\langle \phi h,X_{t}\rangle\log^{+}\langle \phi h,X_{t}\rangle\right)h(x)m(dx)<{\infty}.\nonumber
\end{equation}
\end{lemma}

\proof Note that for $\phi\in\mathcal{B}^{+}_{b}(E)$ and $x\in E$,
\begin{eqnarray}
\p_{\delta_{x}}\left(\langle \phi h,X_{t}\rangle\log^{+}\langle \phi h,X_{t}\rangle\right)
&\le&e^{-\lambda_{1}t}\|\phi\|_{\infty}\p_{\delta_{x}}\left(e^{\lambda_{1}t}\langle h,X_{t}\rangle\log^{+}\langle \phi h,X_{t}\rangle\right)\nonumber\\
&=&e^{-\lambda_{1}t}\|\phi\|_{\infty}h(x)\q_{\delta_{x}}\left(\log^{+}\langle \phi h,X_{t}\rangle\right)\nonumber\\
&=&e^{-\lambda_{1}t}\|\phi\|_{\infty}h(x)\qq_{h\delta_{x}}\left(\log^{+}\langle \phi h,\Gamma_{t}\rangle\right).\label{5.1}
\end{eqnarray}
Using the fact $\log^{+}(a+b)\le \log^{+}a+\log^{+}b+\log 2$ for $a,b\ge 0$ and Jensen inequality, we have
\begin{eqnarray}
\qq_{h\delta_{x}}\left(\log^{+}\langle \phi h,\Gamma_{t}\rangle\right)
&=&\qq_{h\delta_{x}}\left[\qq_{h\delta_{x}}\left(\log \left(\langle \phi h,\Gamma_{t}\rangle \vee 1\right)|\mathcal{G}\right)\right]\nonumber\\
&\le&\qq_{h\delta_{x}}\left[\log\qq_{h\delta_{x}}\left(\langle \phi h,\Gamma_{t}\rangle \vee 1|\mathcal{G}\right)\right]\nonumber\\
&\le&\qq_{h\delta_{x}}\left[\log^{+}\left(\qq_{h\delta_{x}}\left(\langle \phi h,\Gamma_{t}\rangle |\mathcal{G}\right)+1\right)\right]\nonumber\\
&\le&\qq_{h\delta_{x}}\left[\log^{+}\qq_{h\delta_{x}}\left(\langle \phi h,\Gamma_{t}\rangle | \mathcal{G}\right)\right]+\log 2.\label{5.2}
\end{eqnarray}
Thus by \eqref{5.1} and \eqref{5.2}, we only need to prove that
\begin{equation}
\int_{E}\qq_{h\delta_{x}}\left[\log^{+}\qq_{h\delta_{x}}\left(\langle \phi h,\Gamma_{t}\rangle|\mathcal{G}\right)\right]h(x)^{2}m(dx)<{\infty}.\label{5.3}
\end{equation}
By the spine decomposition formula \eqref{sd2}, we have under $\qq_{h\delta_{x}}$
\begin{eqnarray}
\qq_{h\delta_{x}}\left(\langle \phi h,\Gamma_{t}\rangle|\mathcal{G}\right)&\le&\|\phi\|_{\infty}\qq_{h\delta_{x}}\left(\langle  h,\Gamma_{t}\rangle|\mathcal{G}\right)\nonumber\\
&=&\|\phi\|_{\infty}e^{-\lambda_{1}t}\left(h(x)+\sum_{s\in D^{n}_{t}}e^{\lambda_{1}s}h(\xi_{s})+\sum_{s\in D^{m}_{t}}e^{\lambda_{1}s}I^{m}_{s}h(\xi_{s})\right).\nonumber
\end{eqnarray}
Using the fact that $\log^{+}(ab)\le \log^{+}a+\log^{+}b$ for $a,b\ge 0$, we have
\begin{eqnarray}
\log^{+}\qq_{h\delta_{x}}\left(\langle \phi h,\Gamma_{t}\rangle|\mathcal{G}\right)
&\le&\log^{+}\left(e^{-\lambda_{1}t}\|\phi\|_{\infty}\right)+\log^{+}h(x)
+\sum_{s\in D^{n}_{t}}\log^{+}\left(e^{\lambda_{1}s}h(\xi_{s})\right)\nonumber\\
&&+\sum_{s\in D^{m}_{t}}\log^{+}\left(e^{\lambda_{1}s}I^{m}_{s}h(\xi_{s})\right)+\sum_{s\in D^{n}_{t}\cup D^{m}_{t}}\log 2.\nonumber
\end{eqnarray}
Assumption 4(i)
implies that $\int_{E}\log^{+}h(x)h(x)^{2}m(dx)<{\infty}$. Hence to prove \eqref{5.3},
it suffices to prove
\begin{equation}\label{5.4}
\int_{E}\qq_{h\delta_{x}}\left[\sum_{s\in D^{n}_{t}}\left(\log^{+}(e^{\lambda_{1}s}h(\xi_{s}))+1\right)+\sum_{s\in D^{m}_{t}}\left(\log^{+}(e^{\lambda_{1}s}I^{m}_{s}h(\xi_{s}))+1\right)\right]\widetilde{m}(dx)<{\infty}.
\end{equation}
Since $\lambda_{1}<0$, we have by Fubini's theorem and the $\widetilde{m}$-symmetry of $(\xi^{h};\Pi^{h})$ that
\begin{eqnarray}
&&\int_{E}\qq_{h\delta_{x}}\left[\sum_{s\in D^{n}_{t}}\left(\log^{+}(e^{\lambda_{1}s}h(\xi_{s}))+1\right)\right]\widetilde{m}(dx)\nonumber\\
&=&\int_{E}\qq_{h\delta_{x}}\left[\int_{0}^{t}2\alpha(\xi_{s})\left(\log^{+}(e^{\lambda_{1}s}h(\xi_{s}))+1\right)\right]\widetilde{m}(dx)\nonumber\\
&=&2\int_{0}^{t}ds\int_{E}\Pi^{h}_{x}\left[\alpha(\xi_{s})\left(\log^{+}(e^{\lambda_{1}s}h(\xi_{s}))+1\right)\right]\widetilde{m}(dx)\nonumber\\
&=&2\int_{0}^{t}ds\int_{E}\alpha(x)\left(\log^{+}(e^{\lambda_{1}s}h(x))+1\right)\widetilde{m}(dx)\nonumber\\
&\le&2\int_{0}^{t}ds\int_{E}\alpha(x)\left(\log^{+}h(x)+1\right)\widetilde{m}(dx)\nonumber\\
&\le&2\|\alpha\|_{\infty}t\left(1+\int_{E}\log^{+}h(x)h(x)^{2}m(dx)\right)<{\infty}.\nonumber
\end{eqnarray}
Similarly by \eqref{llogl1}, we have
\begin{eqnarray}
&&\int_{E}\qq_{h\delta_{x}}\left[\sum_{s\in D^{m}_{t}}\left(\log^{+}(e^{\lambda_{1}s}I^{m}_{s}h(\xi_{s}))+1\right)\right]\widetilde{m}(dx)\nonumber\\
&=&\int_{E}\Pi^{h}_{x}\left[\int_{0}^{t}ds\left(\int_{(0,{\infty})}\log^{+}(e^{\lambda_{1}s}rh(\xi_{s}))r\pi(\xi_{s},dr)+1\right)\right]\widetilde{m}(dx)\nonumber\\
&=&t+\int_{0}^{t}ds\int_{E}\Pi^{h}_{x}\left[\int_{(0,{\infty})}r\log^{+}(e^{\lambda_{1}s}rh(\xi_{s}))\pi(\xi_{s},dr)\right]\widetilde{m}(dx)\nonumber\\
&=&t+\int_{0}^{t}ds\langle \int_{(0,{\infty})}r\log^{+}(e^{\lambda_{1}s}rh(\cdot))\pi(\cdot,dr),h^{2}\rangle\nonumber\\
&\le&t+t\langle \int_{(0,{\infty})}r\log^{+}(rh(\cdot))\pi(\cdot,dr),h^{2}\rangle\nonumber\\
&\le&t+t\langle \int_{(0,{\infty})}r\log^{*}(rh(\cdot))\pi(\cdot,dr),h^{2}\rangle<{\infty}.\nonumber
\end{eqnarray}
This completes the proof for \eqref{5.4}.
\qed

\medskip

\begin{lemma}\label{lem3.3}
For any $\mu\in\mc$ and $\phi\in\mathcal{B}^{+}_{b}(E)$,
\begin{equation}
\lim_{s\to{\infty}}\lim_{t\to{\infty}}e^{\lambda_{1}(s+t)}\langle \phi h,X^{*}_{s+t}+I^{*,t}_{s}\rangle=0\quad\mbox{ in }L^{1}(\pp_{\mu}).\nonumber
\end{equation}
\end{lemma}
\proof Recall that given $\mathcal{F}_{t}$, $(X^{*}_{s+t}+I^{*,t}_{s})_{s\ge 0}$ is equal in distribution to $((X^{*}_{s})_{s\ge 0};\pp_{X_{t}})$. Thus by Markov property,
\begin{eqnarray}
e^{\lambda_{1}(s+t)}\pp_{\mu}\left(\langle \phi h,X^{*}_{s+t}+I^{*,t}_{s}\rangle \right)
&=&e^{\lambda_{1}(s+t)}\pp_{\mu}\left[\pp_{\mu}\left(\langle \phi h,X^{*}_{s+t}+I^{*,t}_{s}\rangle|\mathcal{F}_{t}\right)\right]\nonumber\\
&=&e^{\lambda_{1}(s+t)}\pp_{\mu}\left[\pp_{X_{t}}\left(\langle \phi h,X^{*}_{s}\rangle\right)\right]\nonumber\\
&=&e^{\lambda_{1}(s+t)}\pp_{\mu}\left[\langle P^{\beta^{*}}_{s}\left(\phi h\right),X_{t}\rangle\right]\nonumber\\
&=&e^{\lambda_{1}t}\pp_{\mu}\left(\langle h\theta^{*}_{\phi h}(s,\cdot),X_{t}\rangle\right)\nonumber\\
&=&e^{\lambda_{1}t}\langle P^{\beta}_{t}\left(h\theta^{*}_{\phi h}(s,\cdot)\right),\mu\rangle\nonumber\\
&=&\langle hP^{h}_{t}\left(\theta^{*}_{\phi h}(s,\cdot)\right),\mu\rangle.\nonumber
\end{eqnarray}
Thus we have
\begin{eqnarray}
\limsup_{t\to{\infty}}e^{\lambda_{1}(s+t)}\pp_{\mu}\left(\langle \phi h,X^{*}_{s+t}+I^{*,t}_{s}\rangle \right)
&=& \limsup_{t\to{\infty}}\langle hP^{h}_{t}\left(\theta^{*}_{\phi h}(s,\cdot)\right),\mu\rangle\nonumber\\
&\le& \langle \limsup_{t\to{\infty}}P^{h}_{t}\left(\theta^{*}_{\phi h}(s,\cdot)\right),h\mu\rangle.\label{5.6}
\end{eqnarray}
We get by \eqref{p3} that for $t>1$ and $x\in E$
$$P^{h}_{t}\left(\theta^{*}_{\phi h}(s,\cdot)\right)(x)\le \int_{E}\theta^{*}_{\phi h}(s,x)\widetilde{m}(dx)+e^{-\lambda_{h}(t-1)}\widetilde{a}_{2}(x)^{1/2}\|\theta^{*}_{\phi h}(s,\cdot)\|_{L^{2}(E,\widetilde{m})}.$$
This implies that
\begin{equation}
\limsup_{t\to{\infty}}P^{h}_{t}\left(\theta^{*}_{\phi h}(s,\cdot)\right)(x)\le \int_{E}\theta^{*}_{\phi h}(s,x)\widetilde{m}(dx)\label{5.7}
\end{equation}
for all $x\in E$.
We know from Lemma \ref{lem1.1} that $\theta^{*}_{\phi h}(s,x)$ is bounded.
and $\lim_{s\to{\infty}}\theta^{*}_{\phi h}(s,x)=0$. Hence by \eqref{5.6}, \eqref{5.7} and the bounded convergence theorem,
\begin{eqnarray}
\limsup_{s\to{\infty}}\limsup_{t\to{\infty}}e^{\lambda_{1}(s+t)}\pp_{\mu}\left(\langle \phi h,X^{*}_{s+t}+I^{*,t}_{s}\rangle \right)
&\le&\limsup_{s\to{\infty}}\int_{E}\theta^{*}_{\phi h}(s,x)\widetilde{m}(dx)\langle h,\mu\rangle\nonumber\\
&\le&\int_{E}\limsup_{s\to{\infty}}\theta^{*}_{\phi h}(s,x)\widetilde{m}(dx)\langle h,\mu\rangle=0.\nonumber
\end{eqnarray}
\qed

\medskip

\begin{lemma}\label{lem3.4}
For any $\mu\in\mc$, $\phi\in\mathcal{B}^{+}_{b}(E)$ and $s>0$,
\begin{equation}
\lim_{t\to{\infty}}e^{\lambda_{1}(s+t)}\left[\sum_{i=1}^{N_{t}}\langle \phi h,I^{i,t}_{s}\rangle-\sum_{i=1}^{N_{t}}\pp_{\mu}\left(\langle\phi h,I^{i,t}_{s}\rangle|\mathcal{F}_{t}\rangle \right)\right]=0\quad\mbox{ in }L^{1}(\pp_{\mu}).\nonumber
\end{equation}
\end{lemma}

\proof Fix $\phi\in\mathcal{B}^{+}_{b}(E)$. For $s,t>0$, define
$$S_{s,t}:=e^{\lambda_{1}(s+t)}\sum_{i=1}^{N_{t}}\langle \phi h,I^{i,t}_{s}\rangle,$$
$$\widehat{S}_{s,t}:=e^{\lambda_{1}(s+t)}\sum_{i=1}^{N_{t}}\langle \phi h,I^{i,t}_{s}\rangle 1_{\{\langle \phi h,I^{i,t}_{s}\rangle<e^{-\lambda_{1}t}\}}.$$
Clearly by definition, $\widehat{S}_{s,t}\le S_{s,t}$. Recall that given $\mathcal{F}_{t}$, $\{I^{i,t}:i=1,\cdots,N_{t}\}$ are independent and $I^{i,t}$ is equal in distribution to $(I;\pp_{\cdot,\delta_{z_{i}(t)}})$. By this, we have
\begin{eqnarray}
\pp_{\mu}\left[\left(\widehat{S}_{s,t}-\pp_{\mu}\left(\widehat{S}_{s,t}|\mathcal{F}_{t}\right)\right)^{2}\right]
&=&\pp_{\mu}\left[\mbox{Var}\left(\widehat{S}_{s,t}|\mathcal{F}_{t}\right)\right]\nonumber\\
&\le&e^{2\lambda_{1}(s+t)}\pp_{\mu}\left[\sum_{i=1}^{N_{t}}
\pp_{\mu}\left(\langle \phi h,I^{i,t}_{s}\rangle^{2}1_{\{\langle \phi h,I^{i,t}_{s}\rangle<e^{-\lambda_{1}t}\}}|\mathcal{F}_{t}\right)\right]\nonumber\\
&=&e^{2\lambda_{1}(s+t)}\pp_{\mu}\left[\sum_{i=1}^{N_{t}}\pp_{\cdot,\delta_{z_{i}(t)}}\left(\langle \phi h,I_{s}\rangle^{2}1_{\{\langle \phi h,I_{s}\rangle<e^{-\lambda_{1}t}\}}\right)\right].\label{5.9}
\end{eqnarray}
Let $g_{s,t}(x):=w(x)\pp_{\cdot,\delta_{x}}\left(\langle \phi h,I_{s}\rangle^{2}1_{\{\langle \phi h,I_{s}\rangle<e^{-\lambda_{1}t}\}}\right)$ for $x\in E$.
Then by \eqref{meanz}
\begin{eqnarray}
\mbox{RHS of }\eqref{5.9}&=&e^{2\lambda_{1}(s+t)}\pp_{\mu}\left(\langle \frac{g_{s,t}}{w},Z_{t}\rangle\right)\nonumber\\
&=&e^{2\lambda_{1}(s+t)}\langle P^{\beta}_{t}(g_{s,t}),\mu\rangle\nonumber\\
&=&e^{\lambda_{1}(t+2s)}\langle P^{h}_{t}(\frac{g_{s,t}}{h}),h\mu\rangle.\label{5.10}
\end{eqnarray}
Clearly by \eqref{2.19},
$$g_{s,t}(x)\le\mathrm{e}^{-\lambda_{1}t}w(x)\pp_{\cdot,\delta_{x}}\left(\langle \phi h,I_{s}\rangle\right)=e^{-\lambda_{1}t}\pp_{\delta_{x}}\left(\langle \phi h,I_{s}\rangle\right)\le e^{-\lambda_{1}t}\pp_{\delta_{x}}\left(\langle \phi h,X_{s}\rangle\right)\le e^{-\lambda_{1}(s+t)}\|\phi\|_{\infty}
h(x)$$
for all $x\in E$. For $x\in E$ and $t>1$, by \eqref{p3},
\begin{eqnarray}
P^{h}_{t}\left(\frac{g_{s,t}}{h}\right)(x)&\le&\langle g_{s,t},h\rangle+e^{-\lambda_{h}(t-1)}\widetilde{a}_{2}(x)^{1/2}\|g_{s,t}\|_{L^{2}(E,m)}\nonumber\\
&\le &\langle g_{s,t},h\rangle+e^{-\lambda_{h}(t-1)-\lambda_{1}(s+t)}\|\phi\|_{\infty}\widetilde{a}_{2}(x)^{1/2}.\label{5.11}
\end{eqnarray}
Hence by \eqref{5.9}-\eqref{5.11}, we get
\begin{equation}
\pp_{\mu}\left[\left(\widehat{S}_{s,t}-\pp_{\mu}\left(\widehat{S}_{s,t}|\mathcal{F}_{t}\right)\right)^{2}\right]
\le e^{\lambda_{1}(t+2s)}\langle h,\mu\rangle \langle g_{s,t},h\rangle+e^{-\lambda_{h}(t-1)+\lambda_{1}s}\|\phi\|_{\infty}\langle \widetilde{a}^{1/2}_{2}h,\mu\rangle.\label{5.12}
\end{equation}
Here $\langle \widetilde{a}^{1/2}_{2}h,\mu\rangle=\int_{E}p^{h}(2,x,x)^{1/2}h(x)\mu(dx)<{\infty}$
since both $x\mapsto p^{h}(2,x,x)$ and $h$ are bounded on the compact support of $\mu$
due to the continuity of $x\mapsto p^{h}(2,x,x)$ and $h$.
Note that by Fubini's theorem
 \begin{eqnarray}
&&\int_{1}^{{\infty}}e^{\lambda_{1}t}\langle g_{s,t},h\rangle dt\nonumber\\
&=&\int_{1}^{{\infty}}e^{\lambda_{1}t}dt\int_{E}h(x)w(x)m(dx)\int_{[0,e^{-\lambda_{1}t})}y^{2}\pp_{\cdot,\delta_{x}}\left(\langle \phi h,I_{s}\rangle\in dy\right)\nonumber\\
&=&\int_{E}h(x)w(x)m(dx)\int_{e^{-\lambda_{1}}}^{{\infty}}\frac{1}{-\lambda_{1}s^{2}}\,ds\int_{[0,s)}y^{2}\pp_{\cdot,\delta_{x}}\left(\langle \phi h,I_{s}\rangle\in dy\right)\nonumber\\
&=&\int_{E}h(x)w(x)m(dx)\int_{[0,{\infty})}y^{2}\pp_{\cdot,\delta_{x}}\left(\langle \phi h,I_{s}\rangle\in dy\right)\int_{e^{-\lambda_{1}}\vee y}^{{\infty}}\frac{1}{-\lambda_{1}s^{2}}\,ds\nonumber\\
&\le&\frac{1}{-\lambda_{1}}\int_{E}h(x)w(x)m(dx)\int_{[0,{\infty})}y\pp_{\cdot,\delta_{x}}\left(\langle \phi h,I_{s}\rangle\in dy\right)\nonumber\\
&=&\frac{1}{-\lambda_{1}}\int_{E}w(x)\pp_{\cdot,\delta_{x}}\left(\langle \phi h,I_{s}\rangle\right)h(x)m(dx).\nonumber
\end{eqnarray}
By \eqref{2.19} the integral in the right hand side equals
\begin{eqnarray}
\int_{E}\pp_{\delta_{x}}\left(\langle \phi h,I_{s}\rangle\right)h(x)m(dx)
&\le&\|\phi\|_{\infty}\int_{E}\pp_{\delta_{x}}\left(\langle  h,X_{s}\rangle\right)h(x)m(dx)\nonumber\\
&=&\|\phi\|_{\infty}e^{-\lambda_{1}s}\int_{E}h(x)^{2}m(dx)<{\infty}.\label{5.22}
\end{eqnarray}
Hence we get $\limsup_{t\to{\infty}}e^{\lambda_{1}t}\langle g_{s,t},h\rangle =0$. This together with \eqref{5.12} yields that
$$\limsup_{t\to{\infty}}\pp_{\mu}\left[\left(\widehat{S}_{s,t}-\pp_{\mu}\left(\widehat{S}_{s,t}|\mathcal{F}_{t}\right)\right)^{2}\right]=0,$$ and consequently
\begin{equation}
\lim_{t\to{\infty}}\left(\widehat{S}_{s,t}-\pp_{\mu}\left(\widehat{S}_{s,t}|\mathcal{F}_{t}\right)\right)=0\quad\mbox{ in }L^{1}(\pp_{\mu}).\label{5.13}
\end{equation}
Recall that under $\pp_{\delta_{z}}$, $Z_{0}$ is a Poisson random measure with intensity $w(x)\delta_{x}(dy)$. Thus we have
\begin{eqnarray}
\pp_{\delta_{x}}\left[\langle \phi h,I_{s}\rangle 1_{\{\langle \phi h,I_{s}\rangle\ge \mathrm{e}^{-\lambda_{1}t}\}}\right]
&=&\pp_{\delta_{x}}\left[\sum_{i=1}^{N_{0}}\langle \phi h,I^{i,0}_{s}\rangle 1_{\{\sum_{i=1}^{N_{0}}\langle \phi h,I^{i,0}_{s}\rangle\ge \mathrm{e}^{-\lambda_{1}t}\}}\right]\nonumber\\
&\ge &\pp_{\delta_{x}}\left[\sum_{i=1}^{N_{0}}\langle \phi h,I^{i,0}_{s}\rangle 1_{\{\langle \phi h,I^{i,0}_{s}\rangle\ge \mathrm{e}^{-\lambda_{1}t}\}}\right]\nonumber\\
&=&\pp_{\delta_{x}}\left[\sum_{i=1}^{N_{0}}\pp_{\cdot,\delta_{z_{i}(0)}}\left[\langle \phi h,I_{s}\rangle 1_{\{\langle \phi h,I_{s}\rangle\ge \mathrm{e}^{-\lambda_{1}t}\}}\right]\right]\nonumber\\
&=&w(x)\pp_{\cdot,\delta_{x}}\left[\langle \phi h,I_{s}\rangle 1_{\{\langle \phi h,I_{s}\rangle\ge \mathrm{e}^{-\lambda_{1}t}\}}\right].\nonumber
\end{eqnarray}
Using this and Markov property, we have
\begin{eqnarray}
\pp_{\mu}\left(S_{s,t}-\widehat{S}_{s,t}\right)&=&\pp_{\mu}\left[\pp_{\mu}\left(S_{s,t}|\mathcal{F}_{t}\right)-
\pp_{\mu}\left(\widehat{S}_{s,t}|\mathcal{F}_{t}\right)\right]\nonumber\\
&=&e^{\lambda_{1}(s+t)}\pp_{\mu}\left[\sum_{i=1}^{N_{t}}\pp_{\cdot,\delta_{z_{i}(t)}}\left(\langle \phi h,I_{s}\rangle 1_{\{\langle \phi h,I_{s}\rangle\ge e^{-\lambda_{1}t}\}}\right)\right]\nonumber\\
&\le&e^{\lambda_{1}(s+t)}\pp_{\mu}\left[\sum_{i=1}^{N_{t}}w(z_{i}(t))^{-1}\pp_{\delta_{z_{i}(t)}}\left(\langle \phi h,I_{s}\rangle 1_{\{\langle \phi h,I_{s}\rangle\ge e^{-\lambda_{1}t}\}}\right)\right]\nonumber\\
&=&e^{\lambda_{1}(s+t)}\pp_{\mu}\left(\langle \frac{f_{s,t}}{w},Z_{t}\rangle\right)\nonumber\\
&=&e^{\lambda_{1}(s+t)}\langle P^{\beta}_{t}f_{s,t},\mu\rangle\nonumber\\
&=&e^{\lambda_{1}s}\langle P^{h}_{t}\left(\frac{f_{s,t}}{h}\right),h\mu\rangle,\label{5.14}
\end{eqnarray}
where $f_{s,t}(x):=\pp_{\delta_{x}}\left(\langle \phi h,I_{s}\rangle 1_{\{\langle \phi h,I_{s}\rangle\ge e^{-\lambda_{1}t}\}}\right)$. We observe that $f_{s,t}(x)\le \pp_{\delta_{x}}\left(\langle \phi h,X_{s}\rangle\right)\le \|\phi\|_{\infty}e^{-\lambda_{1}s}h(x)$ for every $x\in E$. By \eqref{p3} again, we get for $x\in E$ and $t>1$,
\begin{eqnarray}
P^{h}_{t}\left(\frac{f_{s,t}}{h}\right)(x)&\le&\langle f_{s,t},h\rangle+e^{-\lambda_{h}(t-1)}\widetilde{a}_{2}(x)^{1/2}\|f_{s,t}\|_{L^{2}(E,m)}\nonumber\\
&\le &\langle f_{s,t},h\rangle+e^{-\lambda_{h}(t-1)-\lambda_{1}s}\|\phi\|_{\infty}\widetilde{a}_{2}(x)^{1/2}.\nonumber
\end{eqnarray}
This together with \eqref{5.14} yields that for $t>1$,
\begin{eqnarray}
\pp_{\mu}\left(S_{s,t}-\widehat{S}_{s,t}\right)&=&\pp_{\mu}\left[\pp_{\mu}\left(S_{s,t}|\mathcal{F}_{t}\right)-
\pp_{\mu}\left(\widehat{S}_{s,t}|\mathcal{F}_{t}\right)\right]\nonumber\\
&\le&e^{\lambda_{1}s}\langle h,\mu\rangle \langle f_{s,t},h\rangle+e^{-\lambda_{h}(t-1)}\|\phi\|_{\infty}\langle \widetilde{a}_{2}^{1/2}h,\mu\rangle.\label{5.15}
\end{eqnarray}
It follows from Fubini's theorem and Lemma \ref{lem3.1} that
\begin{eqnarray}
\int_{0}^{{\infty}}\langle f_{s,t},h\rangle dt
&=&\int_{0}^{{\infty}}dt\int_{E}h(x)m(dx)\int_{[e^{-\lambda_{1}t},{\infty})}y\pp_{\delta_{x}}\left(\langle \phi h,I_{s}\rangle \in dy\right)\nonumber\\
&=&\int_{E}h(x)m(dx)\int_{1}^{{\infty}}\frac{1}{-\lambda_{1}s}\,ds\int_{[s,{\infty})}y\pp_{\delta_{x}}\left(\langle \phi h,I_{s}\rangle \in dy\right)\nonumber\\
&=&\frac{1}{-\lambda_{1}}\int_{E}h(x)m(dx)\int_{[1,{\infty})}y\pp_{\delta_{x}}\left(\langle \phi h,I_{s}\rangle \in dy\right)\int_{1}^{y}\frac{1}{s}\,ds\nonumber\\
&=&\frac{1}{-\lambda_{1}}\int_{E}h(x)m(dx)\int_{[1,{\infty})}y\log y\pp_{\delta_{x}}\left(\langle \phi h,I_{s}\rangle \in dy\right)\nonumber\\
&=&\frac{1}{-\lambda_{1}}\int_{E}\pp_{\delta_{x}}\left(\langle \phi h,I_{s}\rangle\log^{+}\langle \phi h,I_{s}\rangle\right)h(x)m(dx)\nonumber\\
&\le&\frac{1}{-\lambda_{1}}\int_{E}\pp_{\delta_{x}}\left(\langle \phi h,X_{s}\rangle\log^{+}\langle \phi h,X_{s}\rangle\right)h(x)m(dx)<{\infty}.\label{5.23}
\end{eqnarray}
Immediately, $\lim_{t\to{\infty}}\langle f_{s,t},h\rangle =0$. Hence by \eqref{5.15} we get $$\lim_{t\to{\infty}}\pp_{\mu}\left(S_{s,t}-\widehat{S}_{s,t}\right)=\lim_{t\to{\infty}}\pp_{\mu}\left[\pp_{\mu}\left(S_{s,t}|\mathcal{F}_{t}\right)-
\pp_{\mu}\left(\widehat{S}_{s,t}|\mathcal{F}_{t}\right)\right]=0.$$ In other words,
\begin{equation}
\lim_{t\to{\infty}}\left(S_{s,t}-\widehat{S}_{s,t}\right)=\lim_{t\to{\infty}}\left(\pp_{\mu}\left(S_{s,t}|\mathcal{F}_{t}\right)-
\pp_{\mu}\left(\widehat{S}_{s,t}|\mathcal{F}_{t}\right)\right)=0\quad\mbox{ in }L^{1}(\pp_{\mu}).\label{5.16}
\end{equation}
Combining \eqref{5.13} and \eqref{5.16}, we get $\lim_{t\to{\infty}}\left(S_{s,t}-\pp_{\mu}\left(S_{s,t}|\mathcal{F}_{t}\right)\right)=0$ in $L^{1}(\pp_{\mu})$. This completes the proof.\qed

\medskip

\begin{lemma}\label{lem3.5}
For any $\mu\in\mc$ and $\phi\in\mathcal{B}^{+}_{b}(E)$,
\begin{equation}
\lim_{s\to{\infty}}\lim_{t\to{\infty}}e^{\lambda_{1}(s+t)}\sum_{i=1}^{N_{t}}\pp_{\mu}\left(\langle \phi h,I^{i,t}_{s}\rangle |\mathcal{F}_{t}\right)=\langle \phi h,h\rangle W^{h}_{\infty}(X)\quad\mbox{ in }L^{1}(\pp_{\mu}).\label{5.17}
\end{equation}
\end{lemma}

\proof In view of Theorem \ref{them1} and Theorem \ref{them2}, under Assumptions 1-3 and
Assumption 4(i),
\eqref{5.17} is equivalent to
\begin{equation}\label{5.18}
\lim_{s\to{\infty}}\lim_{t\to{\infty}}\left[e^{\lambda_{1}(s+t)}\sum_{i=1}^{N_{t}}\pp_{\mu}\left(\langle \phi h,I^{i,t}_{s}\rangle |\mathcal{F}_{t}\right)-\langle \phi h,h\rangle W^{h/w}_{t}(Z)\right]=0\quad\mbox{ in }L^{1}(\pp_{\mu}).
\end{equation}
Fix $\mu\in\mc$ and $\phi\in\mathcal{B}^{+}_{b}(E)$.
Note that
for $x\in E$ and $s>0$,
\begin{eqnarray}
\pp_{\delta_{x}}\left(\langle \phi h,I_{s}\rangle\right)&=&\pp_{\delta_{x}}\left(\langle \phi h,X_{s}\rangle\right)-\pp_{\delta_{x}}\left(\langle \phi h,X^{*}_{s}\rangle\right)\nonumber\\
&=&P^{\beta}_{s}(\phi h)(x)-P^{\beta^{*}}_{s}(\phi h)(x).\nonumber
\end{eqnarray}
By this and Markov property, we have
\begin{eqnarray}
e^{\lambda_{1}(s+t)}\sum_{i=1}^{N_{t}}\pp_{\mu}\left(\langle \phi h,I^{i,t}_{s}\rangle |\mathcal{F}_{t}\right)
&=&e^{\lambda_{1}(s+t)}\sum_{i=1}^{N_{t}}\pp_{\cdot,\delta_{z_{i}(t)}}\left(\langle \phi h,I_{s}\rangle\right)\nonumber\\
&=&e^{\lambda_{1}(s+t)}\sum_{i=1}^{N_{t}}w(z_{i}(t))^{-1}\pp_{\delta_{z_{i}(t)}}\left(\langle \phi h,I_{s}\rangle\right)\nonumber\\
&=&e^{\lambda_{1}(s+t)}\langle \frac{1}{w}\left(P^{\beta}_{s}(\phi h)-P^{\beta^{*}}_{s}(\phi h)\right),Z_{t}\rangle\nonumber\\
&=&e^{\lambda_{1}t}\langle \frac{h}{w}\left(P^{h}_{s}\phi-\theta^{*}_{\phi h}(s,\cdot)\right),Z_{t}\rangle.\label{5.19}
\end{eqnarray}
Let $g(s,x):=\left|P^{h}_{s}\phi(x)-\theta^{*}_{\phi h}(s,x)-\langle \phi h,h\rangle\right|$ for $s>0$ and $x\in E$.
Clearly $g$ is bounded from above by $3\|\phi\|_{\infty}$. By \eqref{5.19},
\begin{eqnarray}
I(s,t)&:=&\pp_{\mu}\left[|e^{\lambda_{1}(s+t)}\sum_{i=1}^{N_{t}}\pp_{\mu}\left(\langle \phi h,I^{i,t}_{s}\rangle |\mathcal{F}_{t}\right)
-\langle \phi h,h\rangle W^{h/w}_{t}(Z)|\right]\nonumber\\
&=&e^{\lambda_{1}t}\pp_{\mu}\left(\left|\langle \frac{h}{w}\left(P^{h}_{s}\phi-\theta^{*}_{\phi h}(s)-\langle \phi h,h\rangle\right),Z_{t}\rangle\right|\right)\nonumber\\
&\le&e^{\lambda_{1}t}\pp_{\mu}\left(\langle \frac{h}{w}g(s,\cdot),Z_{t}\rangle\right)\nonumber\\
&=&e^{\lambda_{1}t}\langle P^{\beta}_{t}(hg(s,\cdot)),\mu\rangle\nonumber\\
&=&\langle P^{h}_{t}(g(s,\cdot)),h\mu\rangle.\label{5.20}
\end{eqnarray}
Since by \eqref{p3}, for $t>1$ and $x\in E$,
\begin{eqnarray}
P^{h}_{t}(g(s,\cdot))(x)&\le&\langle g(s,\cdot),h^{2}\rangle+e^{-\lambda_{h}(t-1)}\widetilde{a}_{2}(x)^{1/2}\|g(s,\cdot)\|_{L^{2}(E,\widetilde{m})}\nonumber\\
&\le& \langle g(s,\cdot),h^{2}\rangle+3e^{-\lambda_{h}(t-1)}\|\phi\|_{\infty}\widetilde{a}_{2}(x)^{1/2},\nonumber
\end{eqnarray}
it follows by \eqref{5.20} that $I(s,t)\le \langle g(s,\cdot),h^{2}\rangle\langle h,\mu\rangle+3e^{-\lambda_{h}(t-1)}\|\phi\|_{\infty}\langle \widetilde{a}_{2}^{1/2}h,\mu\rangle$. Hence we get
\begin{equation}\label{5.21}
\limsup_{t\to{\infty}}I(s,t)\le \langle g(s,\cdot),h^{2}\rangle\langle h,\mu\rangle.
\end{equation}
Furthermore, by \eqref{p3} and Lemma \ref{lem1.1}, for $s>1$ and $x\in E$,
\begin{eqnarray}
g(s,x)&\le&\left|P^{h}_{s}\phi-\langle \phi h,h\rangle\right|+\theta^{*}_{\phi h}(s,x)\nonumber\\
&\le &e^{-\lambda_{h}(t-1)}\widetilde{a}_{2}(x)^{1/2}\|\phi\|_{L^{2}(E,\widetilde{m})}+\theta^{*}_{\phi h}(s,x)\to 0 \mbox{ as }s\to{\infty}.\nonumber
\end{eqnarray}
Thus $\lim_{s\to{\infty}}\langle g(s,\cdot),h^{2}\rangle=0$ by the bounded convergence theorem.
This and \eqref{5.21} yields that $\lim_{s\to{\infty}}\lim_{t\to {\infty}}I(s,t)=0$. Hence we prove \eqref{5.18}.\qed

\medskip

\noindent\textit{Proof of Theorem \ref{them3}:} We know by \cite[Lemma 3.2]{EKW} that
if for any
$\mu\in\mc$ and $f\in\mathcal{B}^{+}(E)$ with $f/h$ bounded,
\begin{equation}
\lim_{t\to{\infty}}e^{\lambda_{1}t}\langle f,X_{t}\rangle =\langle f,h\rangle W^{h}_{\infty}(X)\quad\mbox{ in }L^{1}(\p_{\mu}),\label{5.5}
\end{equation}
then the convergence in \eqref{5.5} holds for any $\mu\in \mh$.
Henceforth we assume $\mu\in\mc$. We consider the skeleton space for convenience. For an arbitrary $f\in \mathcal{B}^{+}(E)$ with $f/h$ bounded, let $\phi(x):=f(x)/h(x)\in\mathcal{B}^{+}_{b}(E)$. For any $s,t>0$, under $\pp_{\mu}$,
\begin{eqnarray}
&&e^{\lambda_{1}(s+t)}\langle f,X_{s+t}\rangle-\langle f,h\rangle W^{h}_{\infty}(X)\nonumber\\
&=&e^{\lambda_{1}(s+t)}\langle \phi h,X^{*}_{s+t}+I^{*,t}_{s}\rangle+\left[e^{\lambda_{1}(s+t)}\sum_{i=1}^{N_{t}}\langle \phi h,I^{i,t}_{s}\rangle-
e^{\lambda_{1}(s+t)}\sum_{i=1}^{N_{t}}\pp_{\mu}\left(\langle \phi h,I^{i,t}_{s}\rangle|\mathcal{F}_{t}\right)\right]\nonumber\\
&&+\left[e^{\lambda_{1}(s+t)}\sum_{i=1}^{N_{t}}\pp_{\mu}\left(\langle \phi h,I^{i,t}_{s}\rangle|\mathcal{F}_{t}\right)-\langle \phi h,h\rangle W^{h}_{\infty}(X)\right].\nonumber
\end{eqnarray}
Thus by letting $t\to {\infty}$ and $s\to {\infty}$, \eqref{5.5} follows immediately from Lemma \ref{lem3.3}, Lemma \ref{lem3.4} and Lemma \ref{lem3.5}.\qed

\medskip

\begin{prop}\label{P:6.5}
Under Assumptions 1-3 and Assumption 4(i),
Assumption 5' is equivalent to Assumption 5.
\end{prop}

\proof
Clearly by the symmetry of $P^B_t$, Assumption 5 implies Assumption 5'. So we
only need to show that Assumption 5' is sufficient for Assumption 5. For $f\in\mathcal{B}^{+}(E)$ with $fw/h$ bounded, let $\phi:=fw/h\in\mathcal{B}^{+}_{b}(E)$. Then
for any $\mu\in\mc$, $\sigma>0$ and $m,n\in\mathbb{N}$,
\begin{eqnarray}
&&e^{\lambda_{1}(m+n)\sigma}\langle f,Z_{(m+n)\sigma}\rangle-\langle f,wh\rangle W^{h/w}_{\infty}(Z)\nonumber\\
&=&e^{\lambda_{1}(m+n)\sigma}\langle \frac{h}{w}\phi,Z_{(m+n)\sigma}\rangle-\langle \phi,h^{2}\rangle W^{h/w}_{\infty}(Z)\nonumber\\
&=&\left[e^{\lambda_{1}(m+n)\sigma}\langle \frac{h}{w}\phi,Z_{(m+n)\sigma}\rangle-e^{\lambda_{1}(m+n)\sigma}\pp_{\mu}\left(\langle \frac{h}{w}\phi,Z_{(m+n)\sigma}\rangle|\mathcal{F}_{n\sigma}\right)\right]\nonumber\\
&&+\left[e^{\lambda_{1}(m+n)\sigma}\pp_{\mu}\left(\langle \frac{h}{w}\phi,Z_{(m+n)\sigma}\rangle|\mathcal{F}_{n\sigma}\right)-\langle \phi,h^{2}\rangle W^{h/w}_{\infty}(Z)\right]\nonumber\\
&=:&I(m,n)+II(m,n).\nonumber
\end{eqnarray}
Note that by the Markov property of $Z$,
$e^{\lambda_{1}(m+n)\sigma}\pp_{\mu}\left(\langle \frac{h}{w}\phi,Z_{(m+n)\sigma}\rangle|\mathcal{F}_{n\sigma}\right)=e^{\lambda_{1}(m+n)\sigma}\pp_{\cdot,Z_{n\sigma}}
\left(\langle \frac{h}{w}\phi,Z_{m\sigma}\rangle\right)=e^{\lambda_{1}(m+n)\sigma}\langle \frac{1}{w}P^{\beta}_{m\sigma}(h\phi),Z_{n\sigma}\rangle=e^{\lambda_{1}n\sigma}\langle \frac{h}{w}P^{h}_{m\sigma}(\phi),Z_{n\sigma}\rangle$.
Thus Assumption 5' implies that $II(m,n)$ converges to $0$ almost surely for some $m\in \mathbb{N}$. Hence the proof is finished if we can show that
for any $m\in\mathbb{N}$,
\begin{equation}
\lim_{n\to{\infty}}I(m,n)=0\quad\pp_{\mu}
\mbox{-a.s.}
\label{6.27}
\end{equation}
It follows from Theorem \ref{them1}, Lemma \ref{lem3.1} and the fact that $\log^{+}ab\le \log^{+}a+\log^{+}b$ for all $a,b\ge 0$ that under Assumptions 1-3 and Assumption 4(i),
\begin{equation}
\int_{E}\pp_{\delta_{x}}\left(\langle \frac{h}{w}\phi,Z_{t}\rangle
\log^{+}\langle \frac{h}{w}\phi,Z_{t}\rangle
\right)h(x)m(dx)<{\infty}\label{6.28}
\end{equation}
for all $t>0$. Recall that under $\pp_{\delta_{x}}$, $N_{0}$ is a random variable with mean $w(x)$.
By this we have
\begin{eqnarray}
\pp_{\delta_{x}}\left(\langle \frac{h}{w}\phi,Z_{t}\rangle
\log^{+}\langle \frac{h}{w}\phi,Z_{t}\rangle
\right)
&=&\pp_{\delta_{x}}\left(\sum_{i=1}^{N_{0}}\langle \frac{h}{w}\phi,Z^{i,0}_{t}\rangle
\log^{+}\langle \frac{h}{w}\phi,\sum_{i=1}^{N_{0}}Z^{i,0}_{t}\rangle\right)\nonumber\\
&\ge&\pp_{\delta_{x}}\left(\sum_{i=1}^{N_{0}}\langle \frac{h}{w}\phi,Z^{i,0}_{t}\rangle
\log^{+}\langle \frac{h}{w}\phi,Z^{i,0}_{t}\rangle\right)\nonumber\\
&=&\pp_{\delta_{x}}\left(N_{0}\pp_{\cdot,\delta_{x}}\left(\langle \frac{h}{w}\phi,Z_{t}\rangle\log^{+}\langle \frac{h}{w}\phi,Z_{t}\rangle\right)\right)\nonumber\\
&=&w(x)\pp_{\cdot,\delta_{x}}\left(\langle \frac{h}{w}\phi,Z_{t}\rangle\log^{+}\langle \frac{h}{w}\phi,Z_{t}\rangle\right).\nonumber
\end{eqnarray}
Thus \eqref{6.28} implies that
\begin{equation}
\int_{E}\pp_{\cdot,\delta_{x}}\left(\langle \frac{h}{w}\phi,Z_{t}\rangle
\log^{+}\langle \frac{h}{w}\phi,Z_{t}\rangle
\right)w(x)h(x)m(dx)<{\infty}\label{6.29}
\end{equation}
Recall that the skeleton $Z$ is a supercriticle branching Markov process with spatial motion $\xi^{w}$ which is symmetric with respect to the measure $w^{2}m$.  \eqref{6.27} follows from \eqref{6.29} (instead of \cite[(5.1)]{CRY}) and Borel-Cantelli lemma
in the same way as \cite[Lemma 5.2]{CRY}. We omit the details here.\qed

\section{Strong law of large numbers}\label{sec7}

To prove the strong law of large numbers, we know by \cite{CRSZ} that it suffices to consider fixed test functions.

\begin{lemma}\label{lem4.6}{\rm [Chen et al. \cite{CRSZ}]}
 Suppose Assumption 2 holds.
 If for any
$\mu\in \mh$ and $\phi\in C^{+}_{0}(E)$,
\begin{equation}
\lim_{t\to{\infty}}e^{\lambda_{1}t}\langle \phi h,X_{t}\rangle=\langle \phi h,h\rangle W^{h}_{\infty}(X)\quad\p_{\mu}
\mbox{-a.s.}
\label{6.1}
\end{equation}
then there exists $\Omega_{0}\subset \Omega$ of $\p_{\mu}$-full probability for every $\mu\in\mh$ such that on $\Omega_{0}$, for every $m$-almost everywhere continuous nonnegative measurable function $f$ with $f/h$ bounded,
\begin{equation}
\lim_{t\to{\infty}}e^{\lambda_{1}t}\langle f,X_{t}\rangle=\langle f,h\rangle W^{h}_{\infty}(X).\nonumber
\end{equation}
\end{lemma}

\begin{lemma}\label{lem4.7}
 Suppose Assumption 2 holds.
If \eqref{6.1} holds
for any
$\mu\in \mc$ and $\phi\in C^{+}_{0}(E)$, it holds
for any
$\mu\in \mh.$
\end{lemma}
\proof The main idea of this proof is borrowed from \cite[Lemma 3.2(i)]{EKW}. Fix $\mu\in\mh$.
Take a sequence of sets $\{B_{k}:k\in\mathbb{N}\}$ such that $B_{0}=\emptyset$, $B_{k}\Subset E$, $B_{k}\subset B_{k+1}$ and $\bigcup_{k=1}^{{\infty}}B_{k}=E$.
Let $\widehat{B}_{k}:=B_{k}\setminus B_{k-1}$. On a suitable probability space with probability measure $\p_{\mu}$, let $\{X^{\widehat{B}_{k}}:k\in\mathbb{N}\}$ be independent $(P_{t},\psi_{\beta})$-superprocesses where $X^{\widehat{B}_{k}}$ is started in $1_{\widehat{B}_{k}}\mu$. Then by the branching property, $X^{B_{k}}:=\sum_{i=1}^{k}X^{\widehat{B}_{i}}$, $X^{E\setminus B_{k}}:=\sum_{i=k+1}^{{\infty}}X^{\widehat{B}_{i}}$ and $X:=X^{B_{k}}+X^{E\setminus B_{k}}$ are $(P_{t},\psi_{\beta})$-superprocesses started in $1_{B_{k}}\mu$, $1_{E\setminus B_{k}}\mu$ and $\mu$ respectively.
In particular, $W^{h}_{t}(X)=W^{h}_{t}(X^{B_{k}})+W^{h}_{t}(X^{E\setminus B_{k}})$ and the martingale limits $W^{h}_{\infty}(X^{E\setminus B_{k}})=\lim_{t\to{\infty}}W^{h}_{t}(X^{E\setminus B_{k}})$ is non-increasing in $k$.
Since $\mu\in \mh$, by Fatou's lemma
\begin{equation}
\p_{\mu}\left(W^{h}_{\infty}(X^{E\setminus B_{k}})\right)\le \liminf_{t\to{\infty}}\p_{\mu}\left(W^{h}_{t}(X^{E\setminus B_{k}})\right)=\langle h, 1_{E\setminus B_{k}}\mu\rangle \to 0\quad\mbox{ as }k\to{\infty}.\nonumber
\end{equation}
Thus $\lim_{k\to{\infty}}W^{h}_{\infty}(X^{E\setminus B_{k}})=0$ in $L^{1}(\p_{\mu})$ and hence $\p_{\mu}$-a.s. by the monotonicity.

For any $\phi\in C^{+}_{0}(E)$,
\begin{eqnarray}
&&\left|e^{\lambda_{1}t}\langle \phi h,X_{t}\rangle-\langle \phi h,h\rangle W^{h}_{\infty}(X)\right|\nonumber\\
&\le&\left|e^{\lambda_{1}t}\langle \phi h,X^{B_{k}}_{t}\rangle-\langle \phi h,h\rangle W^{h}_{\infty}(X^{B_{k}})\right|
+e^{\lambda_{1}t}\langle \phi h,X^{E\setminus B_{k}}_{t}\rangle+\langle \phi h,h\rangle W^{h}_{\infty}(X^{E\setminus B_{k}})\nonumber\\
&\le&\left|e^{\lambda_{1}t}\langle \phi h,X^{B_{k}}_{t}\rangle-\langle \phi h,h\rangle W^{h}_{\infty}(X^{B_{k}})\right|
+\|\phi\|_{\infty}W^{h}_{t}(X^{E\setminus B_{k}})+\langle \phi h,h\rangle W^{h}_{\infty}(X^{E\setminus B_{k}}).\nonumber
\end{eqnarray}
 Since $1_{B_{k}}\mu\in \mc$, the first term converges to $0$ as $t\to +\infty$ by our assumption.
 Therefore we get \eqref{6.1} for $\mu\in\mh$ by first letting $t\to{\infty}$ and then $k\to{\infty}$.\qed

\subsection{Strong law of large numbers along lattice times}

\begin{lemma}\label{lem4.1}
Suppose Assumptions 1-3
and Assumption 4(i) hold.
Then for any
$\mu\in\mc$, $\phi\in\mathcal{B}^{+}_{b}(E)$, $\sigma>0$ and $m\in\mathbb{N}$,
\begin{equation}
\lim_{n\to{\infty}}e^{\lambda_{1}(m+n)\sigma}\left[\sum_{i=1}^{N_{n\sigma}}\langle \phi h,I^{i,n\sigma}_{m\sigma}\rangle-
\sum_{i=1}^{N_{n\sigma}}\pp_{\mu}\left(\langle \phi h,I^{i,n\sigma}_{m\sigma}\rangle|\mathcal{F}_{n\sigma}\right)\right]=0\quad \pp_{\mu}\mbox{-a.s.}\nonumber
\end{equation}
\end{lemma}

\proof In this proof, we adopt the same notation as defined in the proof of Lemma \ref{lem3.4}. For $m,n\in \mathbb{N}$, let
$$S_{m\sigma,n\sigma}:=e^{\lambda_{1}(m+n)\sigma}\sum_{i=1}^{N_{n\sigma}}\langle \phi h,I^{i,n\sigma}_{m\sigma}\rangle,$$
$$\widehat{S}_{m\sigma,n\sigma}:=e^{\lambda_{1}(m+n)\sigma}\sum_{i=1}^{N_{n\sigma}}\langle \phi h,I^{i,n\sigma}_{m\sigma}\rangle 1_{\{\langle \phi h,I^{i,n\sigma}_{m\sigma}\rangle<e^{-\lambda_{1}n\sigma}\}}.$$
We have proved in \eqref{5.12} that for $n\sigma>1$,
\begin{eqnarray}
\pp_{\mu}\left[\left(\widehat{S}_{m\sigma,n\sigma}-\pp_{\mu}\left(\widehat{S}_{m\sigma,n\sigma}|\mathcal{F}_{n\sigma}\right)\right)^{2}\right]
&\le& e^{\lambda_{1}(n\sigma+2m\sigma)}\langle g_{m\sigma,n\sigma},h\rangle \langle h,\mu\rangle\nonumber\\
&&+e^{-\lambda_{h}(n\sigma-1)+\lambda_{1}m\sigma}\|\phi\|_{\infty}\langle \widetilde{a}^{1/2}_{2}h,\mu\rangle\label{6.2}
\end{eqnarray}
where $g_{m\sigma,n\sigma}(x):=w(x)\pp_{\cdot,\delta_{x}}\left(\langle \phi h,I_{m\sigma}\rangle^{2}1_{\{\langle \phi h,I_{m\sigma}\rangle<e^{-\lambda_{1}n\sigma}\}}\right)$ for $x\in E$. Note that by \eqref{5.22}
\begin{eqnarray}
\sum_{n=1}^{{\infty}}e^{\lambda_{1}n\sigma}\langle g_{m\sigma,n\sigma},h\rangle
&=&\sum_{n=1}^{{\infty}}e^{\lambda_{1}n\sigma}\int_{E}h(x)w(x)\pp_{\cdot,\delta_{x}}\left(\langle \phi h,I_{m\sigma}\rangle^{2}1_{\{\langle \phi h,I_{m\sigma}\rangle<e^{-\lambda_{1}n\sigma}\}}\right)m(dx)\nonumber\\
&\le&\int_{1}^{{\infty}}e^{\lambda_{1}\sigma(s-1)}\int_{E}h(x)w(x)\pp_{\cdot,\delta_{x}}\left(\langle \phi h,I_{m\sigma}\rangle^{2}1_{\{\langle \phi h,I_{m\sigma}\rangle<e^{-\lambda_{1}s\sigma}\}}\right)m(dx)\nonumber\\
&=&e^{-\lambda_{1}\sigma}\int_{1}^{{\infty}}e^{\lambda_{1}s\sigma}\langle g_{m\sigma,s\sigma},h\rangle ds<{\infty}.\nonumber
\end{eqnarray}
This together with \eqref{6.2} yields that $\sum_{n=1}^{{\infty}}\pp_{\mu}\left[\left(\widehat{S}_{m\sigma,n\sigma}-\pp_{\mu}\left(\widehat{S}_{m\sigma,n\sigma}|\mathcal{F}_{n\sigma}\right)\right)^{2}\right]
<{\infty}$. Thus by Borel-Cantelli lemma,
\begin{equation}\label{6.3}
\lim_{n\to{\infty}}\left(\widehat{S}_{m\sigma,n\sigma}-\pp_{\mu}\left(\widehat{S}_{m\sigma,n\sigma}|\mathcal{F}_{n\sigma}\right)\right)=0\quad\pp_{\mu}\mbox{-a.s.}
\end{equation}
Similarly, we have proved in $\eqref{5.15}$ that for $n\sigma>1$,
\begin{eqnarray}
\pp_{\mu}\left(S_{m\sigma,n\sigma}-\widehat{S}_{m\sigma,n\sigma}\right)
&=&\pp_{\mu}\left[\pp_{\mu}\left(S_{m\sigma,n\sigma}|\mathcal{F}_{n\sigma}\right)
-\pp_{\mu}\left(\widehat{S}_{m\sigma,n\sigma}|\mathcal{F}_{n\sigma}\right)\right]\nonumber\\
&\le&e^{\lambda_{1}m\sigma}\langle f_{m\sigma,n\sigma},h\rangle \langle h,\mu\rangle+
e^{-\lambda_{h}(n\sigma-1)}\|\phi\|_{\infty}\langle \widetilde{a}^{1/2}_{2}h,\mu\rangle,\label{6.4}
\end{eqnarray}
where $f_{m\sigma,n\sigma}(x):=\pp_{\delta_{x}}\left(\langle \phi h,I_{m\sigma}\rangle 1_{\{\langle \phi h,I_{m\sigma}\rangle\ge e^{-\lambda_{1}n\sigma}\}}\right)$. Note that by \eqref{5.23},
\begin{eqnarray}
\sum_{n=1}^{{\infty}}\langle f_{m\sigma,n\sigma},h\rangle
&=&\sum_{n=1}^{{\infty}}\int_{E}h(x)\pp_{\delta_{x}}\left(\langle \phi h,I_{m\sigma}\rangle 1_{\{\langle \phi h,I_{m\sigma}\rangle\ge e^{-\lambda_{1}n\sigma}\}}\right)m(dx)\nonumber\\
&\le&\int_{0}^{{\infty}}ds\int_{E}h(x)\pp_{\delta_{x}}\left(\langle \phi h,I_{m\sigma}\rangle 1_{\{\langle \phi h,I_{m\sigma}\rangle\ge e^{-\lambda_{1}s\sigma}\}}\right)m(dx)\nonumber\\
&=&\int_{0}^{{\infty}}\langle f_{m\sigma,s\sigma},h\rangle ds<{\infty}.\label{6.5}
\end{eqnarray}
Hence by \eqref{6.4} and \eqref{6.5},
$$\sum_{n=1}^{{\infty}}\pp_{\mu}\left(S_{m\sigma,n\sigma}-\widehat{S}_{m\sigma,n\sigma}\right)
=\pp_{\mu}\left[\pp_{\mu}\left(S_{m\sigma,n\sigma}|\mathcal{F}_{n\sigma}\right)
-\pp_{\mu}\left(\widehat{S}_{m\sigma,n\sigma}|\mathcal{F}_{n\sigma}\right)\right]<{\infty}.
$$
It follows by Borel-Cantelli lemma that
\begin{equation}
\lim_{n\to{\infty}}\left(S_{m\sigma,n\sigma}-\widehat{S}_{m\sigma,n\sigma}\right)
=\lim_{n\to{\infty}}\left(\pp_{\mu}\left(S_{m\sigma,n\sigma}|\mathcal{F}_{n\sigma}\right)-\pp_{\mu}\left(\widehat{S}_{m\sigma,n\sigma}
|\mathcal{F}_{n\sigma}\right)\right)=0\quad\pp_{\mu}\mbox{-a.s.}\nonumber
\end{equation}
This together with \eqref{6.3} yield that $\lim_{n\to{\infty}}\left(S_{m\sigma,n\sigma}-\pp_{\mu}\left(S_{m\sigma,n\sigma}|\mathcal{F}_{n\sigma}\right)\right)=0$ $\pp_{\mu}$-a.s.\qed

\medskip

\begin{lemma}\label{lem4.2}
Suppose Assumptions 1-3 and Assumption 5 hold. Then for any $\mu\in\mc$, $\sigma>0$ and $\phi\in\mathcal{B}^{+}_{b}(E)$,
\begin{equation}
\lim_{m\to{\infty}}\lim_{n\to{\infty}}e^{\lambda_{1}(m+n)\sigma}\sum_{i=1}^{N_{n\sigma}}\pp_{\mu}\left(\langle \phi h,I^{i,n\sigma}_{m\sigma}\rangle|\mathcal{F}_{n\sigma}\right)=\langle \phi h,h\rangle W^{h}_{\infty}(X)\quad\pp_{\mu}\mbox{-a.s.}\nonumber
\end{equation}
\end{lemma}

\proof For $m,n\in\mathbb{N}$ and $\sigma>0$,  let $I(m,n):=e^{\lambda_{1}(m+n)\sigma}\sum_{i=1}^{N_{n\sigma}}\pp_{\mu}\left(\langle \phi h,I^{i,n\sigma}_{m\sigma}\rangle|\mathcal{F}_{n\sigma}\right)$. Then we have that (see, \eqref{5.19} in the proof of Lemma \ref{lem3.5})
$$I(m,n)=e^{\lambda_{1}n\sigma}\langle \frac{h}{w}\,P^{h}_{m\sigma}\phi, Z_{n\sigma}\rangle -e^{\lambda_{1}n\sigma}\langle \frac{h}{w}\theta^{*}_{\phi h}(m\sigma, \cdot),Z_{n\sigma}\rangle.$$
Since $P^{h}_{m\sigma}\phi,\ \theta^{*}_{\phi h}(m\sigma, \cdot)\in \mathcal{B}^{+}_{b}(E)$, it follows by Assumption 5 that
\begin{equation}
\lim_{n\to{\infty}}I(m,n)=\langle P^{h}_{m\sigma}\phi,h^{2}\rangle W^{h/w}_{\infty}(Z)-\langle \theta^{*}_{\phi h}(m\sigma, \cdot),h^{2}\rangle W^{h/w}_{\infty}(Z)\quad\pp_{\mu}
\mbox{-a.s.}
\label{8.8}
\end{equation}
Recall that $h^{2}m$ is the invariant probability measure for $(\xi^{h};\Pi^{h})$. We have $\langle P^{h}_{m\sigma}\phi,h^{2}\rangle=\langle \phi h,h\rangle$. Moreover by Lemma \ref{lem1.1} and the bounded convergence theorem, $\lim_{m\to{\infty}}\langle \theta^{*}_{\phi h}(m\sigma, \cdot),h^{2}\rangle=0$. Therefore by Theorem \ref{them1} and \eqref{8.8}
$$\lim_{m\to{\infty}}\lim_{n\to{\infty}}I(m,n)=\langle \phi h,h\rangle W^{h/w}_{\infty}(Z)=\langle \phi h,h\rangle W^{h}_{\infty}(X) \quad\pp_{\mu}
\mbox{-a.s.}
$$\qed

\medskip

\begin{lemma}\label{lem4.3}
Suppose Assumptions 1-3, 4(i) and 5 hold.
Then for any
$\mu\in\mc$, $\sigma>0$ and $\phi\in\mathcal{B}^{+}_{b}(E)$,
$$\lim_{n\to{\infty}}e^{\lambda_{1}n\sigma}\langle \phi h,X_{n\sigma}\rangle=\langle \phi h,h\rangle W^{h}_{\infty}(X)\quad\pp_{\mu}
\mbox{-a.s.}
$$
\end{lemma}

\proof Note that for any $m,n\in\mathbb{N}$, $\sigma>0$ and $\phi\in \mathcal{B}^{+}_{b}(E)$,
\begin{eqnarray}
&&e^{\lambda_{1}(m+n)\sigma}\langle \phi h,X_{(m+n)\sigma}\rangle\nonumber\\
&\ge&e^{\lambda_{1}(m+n)\sigma}\sum_{i=1}^{N_{n\sigma}}\langle \phi h,I^{i,n\sigma}_{m\sigma}\rangle\nonumber\\
&=&e^{\lambda_{1}(m+n)\sigma}\sum_{i=1}^{N_{n\sigma}}\langle \phi h,I^{i,n\sigma}_{m\sigma}\rangle-e^{\lambda_{1}(m+n)\sigma}\sum_{i=1}^{N_{n\sigma}}\pp_{\mu}\left(\langle \phi h,I^{i,n\sigma}_{m\sigma}\rangle|\mathcal{F}_{n\sigma}\right)+e^{\lambda_{1}(m+n)\sigma}\sum_{i=1}^{N_{n\sigma}}\pp_{\mu}
\left(\langle \phi h,I^{i,n\sigma}_{m\sigma}\rangle|\mathcal{F}_{n\sigma}\right).\nonumber
\end{eqnarray}
Then by Lemma \ref{lem4.1} and Lemma \ref{lem4.2}, we have
$$\liminf_{n\to{\infty}}e^{\lambda_{1}n\sigma}\langle \phi h,X_{n\sigma}\rangle\ge \langle \phi h,h\rangle W^{h}_{\infty}(X)\quad\pp_{\mu}\mbox{-a.s}.$$
Let $c:=\|\phi\|_{\infty}$. Since $(c-\phi)\in\mathcal{B}^{+}_{b}(E)$, the same argument can be applied to $c-\phi$, and we conclude that under $\pp_{\mu}$,
\begin{eqnarray}
\limsup_{n\to{\infty}}e^{\lambda_{1}n\sigma}\langle \phi h,X_{n\sigma}\rangle
&=&\limsup_{n\to{\infty}}e^{\lambda_{1}n\sigma}\langle ch-(c-\phi)h,X_{n\sigma}\rangle\nonumber\\
&=&cW^{h}_{\infty}(X)-\liminf_{n\to{\infty}}e^{\lambda_{1}n\sigma}\langle (c-\phi)h,X_{n\sigma}\rangle\nonumber\\
&\le&cW^{h}_{\infty}(X)-\langle (c-\phi)h,h\rangle W^{h}_{\infty}(X)=\langle \phi h,h\rangle W^{h}_{\infty}(X).\nonumber
\end{eqnarray}
Hence we complete the proof.\qed

\subsection{From lattice times to continuous time}

In this section we extend the convergence along lattice times in Lemma \ref{lem4.3} to convergence along continuous time and
then prove Theorem \ref{them4}.
Let $\{ U^{\kappa}; \kappa >0\}$ be the resolvent of the semigroup $\{P^{h}_{t}; t\geq 0\}$,
that is,
$$U^{\K}f(x):=\int_{0}^{{\infty}}e^{-\K t}P^{h}_{t}f(x)dt\quad\mbox{ for }x\in E\mbox{ and }f\in \mathcal{B}^{+}_{b}(E).$$
Under Assumption 6, for every $f\in C_{0}(E)$, by the dominated convergence theorem,
\begin{equation}
\|\K U^{\K}f-f\|_{\infty}\le \int_{0}^{{\infty}}\K e^{-\K t }\|P^{h}_{t}f-f\|_{\infty}dt
=\int_{0}^{{\infty}} e^{- s }\|P^{h}_{s/\K}f-f\|_{\infty}ds\to 0\mbox{ as }\K \to{\infty}.\label{6.9}
\end{equation}

\begin{lemma}\label{lem4.4}
Suppose Assumptions 1-5 hold.
Then for any $\mu\in\mc$, $\K>0$ and $\phi\in\mathcal{B}^{+}_{b}(E)$,
\begin{equation}
\lim_{t\to{\infty}}e^{\lambda_{1}t}\langle (\K U^{\K}\phi)h,X_{t}\rangle =\langle \phi h,h\rangle W^{h}_{\infty}(X) \quad \p_{\mu}
\mbox{-a.s.}
\nonumber
\end{equation}
\end{lemma}

\proof The idea of this proof is similar to that of \cite[Proposition 3.14]{EKW}. The main difference is that here we use the $L^{1}$-convergence instead of the $L^{p}$-convergence used there. We consider in the skeleton space for convenience. Since $\K U^{\K}1=1$ and $\K U^{\K}\phi\in\mathcal{B}^{+}_{b}(E)$, the same argument that led to Lemma \ref{lem4.3} can be applied here and it suffices to prove that
\begin{equation}
\liminf_{t\to{\infty}}e^{\lambda_{1}t}\langle (\K U^{\K}\phi)h,X_{t}\rangle
\ge
\langle \phi h,h\rangle W^{h}_{\infty}(X) \quad \pp_{\mu}
\mbox{-a.s.}
\label{6.10}
\end{equation}
Fix $\mu\in\mc$ and $\phi\in\mathcal{B}^{+}_{b}(E)$. Let $g(x):=\K U^{\K}\phi(x)$ and $g_{B}(x):=g(x)1_{B}(x)$ with $B\Subset E$. Note that $\langle (\K U^{\K}\phi)h,h\rangle =\langle \phi h,h\rangle$ since $h^{2}m$ is the invariant probability measure for $(\xi^{h};\Pi^{h})$. If we can show that
\begin{equation}
\liminf_{t\to{\infty}}e^{\lambda_{1}t}\langle gh,X_{t}\rangle
\ge
\langle g_{B} h,h\rangle W^{h}_{\infty}(X) \quad \pp_{\mu}
\mbox{-a.s.}
\label{6.11}
\end{equation}
for any $B\Subset E$, then choose an increasing sequence $B_{n}\Subset E$ with $\bigcup_{n=1}^{{\infty}}B_{n}=E$, \eqref{6.10} follows from the monotone convergence theorem.
Suppose $\sigma>0$. For any $t\in [n\sigma,(n+1)\sigma]$,
\begin{eqnarray}
&&e^{\lambda_{1}t}\langle gh,X_{t}\rangle-\langle g_{B} h,h\rangle W^{h}_{\infty}(X)\nonumber\\
&\ge&\left[e^{\lambda_{1}t}\langle gh,X_{t}\rangle-e^{\lambda_{1}(n+1)\sigma}\pp_{\mu}\left(\langle
g h,X_{(n+1)\sigma}\rangle | \mathcal{F}_{t}\right)\right]\nonumber\\
&+&\left[e^{\lambda_{1}(n+1)\sigma}\pp_{\mu}\left(\langle
g_{B} h,X_{(n+1)\sigma}\rangle | \mathcal{F}_{t}\right)-e^{\lambda_{1}(n+1)\sigma}\pp_{\mu}\left(\langle
g_{B} h,X_{(n+1)\sigma}\rangle | \mathcal{F}_{n\sigma}\right)\right]\nonumber\\
&+&\left[e^{\lambda_{1}(n+1)\sigma}\pp_{\mu}\left(\langle
g_{B} h,X_{(n+1)\sigma}\rangle | \mathcal{F}_{n\sigma}\right)-\langle g_{B} h,h\rangle W^{h}_{\infty}(X)\right]\nonumber\\
&=:&\theta_{1,g}(n,\sigma,t)+\theta_{2,g_{B}}(n,\sigma,t)+\theta_{3,g_{B}}(n,\sigma).\nonumber
\end{eqnarray}
To prove \eqref{6.11}, it suffices to prove that under $\pp_{\mu}$,
\begin{description}
\item{(1)} $\limsup_{\sigma\to 0}\limsup_{n\to{\infty}}\sup_{t\in [n\sigma,(n+1)\sigma]}|\theta_{1,g}(n,\sigma,t)|=0.$
\item{(2)} $\limsup_{n\to{\infty}}\sup_{t\in [n\sigma,(n+1)\sigma]}|\theta_{2,g_{B}}(n,\sigma,t)|=0.$
\item{(3)} $\lim_{n\to{\infty}}|\theta_{3,g_{B}}(n,\sigma)|=0.$
\end{description}
We begin with the proof of $(1)$.
By the Markov property of $X$,
we have
\begin{eqnarray}
\theta_{1,g}(n,\sigma,t)&=&e^{\lambda_{1}t}\langle g h,X_{t}\rangle-e^{\lambda_{1}(n+1)\sigma}\pp_{X_{t}}\left(
\langle gh,X_{(n+1)\sigma-t}\right)\nonumber\\
&=&e^{\lambda_{1}t}\langle g h,X_{t}\rangle-e^{\lambda_{1}(n+1)\sigma}\langle P^{\beta}_{(n+1)\sigma-t}(gh),X_{t}\rangle\nonumber\\
&=&e^{\lambda_{1}t}\langle g h,X_{t}\rangle-e^{\lambda_{1}t}\langle hP^{h}_{(n+1)\sigma-t}(g),X_{t}\rangle.\label{6.12}
\end{eqnarray}
Note that for any $s>0$ and $x\in E$,
\begin{eqnarray}
|g(x)-P^{h}_{s}g(x)|&=&\left|\int_{0}^{{\infty}}\K e^{-\K t}P^{h}_{t}\phi(x)dt-\int_{0}^{{\infty}}\K e^{-\K t}P^{h}_{t+s}\phi(x)dt\right|\nonumber\\
&=&\left|\int_{0}^{{\infty}}\K e^{-\K t}P^{h}_{t}\phi(x)dt-\int_{s}^{{\infty}}\K e^{-\K t}P^{h}_{t}\phi(x)dt\right|\nonumber\\
&\le&2(1-e^{-\K s})\|\phi\|_{\infty}.\nonumber
\end{eqnarray}
It follows by this and \eqref{6.12} that
\begin{equation}
\sup_{t\in[n\sigma,(n+1)\sigma]}|\theta_{1,g}(n,\sigma,t)|\le 2(1-e^{-\K \sigma})\|\phi\|_{\infty}\sup_{t\in[n\sigma,(n+1)\sigma]}W^{h}_{t}(X).\label{6.13}
 \end{equation}
Since $W^{h}_{t}(X)$ converges to a finite limit almost surely, we obtain (1) from \eqref{6.13} by letting $n\to{\infty}$ and $\sigma\to 0$.
For the proof of (3), note that by the Markov property of $X$,
\begin{eqnarray}
\theta_{3,g}(n,\sigma)&=&e^{\lambda_{1}(n+1)\sigma}\pp_{X_{n\sigma}}\left(\langle g_{B}h,X_{\sigma}\rangle\right)
-\langle g_{B}h,h\rangle W^{h}_{\infty}(X)\nonumber\\
&=&e^{\lambda_{1}(n+1)\sigma}\langle P^{\beta}_{\sigma}(g_{B}h),X_{n\sigma}\rangle
-\langle g_{B}h,h\rangle W^{h}_{\infty}(X)\nonumber\\
&=&e^{\lambda_{1}n\sigma}\langle hP^{h}_{\sigma}(g_{B}),X_{n\sigma}\rangle
-\langle g_{B}h,h\rangle W^{h}_{\infty}(X).\label{6.14}
\end{eqnarray}
Since $P^{h}_{\sigma}g_{B}\in\mathcal{B}^{+}_{b}(E)$, it follows from Lemma \ref{lem4.3} that
$$\lim_{n\to{\infty}}e^{\lambda_{1}n\sigma}\langle hP^{h}_{\sigma}(g_{B}),X_{n\sigma}\rangle=
\langle hP^{h}_{\sigma}(g_{B}),h\rangle W^{h}_{\infty}(X)=\langle g_{B} h, h\rangle W^{h}_{\infty}(X)\quad\pp_{\mu}
\mbox{-a.s.}
$$
By this and \eqref{6.14} we obtain (3). It remains to prove (2). $\theta_{2,g_{B}}(n,\sigma,t)$ can be written as $$\theta_{2,g_{B}}(n,\sigma,t)=\theta^{(1)}_{2,g_{B}}(n,\sigma,t)+\theta^{(2)}_{2,g_{B}}(n,\sigma,t)+\theta^{(3)}_{2,g_{B}}(n,\sigma,t)$$
where
\begin{eqnarray}
\theta^{(1)}_{2,g_{B}}(n,\sigma,t)&:=&e^{\lambda_{1}(n+1)\sigma}\pp_{\mu}\left(\langle g_{B}h,X^{*}_{(n+1)\sigma}+I^{*,n\sigma}_{\sigma}\rangle |\mathcal{F}_{t}\right)\nonumber\\
&&-e^{\lambda_{1}(n+1)\sigma}\pp_{\mu}\left(\langle g_{B}h,X^{*}_{(n+1)\sigma}+I^{*,n\sigma}_{\sigma}\rangle |\mathcal{F}_{n\sigma}\right),\nonumber\\
\theta^{(2)}_{2,g_{B}}(n,\sigma,t)&:=&e^{\lambda_{1}(n+1)\sigma}\pp_{\mu}\left(\sum_{i=1}^{N_{n\sigma}}\langle g_{B}h,I^{i,n\sigma}_{\sigma}\rangle 1_{\{\langle g_{B}h,I^{i,n\sigma}_{\sigma}\rangle<e^{-\lambda_{1}n\sigma}\}}|\mathcal{F}_{t}\right)\nonumber\\
&&-e^{\lambda_{1}(n+1)\sigma}\pp_{\mu}\left(\sum_{i=1}^{N_{n\sigma}}\langle g_{B}h,I^{i,n\sigma}_{\sigma}\rangle 1_{\{\langle g_{B}h,I^{i,n\sigma}_{\sigma}\rangle<e^{-\lambda_{1}n\sigma}\}} |\mathcal{F}_{n\sigma}\right),\nonumber\\
\theta^{(3)}_{2,g_{B}}(n,\sigma,t)&:=&e^{\lambda_{1}(n+1)\sigma}\pp_{\mu}\left(\sum_{i=1}^{N_{n\sigma}}\langle g_{B}h,I^{i,n\sigma}_{\sigma}\rangle 1_{\{\langle g_{B}h,I^{i,n\sigma}_{\sigma}\rangle\ge e^{-\lambda_{1}n\sigma}\}}|\mathcal{F}_{t}\right)\nonumber\\
&&-e^{\lambda_{1}(n+1)\sigma}\pp_{\mu}\left(\sum_{i=1}^{N_{n\sigma}}\langle g_{B}h,I^{i,n\sigma}_{\sigma}\rangle 1_{\{\langle g_{B}h,I^{i,n\sigma}_{\sigma}\rangle\ge e^{-\lambda_{1}n\sigma}\}} |\mathcal{F}_{n\sigma}\right).\nonumber
\end{eqnarray}
Clearly for $i=1,2,3$, $\{\theta^{(i)}_{2,g_{B}}(n,\sigma,t):t\in [n\sigma,(n+1)\sigma]\}$ are $\pp_{\mu}$-martingales.
Let $$S^{g_{B}}_{\sigma,n\sigma}:=e^{\lambda_{1}(n+1)\sigma}\sum_{i=1}^{N_{n\sigma}}\langle g_{B}h,I^{i,n\sigma}_{\sigma}\rangle$$
and
$$\widehat{S}^{g_{B}}_{\sigma,n\sigma}:=e^{\lambda_{1}(n+1)\sigma}\sum_{i=1}^{N_{n\sigma}}\langle g_{B}h,I^{i,n\sigma}_{\sigma}\rangle 1_{\{\langle g_{B}h,I^{i,n\sigma}_{\sigma}\rangle<e^{-\lambda_{1}n\sigma}\}}.$$
Then $\theta^{(2)}_{2,g_{B}}(n,\sigma,t)=\pp_{\mu}\left(\widehat{S}^{g_{B}}_{\sigma,n\sigma}|\mathcal{F}_{t}\right)
-\pp_{\mu}\left(\widehat{S}^{g_{B}}_{\sigma,n\sigma}|\mathcal{F}_{n\sigma}\right)$, and
$$
\theta^{(3)}_{2,g_{B}}(n,\sigma,t)=\pp_{\mu}\left(S^{g_{B}}_{\sigma,n\sigma}-\widehat{S}^{g_{B}}_{\sigma,n\sigma}|\mathcal{F}_{t}\right)
-\pp_{\mu}\left(S^{g_{B}}_{\sigma,n\sigma}-\widehat{S}^{g_{B}}_{\sigma,n\sigma}|\mathcal{F}_{n\sigma}\right).
$$
Suppose ${\eps}>0$. By Doob's maximal inequality and Jensen's inequality, we have
\begin{eqnarray}
&&\pp_{\mu}\left(\sup_{t\in [n\sigma,(n+1)\sigma]}|\theta^{(2)}_{2,g_{B}}(n,\sigma,t)|>{\eps}\right)\nonumber\\
&\le&\frac{1}{{\eps}^{2}}\pp_{\mu}\left(|\theta^{(2)}_{2,g_{B}}(n,\sigma,(n+1)\sigma)|^{2}\right)\nonumber\\
&=&\frac{1}{{\eps}^{2}}\pp_{\mu}\left[\left|\pp_{\mu}\left(\widehat{S}^{g_{B}}_{\sigma,n\sigma}-\pp_{\mu}\left(\widehat{S}^{g_{B}}_{\sigma,n\sigma}|\mathcal{F}_{n\sigma}\right)|\mathcal{F}_{(n+1)\sigma}\right)\right|^{2}\right]\nonumber\\
&\le&\frac{1}{{\eps}^{2}}\pp_{\mu}\left(\left|\widehat{S}^{g_{B}}_{\sigma,n\sigma}-\pp_{\mu}\left(\widehat{S}^{g_{B}}_{\sigma,n\sigma}|\mathcal{F}_{n\sigma}\right)\right|^{2}\right).\nonumber
\end{eqnarray}
We have proved in the proof of Lemma \ref{lem4.1} that $\sum_{n=1}^{{\infty}}\pp_{\mu}\left(|\widehat{S}^{g_{B}}_{\sigma,n\sigma}-\pp_{\mu}\left(\widehat{S}^{g_{B}}_{\sigma,n\sigma}
|\mathcal{F}_{n\sigma}\right)|^{2}\right)<{\infty}$. Thus $$\sum_{n=1}^{{\infty}}\pp_{\mu}\left(\sup_{t\in [n\sigma,(n+1)\sigma]}|\theta^{(2)}_{2,g_{B}}(n,\sigma,t)|>{\eps}\right)<{\infty}$$
and consequently by Borel-Cantelli lemma
\begin{equation}
\lim_{n\to{\infty}}\sup_{t\in[n\sigma,(n+1)\sigma]}|\theta^{(2)}_{2,g_{B}}(n,\sigma,t)|=0\quad\pp_{\mu}\mbox{-a.s.}\label{6.25}
\end{equation}
Similarly by Doob's maximal inequality and Jensen's inequality,
\begin{eqnarray}
&&\pp_{\mu}\left(\sup_{t\in[n\sigma,(n+1)\sigma]}|\theta^{(3)}_{2,g_{B}}(n,\sigma,t)|>{\eps}\right)\nonumber\\
&\le&\frac{1}{{\eps}}\pp_{\mu}\left(|\theta^{(3)}_{2,g_{B}}(n,\sigma,(n+1)\sigma)|\right)\nonumber\\
&=&\frac{1}{{\eps}}\pp_{\mu}\left[\left|\pp_{\mu}\left(S^{g_{B}}_{\sigma,n\sigma}-\widehat{S}^{g_{B}}_{\sigma,n\sigma}
-\pp_{\mu}\left(S^{g_{B}}_{\sigma,n\sigma}-\widehat{S}^{g_{B}}_{\sigma,n\sigma}|\mathcal{F}_{n\sigma}\right)|\mathcal{F}_{(n+1)\sigma}\right)\right|\right]\nonumber\\
&\le&\frac{1}{{\eps}}\left[\pp_{\mu}\left(S^{g_{B}}_{\sigma,n\sigma}-\widehat{S}^{g_{B}}_{\sigma,n\sigma}\right)+
\pp_{\mu}\left(\pp_{\mu}\left(S^{g_{B}}_{\sigma,n\sigma}-\widehat{S}^{g_{B}}_{\sigma,n\sigma}|\mathcal{F}_{n\sigma}\right)\right)\right]\nonumber\\
&=&\frac{2}{{\eps}}\pp_{\mu}\left(S^{g_{B}}_{\sigma,n\sigma}-\widehat{S}^{g_{B}}_{\sigma,n\sigma}\right).\nonumber
\end{eqnarray}
We have showed in the proof of Lemma \ref{lem4.1} that $\sum_{n=1}^{{\infty}}\pp_{\mu}\left(S^{g_{B}}_{\sigma,n\sigma}-\widehat{S}^{g_{B}}_{\sigma,n\sigma}\right)<{\infty}$. Thus by Borel-Cantelli lemma
\begin{equation}
\lim_{n\to{\infty}}\sup_{t\in[n\sigma,(n+1)\sigma]}|\theta^{(3)}_{2,g_{B}}(n,\sigma,t)|=0\quad\pp_{\mu}
\mbox{-a.s.}
\label{6.24}
\end{equation}
Using Doob's maximal inequality and Jensen's inequality again, we get
\begin{equation}\label{6.17}
\pp_{\mu}\left(\sup_{t\in[n\sigma,(n+1)\sigma]}|\theta^{(1)}_{2,g_{B}}(n,\sigma,t)|>{\eps}\right)
\le \frac{1}{{\eps}^{2}}e^{2\lambda_{1}(n+1)\sigma}\pp_{\mu}\left[\mbox{Var}\left(\langle g_{B}h,X^{*}_{(n+1)\sigma}+I^{*,n\sigma}_{\sigma}\rangle |\mathcal{F}_{n\sigma}\right)\right].
\end{equation}
Recall that given $\mathcal{F}_{t}$, $(X^{*}_{s+t}+I^{*,t}_{s})_{s\ge 0}$ is equal in distribution to $((X^{*}_{s})_{s\ge 0};\pp_{X_{t}})$.
It is known that for $\nu\in\mf$ and $f\in \mathcal{B}^{+}_{b}(E)$, the second moment of $\langle f,X_{t}^{*}\rangle$ can be expressed as
$$\mbox{Var}_{\nu}\left(\langle f,X^{*}_{t}\rangle\right)=\int_{0}^{t}\langle P^{\beta^{*}}_{s}\left[(2\alpha+b^{*})(P^{\beta^{*}}_{t-s}f)^{2}\right],\nu\rangle ds,$$
where $b^{*}(x):=\int_{(0,{\infty})}r^{2}\pi^{*}(x,dr)=\int_{(0,{\infty})}r^{2}e^{-w(x)r}\pi(x,dr)$.
Thus we can continue the calculation in \eqref{6.17} to get
\begin{eqnarray}
&&\pp_{\mu}\left(\sup_{t\in[n\sigma,(n+1)\sigma]}|\theta^{(1)}_{2,g_{B}}(n,\sigma,t)|>{\eps}\right)\nonumber\\
&\le&\frac{1}{{\eps}^{2}}e^{2\lambda_{1}(n+1)\sigma}\pp_{\mu}\left[\mbox{Var}_{X_{n\sigma}}\left(\langle g_{B}h,X^{*}_{\sigma}\rangle\right)\right]\nonumber\\
&=&\frac{1}{{\eps}^{2}}e^{2\lambda_{1}(n+1)\sigma}\pp_{\mu}\left[\int_{0}^{\sigma}\langle P^{\beta^{*}}_{s}\left[(2\alpha+b^{*})(P^{\beta^{*}}_{\sigma-s}(g_{B}h))^{2}\right],X_{n\sigma} \rangle ds\right]\nonumber\\
&=&\frac{1}{{\eps}^{2}}e^{2\lambda_{1}(n+1)\sigma}\int_{0}^{\sigma}\langle P^{\beta}_{n\sigma}P^{\beta^{*}}_{s}\left[(2\alpha+b^{*})(P^{\beta^{*}}_{\sigma-s}(g_{B}h))^{2}\right],\mu \rangle ds\nonumber\\
&\le&\frac{1}{{\eps}^{2}}e^{2\lambda_{1}(n+1)\sigma}\int_{0}^{\sigma}\langle P^{\beta}_{n\sigma+s}\left[(2\alpha+b^{*})(P^{\beta^{*}}_{\sigma-s}(g_{B}h))^{2}\right],\mu \rangle ds\nonumber\\
&=&\frac{1}{{\eps}^{2}}e^{\lambda_{1}n\sigma+2\lambda_{1}\sigma}\int_{0}^{\sigma}e^{-\lambda_{1}s}\langle P^{h}_{n\sigma+s}\left[\frac{2\alpha+b^{*}}{h}(P^{\beta^{*}}_{\sigma-s}(g_{B}h))^{2}\right],h\mu \rangle ds.\label{6.19}
\end{eqnarray}
Since $g_{B}$ is compactly supported, for $s\in [0,\sigma)$,
$\|P^{\beta^{*}}_{\sigma-s}(g_{B}h)\|_{\infty}\le \|P^{\beta}_{\sigma-s}(g_{B}h)\|_{\infty}\le e^{\|\beta\|_{\infty}\sigma}\|g_{B}h\|_{\infty}$.
In addition, for each $x\in E$,
$P^{\beta^{*}}_{\sigma-s}(g_{B}h)(x)=e^{-\lambda_{1}(\sigma-s)}h(x)\theta^{*}_{g_{B}h}(\sigma-s,x)\le e^{-\lambda_{1}\sigma}\|\theta^{*}_{g_{B}h}\|_{\infty}h(x)$. Thus we get for $s\in [0,\sigma)$,
\begin{equation}
\left(P^{\beta^{*}}_{\sigma-s}(g_{B}h)(x)\right)^{2}\le c_{1}\wedge c_{2}h(x)\wedge c_{3}h(x)^{2}\quad\mbox{ for all }x\in E.\label{6.21}
\end{equation}
Here $c_{1},c_{2},c_{3}$ are positive constants. We have by \eqref{6.21}
\begin{eqnarray}
\int_{0}^{\sigma}e^{-\lambda_{1}s}\langle P^{h}_{n\sigma+s}\left[\frac{2\alpha}{h}(P^{\beta^{*}}_{\sigma-s}(g_{B}h))^{2}\right],h\mu \rangle ds
&\le&\int_{0}^{\sigma}e^{-\lambda_{1}s}\langle P^{h}_{n\sigma+s}\left(2c_{2}\alpha\right),h\mu\rangle ds\nonumber\\
&\le&2c_{2}\|\alpha\|_{\infty}\langle h,\mu\rangle \int_{0}^{\sigma}e^{-\lambda_{1}s}ds.\label{6.22}
\end{eqnarray}
On the other hand,
\begin{eqnarray}
\int_{0}^{\sigma}e^{-\lambda_{1}s}\langle P^{h}_{n\sigma+s}\left[\frac{b^{*}}{h}(P^{\beta^{*}}_{\sigma-s}(g_{B}h))^{2}\right],h\mu \rangle ds
&\le&\int_{0}^{\sigma}e^{-\lambda_{1}s}\langle P^{h}_{n\sigma+s}\left(c_{1}\frac{b^{*}}{h}\wedge c_{2}b^{*}\wedge c_{3}b^{*}h\right),h\mu\rangle ds.\nonumber
\end{eqnarray}
Let $l(x):=c_{1}\frac{b^{*}}{h}(x)\wedge c_{2}b^{*}(x)\wedge c_{3}b^{*}(x)h(x)$ for $x\in E$. Assumption 4(ii) implies that $l(x)\in L^{2}(E,\widetilde{m})$. Thus by \eqref{p3} for $n\sigma>1$ and $x\in E$,
$$P^{h}_{n\sigma+s}(l)(x)\le \langle l,h^{2}\rangle+e^{-\lambda_{h}(n\sigma+s-1)}\widetilde{a}_{2}(x)^{1/2}\|l\|_{L^{2}(E,\widetilde{m})}.$$
Consequently for $n\sigma>1$,
\begin{eqnarray}
&&\int_{0}^{\sigma}e^{-\lambda_{1}s}\langle P^{h}_{n\sigma+s}\left[\frac{b^{*}}{h}(P^{\beta^{*}}_{\sigma-s}(g_{B}h))^{2}\right],h\mu \rangle ds\nonumber\\
&\le&\int_{0}^{\sigma}e^{-\lambda_{1}s}\langle P^{h}_{n\sigma+s}(l),h\mu \rangle ds\nonumber\\
&\le& \langle l,h^{2}\rangle \langle h,\mu\rangle \int_{0}^{\sigma}e^{-\lambda_{1}s}ds+e^{-\lambda_{h}(n\sigma-1)}\|l\|_{L^{2}(E,\widetilde{m})}\langle \widetilde{a}^{1/2}_{2}h,\mu\rangle\int_{0}^{\sigma}e^{-(\lambda_{1}+\lambda_{h})s}ds.\label{6.23}
\end{eqnarray}
Since $\lambda_{1}<0$, it follows from \eqref{6.22} and \eqref{6.23} that
$$\sum_{n=1}^{\infty}e^{\lambda_{1}n\sigma}\int_{0}^{\sigma}e^{-\lambda_{1}s}\langle P^{h}_{n\sigma+s}\left[\frac{2\alpha+b^{*}}{h}(P^{\beta^{*}}_{\sigma-s}(g_{B}h))^{2}\right],h\mu \rangle ds<{\infty}.$$
Thus by \eqref{6.19}, we get $\sum_{n=1}^{{\infty}}\pp_{\mu}\left(\sup_{t\in[n\sigma,(n+1)\sigma]}|\theta^{(1)}_{2,g_{B}}(n,\sigma,t)|>{\eps}\right)<{\infty}$.
By Borel-Cantelli Lemma
\begin{equation}
\lim_{n\to{\infty}}\sup_{t\in[n\sigma,(n+1)\sigma]}|\theta^{(1)}_{2,g_{B}}(n,\sigma,t)|=0\quad\pp_{\mu}
\mbox{-a.s.}
\label{6.26}
\end{equation}
We obtain (2) by \eqref{6.25}, \eqref{6.24} and \eqref{6.26}. This completes the proof.\qed

\medskip

\begin{lemma}\label{lem4.5}
Suppose Assumptions 1-6 hold. Then for any
$\mu\in \mc$ and $\phi\in C^{+}_{0}(E)$,
\begin{equation}
\lim_{t\to{\infty}}e^{\lambda_{1}t}\langle \phi h,X_{t}\rangle =\langle \phi h,h\rangle W^{h}_{\infty}(X)\quad\pp_{\mu}
\mbox{-a.s.}
\label{6.20}
\end{equation}
\end{lemma}

\proof For any $\K>0$, we have
\begin{eqnarray}
\left|e^{\lambda_{1}t}\langle \phi h,X_{t}\rangle-\langle \phi h,h\rangle W^{h}_{\infty}(X)\right|
&\le&e^{\lambda_{1}t}\langle |\K U^{\K}\phi-\phi|h,X_{t}\rangle+\left|e^{\lambda_{1}t}\langle (\K U^{\K} \phi) h,X_{t}\rangle-\langle \phi h,h\rangle W^{h}_{\infty}(X)\right|\nonumber\\
&\le&\|\K U^{\K}\phi-\phi\|_{\infty}W^{h}_{t}(X)+\left|e^{\lambda_{1}t}\langle (\K U^{\K} \phi) h,X_{t}\rangle-\langle \phi h,h\rangle W^{h}_{\infty}(X)\right|.\nonumber
\end{eqnarray}
In view of \eqref{6.9} and Lemma \ref{lem4.4}, we conclude \eqref{6.20} by letting $t\to{\infty}$ and $\K\to{\infty}$.\qed

\medskip

\noindent\textit{Proof of Theorem \ref{them4}.} This theorem follows immediately from Lemma \ref{lem4.6}, Lemma \ref{lem4.7} and Lemma \ref{lem4.5}.\qed

\section{Examples}\label{sec8}

In this section, we give examples of superprocesses where Assumptions 1-6 are satisfied and Theorem \ref{them1}-Theorem \ref{them4} hold.
 Since the paper is already long, we leave the detailed verifications of assumptions for these examples in the
Appendix of this paper.
We notice that Example \ref{E:8.1} and \ref{E:8.2} are also studied in \cite{EKW}. In Example \ref{E:8.5} we consider a class of super $\alpha$-stable-like  processes which includes super $\alpha$-stable processes as a special case.

\medskip

\begin{example}\label{E:8.1} \rm
 Suppose $E=\R^{d}$ ($d\ge 1$) and $\xi=(\xi_{t};\Pi_{x},x\in \R^{d})$ is an inward Ornstein-Uhlenbeck (OU) process on $\mathbb{R}^{d}$ with infinitesimal generator
$$\mathcal{L}:=\frac{1}{2}\sigma^{2}\Delta-cx\cdot\nabla \quad\mbox{ on }\R^{d}$$
where $\sigma,c>0$.
Without loss of generality, we assume that $\sigma=1$.
Let $dx$ denote the Lebesgue measure on $\R^{d}$ and $m(dx):=\left(\frac{c}{\pi}\right)^{d/2}e^{-c|x|^{2}}dx$. Then $\xi$ is symmetric with respect to the probability measure $m$. Let $\psi(\lambda)$ be a spatially independent branching mechanism
given by \eqref{psi2}
with $\beta>0$, $\psi({\infty})={\infty}$ and
$\int_{(1,{\infty})}y^{p}\pi(dy)<{\infty}$ for some $p\in (1,2]$.
For this $(\xi,\psi)$-superdiffusion
Assumptions 1-6 are
satisfied with $\lambda_{1}=-\beta$, $h(x)\equiv 1$ and
$w(x)\equiv
z_{\psi}
$ where $z_{\psi}:=\sup\{\lambda\ge 0:\ \psi(\lambda)\le 0\}\in (0,{\infty})$.
\end{example}

\medskip

\begin{example} \label{E:8.2} \rm
 Suppose $E=\mathbb{R}^{d}$ ($d\ge 1$) and $\xi:=(\xi_{t};\Pi_{x},x\in\R^{d})$ is an outward OU process on $\R^{d}$ with infinitesimal generator
$$\mathcal{L}:=\frac{1}{2}\sigma^{2}\Delta+cx\cdot\nabla \quad\mbox{ on }\R^{d}$$
where $\sigma,c>0$. Without loss of generality, we assume $\sigma=1$.  Let $m(dx):=\left(\frac{c}{\pi}\right)^{-d/2}e^{c|x|^{2}}dx$. Then $\xi$ is symmetric with respect to $m$. Let $\psi(\lambda)$ be a spatially independent branching mechanism
given by \eqref{psi2}
with $\beta>c d$, $\psi({\infty})={\infty}$ and
$\int_{(1,{\infty})}y^{p}\pi(dy)<{\infty}$ for some $p\in (1,2]$.
Then Assumptions 1-6 are satisfied with $\lambda_{1}=cd-\beta$, $h(x)=\left(\frac{c}{\pi}\right)^{d/2}e^{-c|x|^{2}}$ and $w(x)\equiv z_{\psi}$ where $z_{\psi}:=\sup\{\lambda\ge 0:\ \psi(\lambda)\le 0\}\in (0,{\infty})$.
\end{example}

\medskip

\begin{example}\label{E:8.3} \rm
 Let $(Y=(Y_{t})_{t\ge 0},\Pi_{x})$ be a diffusion on $\mathbb{R}^{d}$ ($d\ge 3$) with generator
$$
\mathcal{A}
:=\rho(x)^{-1}\sum_{i,j=1}^{d}\frac{\partial}{\partial x_{i}}\left(\rho(x)a_{ij}(x)\frac{\partial}{\partial x_{j}}\right)\quad\mbox{ on }\mathbb{R}^{d},$$
where the diffusion matrix $A(x)=(a_{ij}(x))_{ij}$ is uniformly elliptic and symmetric with $a_{ij}\in C^{1}_{b}(\mathbb{R}^{d})$ and function $\rho\in C^{1}_{b}(\mathbb{R}^{d})$ is bounded between two positive constants.
Here $C^{1}_{b}(\mathbb{R}^{d})$ denotes the space of bounded continuous functions on $\R^{d}$ whose first order derivatives are bounded and continuous.
Clearly this includes Brownian motion as a special case.
Note that we can rewrite
$\mathcal{A}$
as
$$
\mathcal{A}
=\sum_{i,j=1}^{d}\frac{\partial}{\partial x_{i}}\left(a_{ij}(x)\frac{\partial}{\partial x_{j}}\right)
+\sum_{j=1}^{d}\left(\sum_{i=1}^{d}a_{ij}(x)\frac{\partial \log \rho}{\partial x_{i}}\right)\frac{\partial}{\partial x_{j}}\quad\mbox{ on }\R^{d}.$$
Define $m(dx):=\rho(x)dx$. Then $Y$ is symmetric with respect to $m$.
Suppose $E$ is a bounded $C^{1,1}$ domain in $\mathbb{R}^{d}$ and $\xi:=Y^{E}$ is the subprocess of $Y$ killing upon leaving $E$. Let $\psi_{\beta}(x,\lambda)
:=
-\beta(x)\lambda+\alpha(x)\lambda^{2}+\int_{(0,{\infty})}\left(e^{-\lambda y}-1+\lambda y\right)\pi(x,dy)$ where $\beta\in C^{1}_{b}(\R^{d})$,
$0<\alpha\in C^{1}_{b}(\R^{d})\cap \mathcal{B}^{+}(\R^{d})$,
$\pi(x,dy)$ satisfies
\begin{equation}\sup_{x\in\R^{d}}\int_{(0,{\infty})}y\log^{*}y\pi(x,dy)<{\infty},\label{8.22}
\end{equation}
and all the first partial derivatives of $\psi_{\beta}(x,\lambda)$ are continuous. We will show in
Appendix \ref{verification}
that
Assumptions 1-6 are satisfied for such $(\xi,\psi_{\beta})$-superdiffusions.

\end{example}

\medskip

\begin{example} \label{E:8.4} \rm
Suppose $E$ is a locally compact separable metric space and $m$ is a $\sigma$-finite nonnegative Radon measure on $E$ with full support. Suppose $\xi$ is
an $m$-symmetric Hunt process on $E$ with transition density function $p(t,x,y)$ which is positive, continuous and symmetric in $(x,y)$.
Define $a_{t}(x):=p(t,x,x)$ for $x\in E$ and $t>0$. We assume that
\begin{description}
\item{(1)} for each $t>0$, $a_{t}(x)\in L^{1}(E,m)$;
\item{(2)} there exists $t_{0}>0$ such that $a_{t_{0}}(x)\in L^{2}(E,m)$.
\end{description}
Define $\hat{a}_{t}(x):=p^{\beta}(t,x,x)$ for $x\in E$ and $t>0$. Under conditions (1) and (2), $\hat{a}_{t}$ satisfies
\begin{description}
\item{(i)} for any $t>0$, $\hat{a}_{t}(x)\in L^{1}(E,m)$;
\item{(ii)} there exists $t_{0}>0$ such that $\hat{a}_{t}(x)\in L^{2}(E,m)$ for all $t\ge t_{0}$.
\end{description}
Property (i) implies that for $t>0$, the Feynman-Kac semigroup
$P^{\beta}_{t}$ is a Hilbert-Schmidt operator in $L^{2}(E,m)$, and hence compact. Let
$\mathcal{L}^{(\beta )}$
be the infinitesimal operator of $P^{\beta}_{t}$, and
$\sigma(\mathcal{L}^{\beta})$ be the spectrum of the self-adjoint operator $\mathcal{L}^{(\beta)}$.
We know by \cite{Pazy} that $\sigma(\mathcal{L}^{(\beta)})$ consists of at most countable eigenvalues. By Jentzsch's theorem,
$\lambda_{1}:=\inf \{-\lambda:\lambda\in\sigma(\mathcal{L}^{(\beta)})\}$
is an eigenvalue of multiplicity 1, and the corresponding eigenfunction $h$ can be chosen to be continuous and positive on $E$ with $\int_{E}h(x)^{2}m(dx)=1$. We assume $\lambda_{1}<0$.
The $h$-transformed semigroup
$P^{h}_{t}$ admits an integral kernel $p^{h}(t,x,y)$ with respect to the measure $\widetilde{m}(dx):=h(x)^{2}m(dx)$, which
is related to $p^{\beta}(t,x,y)$ by
\eqref{density-trans}.
Define $\widetilde{a}_{t}(x):=p^{h}(t,x,x)=e^{\lambda_{1}t}\frac{\hat{a}_{t}(x)}{h(x)^{2}}$ for $x\in E$.
We assume in addition that
\begin{description}
\item{(3)} for all $f\in C_{0}(E)$, $\|P^{h}_{t}f-f\|_{\infty}\to 0$ as $t\to 0$.
\end{description}
Let $X$ be a $(P_{t},\psi_{\beta})$-superprocess with spatial motion $\xi$ satisfying conditions (1)-(3)
and the branching mechanism $\psi_{\beta}(x,\lambda)$ given by \eqref{2.1}
with $\beta(x)\in\mathcal{B}_{b}(E)$, $\alpha\in\mathcal{B}^{+}_{b}(E)$ and $\pi(x,dy)$ satisfying
\begin{equation}
\sup_{x\in E}\int_{(0,{\infty})}y^{2}\pi(x,dy)<{\infty},\label{8.1}
\end{equation}
We assume in addition that $\psi_{\beta}$ satisfies
condition \eqref{condi1} in  Subsection \ref{sup-Br}.
Then  our Assumptions 1-6 are satisfied by this class of superprocesses with $\lambda_{1},h$ defined as above and $w(x):=-\log \p_{\delta_{x}}\left(\exists t\ge 0:\ \langle 1,X_{t}\rangle =0\right)$.
This example covers Example 4.1-4.5 in \cite{CRSZ}.
\end{example}

\medskip

\begin{example}\label{E:8.5} \rm
Suppose $E=\mathbb{R}^{d}\ (d\ge 1)$ and $m(dx)=dx$ is the Lebesgue measure on $\mathbb{R}^{d}$. Suppose $(\xi;\Pi_{x},x\in\mathbb{R}^{d})$ is
an $\alpha$-stable-like process
on $\mathbb{R}^{d}$ with $\alpha\in (0,2)$.
An $\alpha$-stable-like process on $\R^{d}$ is a symmetric Feller process on $\R^{d}$ whose Dirichlet form $(\mathcal{E},\mathcal{F})$ on $L^{2}(\R^{d};dx)$ is given by
$$\mathcal{F}=\left\{u\in L^{2}(\R^{d};dx):\ \int\int_{\R^{d}\times \R^{d}}\frac{(u(x)-u(y))^{2}}{|x-y|^{d+\alpha}}dxdy<{\infty}\right\},$$
$$\mathcal{E}(u,u)=\int\int_{\R^{d}\times \R^{d}}(u(x)-u(y))^{2}\frac{c(x,y)}{|x-y|^{d+\alpha}}dxdy,$$
where $c(x,y)$ is a symmetric function on $\R^{d}\times \R^{d}$ that is bounded between two positive constants.
Clearly this includes symmetric $\alpha$-stable process on $\R^{d}$ as a special case.
Let $\psi_{\beta}(x,\lambda)$ be a branching mechanism given by \eqref{2.1}
where $\beta\in \mathcal{B}_{b}(\mathbb{R}^{d})$ has a compact support,
 $\alpha\in \mathcal{B}^{+}_{b}(\R^{d})$
and $\pi$ satisfies \eqref{8.1}.
Let $X$ be a $(\xi,\psi_{\beta})$-superprocess.
Here we assume in addition that $\psi_{\beta}$ satisfies
condition \eqref{condi1} of Section \ref{sup-Br}.
It is proved in \cite[Theorem 4.14]{CK} that
the transition density function $p(t,x,y)$ of a
$\alpha$-stable-like process
is bounded and continuous in $(x,y)$ for each $t>0$ and
\begin{equation}
p(t,x,y)\asymp t^{-d/\alpha}\wedge \frac{t}{|x-y|^{d+\alpha}}\quad\mbox{ for all } x,y\in \mathbb{R}^{d}\mbox{ and }t>0.\label{7.1}
\end{equation}
 It follows that for each $t>0$, the transition semigroup $P_{t}$ of $\xi$ maps bounded functions to continuous functions, and $P_{t}$ is bounded from $L^{p}(\R^{d},dx)$ to $L^{q}(\R^{d},dx)$ for $1\le p\le q\le {\infty}$.
Furthermore, for each $t>0$, $P^{\beta}_{t}$ maps bounded functions to continuous functions, and $P^{\beta}_{t}$ is bounded from $L^{p}(\R^{d},dx)$ to $L^{q}(\R^{d},dx)$ for $1\le p\le q\le {\infty}$.
Define $\mathcal{E}_{1}(u,u):=\mathcal{E}(u,u)+\int_{\R^{d}}u(x)^{2}dx$ for $u\in\mathcal{F}$, and $\lambda_{0}:=
\inf\{\mathcal{E}(u,u):\ u\in \mathcal{F}\mbox{ with }\int_{\R^{d}}u(x)^{2}dx=1\}\ge 0$. The semigroup $P^{\beta}_{t}$ associates with a quadratic form $(\mathcal{E}^{(\beta ) },\mathcal{F})$ where
$$\mathcal{E}^{(\beta ) }(u,u)=\mathcal{E}(u,u)-\int_{\R^{d}}u(x)^{2}\beta(x)dx\quad\mbox{ for }u\in\mathcal{F}.$$
Define $\lambda_{1}:=
\inf\{\mathcal{E}^{(\beta ) }(u,u):\ u\in \mathcal{F}\mbox{ with }\int_{\R^{d}}u(x)^{2}dx=1\}\ge-
\|\beta^{+}\|_{\infty}
.$
We use $\sigma(\mathcal{E}^{(\beta ) })$ to denote the spectrum of
$- {\mathcal L}^{(\beta)}$, where ${\mathcal L}^{(\beta)}$ is the self-adjoint operator associated with $(\mathcal{E}^{(\beta ) }, \mathcal{F})$.
Clearly the the Dirichlet form $(\mathcal{E},\mathcal{F})$ is comparable to that of the symmetric $\alpha$-stable process.
Thus by applying the same argument as in \cite{T}, we can show that
the embedding of $(\mathcal{F},\mathcal{E}_{1})$ into $L^{2}(\R^{d},\beta(x)dx)$ is compact. Then by Friedrichs theorem, the spectrum of $\sigma(\mathcal{E}^{(\beta ) })$ less than $\lambda_{0}$ consists of only isolated eigenvalues with finite multiplicities. We assume $\lambda_{1}<0$, and hence $\lambda_{1}<\lambda_{0}$ is automatically true. Let $h$ be the normalized nonnegative $L^{2}$-eigenfunction corresponding to $\lambda_{1}$. It holds that $h=e^{\lambda_{1}t}P^{\beta}_{t}h$ on $E$. Since $P^{\beta}_{t}$ maps $L^{2}(\R^{d},dx)$ into $L^{\infty}(\R^{d},dx)$ and maps $\mathcal{B}_{b}(\R^{d})$ into $C(\R^{d})$, $h$ is a bounded continuous function on
$\R^{d}$.
Moreover $h$ is positive everywhere by the irreducibility of $\xi$ and the positivity of $\exp\left(\int_{0}^{t}\beta(\xi_{s})ds\right)$.
For such $(\xi,\psi_{\beta})$-superprocess, Assumptions 1-6 are satisfied with $\lambda_{1}$, $h$ defined as above and $w(x):=-\log \p_{\delta_{x}}\left(\exists t\ge 0,\ \langle 1,X_{t}\rangle =0\right)$.
\end{example}

\bigskip

\appendix
\section{Appendix}

\subsection{Kuznetsov measure for superdiffusions}\label{sup-diffusion}

Suppose the operator $\mathcal{L}$ is defined by the formula:
\begin{equation}\nonumber
\mathcal{L}:=\sum_{i,j=1}^{d}\frac{\partial}{\partial x_{j}}\,a_{ij}(x)\frac{\partial}{\partial x_{j}}+\sum_{i=1}^{d}b_{i}(x)\frac{\partial}{\partial x_{i}}\quad\mbox{ on }\mathbb{R}^{d},
\end{equation}
where the matrix $A(x):=(a_{ij}(x))_{i,j}$ is symmetric and positive definite, and all $a_{ij}$ and $b_{i}$ are bounded
and (globally) H\"older continuous on $\R^d$.
It is known (cf. \cite[Chapter 2]{D02}) that there exists a diffusion process on $\mathbb{R}^{d}$, called the $\mathcal{L}$-diffusion, whose transition density function is a fundamental solution of the equation $\partial{u}/\partial t=\mathcal{L}u$.
In this subsection, we assume $E\subset \mathbb{R}^{d}$ is a nonempty domain and $\xi$ is the subprocess of the $\mathcal{L}$-diffusion on $E$.
We assume that the branching mechanism $\psi_{\beta}(x,{\lambda})$ is given by \eqref{2.1} and
all the first partial derivatives of $\psi_{\beta}(x,{\lambda})$ are continuous.
This $(\xi,\psi_{\beta})$-superprocess is also called $(\xi,\psi_{\beta})$-superdiffusion.
Let $X$ be a $(\xi,\psi_{\beta})$-superdiffusion.
We are concerned with the set $E_{0}=\{x\in E: \mbox{Kuznetsov measure }\mathbb{N}_{x}\mbox{ exists}\}$.
By the argument in  Subsection \ref{sup-Br},
we know that \eqref{condi0} is a sufficient condition for the existence of $\mathbb{N}_{x}$.
We recall that for every open set $D\subset E$ and $t\ge 0$, there is a random measure $\widetilde{X}^{D}_{t}$
concentrated on the boundary of $[0,t)\times D$
such that \eqref{2.4} and \eqref{2.5} hold for every $\mu\in\mf$ and every $\widetilde{f}\in\mathcal{B}^{+}_{b}([0,t]\times E)$.
We also recall that
$X_{t}$ is the projection of $\widetilde{X}^{E}_{t}$ on $E$.
For every $x\in E$, let $\mathbb{O}_{x}$ be the set of open sets in $E$
that contains $x$.
Then condition \eqref{condi0} is satisfied if
\begin{equation}
\p_{\delta_{x}}\left(\widetilde{X}^{D}_{t}=0\right)>0\quad\mbox{ for all }t>0\mbox{ and }D\in \mathbb{O}_{x}.\label{condi3}
\end{equation}
 It follows from the Markov property of branching exit Markov systems that for any $D_{1},D_{2}\in\mathbb{O}_{x}$ with $D_{1}\subset D_{2}$, and any $0<s\le t<{\infty}$,
\begin{equation}
\{\widetilde{X}^{D_{1}}_{s}=0\}\subset\{\widetilde{X}^{D_{2}}_{t}=0\}\quad\p_{\delta_{x}}\mbox{-a.s.}\nonumber
\end{equation}
Thus condition \eqref{condi3} is equivalent to that:
There is $\delta>0$ sufficiently small such that
\begin{equation}
\p_{\delta_{x}}\left(\widetilde{X}^{
B_{r}
}_{t}=0\right)>0\quad\mbox{ for all }t>0\mbox{ and }0<r\le\delta.
\nonumber
\end{equation}
Here $B_{r}:=B(x,r)$ denotes the ball centered at $x$ with radius $r$.

In the remainder of this subsection, we assume
$\psi$ is  a spatially independent branching mechanism
given by \eqref{psi2}.
Suppose
$((Y_{t})_{t\ge 0};\mathbf{P}_{y},y\in\mathbb{R}^{+})$
is a one dimensional CB process with branching mechanism $\psi$ and $\mathbf{P}_{y}(Y_{0}=y)=1$. The process $Y$ is called \textit{subcritical, critical} or \textit{supercritical} according to $\psi'(0+)>$, $=$ or $<0$. It is well-known that $\mathbf{P}_y(\lim_{t\to\infty}Y_t=0)=e^{-y z_\psi},$  where $z_\psi:=\sup\{{\lambda}\ge 0, \psi({\lambda})\le 0\}\in [0,{\infty}]$.
Moreover, by the Markov property of $Y$, we have
$$\mathbf{P}_{y}\left(e^{-z_{\psi}Y_{t}}\right)=e^{-yz_{\psi}}\quad\mbox{ for all }y\ge 0\mbox{ and }t\ge 0.$$
In the subcritical and critical cases,
$z_{\psi}=0$, while in the supercritical case, $z_{\psi}\in (0,{\infty})$ if $\psi({\infty})={\infty}$,
and otherwise $z_{\psi}={\infty}$.
Given $z_{\psi}\in (0,{\infty})$,
conditioned on the event $\{\lim_{t\to\infty}Y_t=0\}$, $Y$ is a CB process with branching mechanism
\begin{equation}
\psi^*({\lambda})=\psi({\lambda}+z_\psi)=\psi'(z_\psi){\lambda}+a{\lambda}^2+\int_{(0,{\infty})}\left(e^{-{\lambda} x}-1+{\lambda} x\right)e^{-xz_\psi}\eta(dx).\label{psi*}
\end{equation}
Since $(\psi^{*})'(0+)=\psi'(z_\psi)> 0$,  the conditioned process is subcritical (cf. Lemma 2 in \cite{BKM}).

\begin{lemma}\label{lemma-supersl}
Suppose $\psi$ given by \eqref{psi2} satisfies $z_{\psi}<{\infty}$ and that
\begin{equation}\label{int-psi}\int^\infty_Nds\left[\int^s_{z_\psi}\psi(u)du\right]^{-1/2}<\infty\quad\mbox{ for some }N>z_\psi.\end{equation}
Then for any bounded open set $Q=(t_1, t_2)\times B$, where $0\le t_1<t_2<\infty$ and  $B:=B(x_0, r)$ with $x_{0}\in E$ and $r>0$, there is a function $v^0(s, x)\in C^{2}(Q)$ such that
\begin{equation}\label{super-solution-onQ}\left\{\begin{array}{rl}&\displaystyle\frac{\partial v^0}{\partial s}+Lv^0\le\psi(v^0), \quad (s,x)\in Q,\\
\\
&\displaystyle v^0(s, x)\ge z_\psi,\quad(s,x)\in Q,\\
\\
&\displaystyle\lim_{ Q\ni(s,x)\to z}v^0(s, x)=\infty,\quad z\in \mathcal{T}_{Q},\end{array}\right.\end{equation}
where $\mathcal{T}_{Q}:=((t_{1},t_{2})\times \partial B)\cup(\{t_2\}\times \overline{B})$
is a total set of $\partial Q$.
\end{lemma}

\proof The idea of this proof is similar to that of \cite[Theorem 5.3.1]{D02}. We only need to prove the result for sufficiently small $t_2-t_1+r$.
Note that \eqref{int-psi} is equivalent to
$$\int^\infty_Nds\left[\int^s_{0}\psi^*(r)dr\right]^{-1/2}<\infty\quad\mbox{ for some }N>0,$$
where $\psi^{*}$ is defined in \eqref{psi*}.
By \cite[Lemma 5.3.1]{D02},
for sufficiently small $t_2-t_1$ there is a nonnegative
solution
$u^*(s)$ of the following problem:
$$\left\{\begin{array}{rl}&\displaystyle\frac{\partial u^*}{\partial s}\le\psi^*(u^*), \quad s\in (t_{1},t_{2}),\\
&\displaystyle
\lim_{(t_{1},t_{2})\ni s\to t_{2}}
u^*(s)=\infty.\end{array}\right.$$
Take $u(s):=z_\psi+u^*(s)$ for $s\in(t_1, t_2)$.
Using the fact that $\psi^*({\lambda})=\psi({\lambda}+z_\psi)$, it is easy to check that $u$ satisfies \eqref{supersolution}:
\begin{equation}\label{supersolution}
\left\{\begin{array}{rl}&\displaystyle\frac{\partial u}{\partial s}\le\psi(u), \quad s\in (t_1,t_2),\\
\\
&\displaystyle u(s)\ge z_\psi,\quad s\in (t_1,t_2).\\
\\
&\displaystyle
\lim_{(t_{1},t_{2})\ni s\to t_2}
u(s)=\infty,
\end{array}\right.\end{equation}
On the other hand, by \cite[Lemma 5.3.4]{D02} and the fact that $\psi^*({\lambda})=\psi({\lambda}+z_\psi)$, for sufficiently small $r$, there is a nonnegative solution $v(x)$ of the following problem:
$$\left\{\begin{array}{rl}&\displaystyle Lv\le\psi(v), \quad x\in B,\\
&\displaystyle
\lim_{B\ni x\to z} v(x)=\infty, \quad z\in\partial B.
\end{array}\right.$$
Hence $v^0(s,x):=u(s)+v(x)$ is a solution of \eqref{super-solution-onQ}.\qed

\begin{lemma} \label{CP} [Comparison Principle] Suppose $\psi$  is given by \eqref{psi2} with $z_{\psi}<{\infty}$ and that $Q$ is a bounded regular open set in $[0,\infty)\times \mathbb{R}^d$. Then $u\le v$ in $Q$ assuming that

(a) $u, v\in C^2(Q)$;

(b)
$
\frac{\partial u}{\partial s}+Lu-\psi(u)\ge \frac{\partial v}{\partial s}+Lv-\psi(v)\quad \mbox{ in } Q;
$

(c) $u$ is bounded from above and $v\ge z_\psi$ in $Q$.

(d) for every $\tilde{z}\in \mathcal{T}_{Q},$ $\limsup_{Q\ni z\to \tilde{z}}[u(z)-v(z)]\le 0.$
\end{lemma}

\proof The proof is similar to that of \cite[Theorem 5.2.3]{D02} by using the fact that $\psi({\lambda})$ is an increasing function in ${\lambda}\in [z_\psi, \infty)$.
We omit the details here.
\qed

\begin{lemma}\label{lem2.8}
 If $x\in E$ satisfies the following condition:
 there exist an open set $D\in\mathbb{O}_{x}$ and a function $\psi$ in the form of \eqref{psi2} with $z_{\psi}<{\infty}$ such that \eqref{int-psi} holds and
 \begin{equation}
 \psi_\beta(x, {\lambda})\ge \psi({\lambda})\quad\mbox{ for all }x\in D \mbox{ and }{\lambda}\ge 0,\label{0.22}
 \end{equation}
then $\p_{\delta_{x}}\left(X_{t}=0\right)>0$ for all $t>0$ and hence $x\in E_{0}$.
 \end{lemma}

\proof   By the argument in the beginning of this subsection, it suffices to prove that for $\delta>0$ small enough such that $B_{\delta}:=B(x,\delta)\subset D$, we have
\begin{equation}
\p_{\delta_{x}}\left(\widetilde{X}^{
B_{r}
}_{t}=0\right)>0\quad\mbox{ for all }t>0\mbox{ and }
B_{r}
:=B(x,r)\mbox{ with }r\in (0,\delta].\label{2.17}
\end{equation}
Fix $r\in (0,\delta]$ and an arbitrary $T>0$. For any $\lambda>0$, let
\begin{equation}
\widetilde{u}^{B_{r}}_{\lambda}(s,y):=-\log\p_{\delta_{y}}\left(\exp\left(-\langle \lambda,\widetilde{X}^{B_{r}}_{s}\rangle\right)\right)\quad\mbox{ for }y\in B_{r}\mbox{ and }s\in [0,T].\label{2.18}
\end{equation}
Since all the first partials of $\psi_{\beta}(x,\lambda)$ are continuous, it follows by \cite[Theorem 5.2.2]{D02} that $v_{\lambda}(s,y):=\widetilde{u}^{B_{r}}_{\lambda}(T-s,y)$ for $(s,y)\in (0,T)\times B_{r}$ is a solution of  the following boundary problem:
$$\left\{\begin{array}{rl}&\displaystyle\frac{\partial v_{\lambda}}{\partial s}+Lv_{\lambda}=\psi_\beta(y,v_{\lambda}(s,y)), \quad (s,y)\in(0,T)\times B_{r},\\
&\displaystyle\lim_{(0,T)\times B_{r}\ni (s,y)\to z}v_{\lambda}(s,y)=\lambda,\quad z\in\partial ((0,T)\times B_{r}).\end{array}\right.$$
Since $\psi_\beta(y,\lambda)\ge\psi(\lambda)$ for
$y\in B_{r}$
and $\lambda\ge 0$, we have
$$\left\{\begin{array}{rl}&\displaystyle\frac{\partial v_{\lambda}}{\partial s}+Lv_{\lambda}\ge\psi(v_{\lambda}), \quad (s,y)\in(0,T)\times B_{r},\\
&\displaystyle\lim_{(0,T)\times B_{r}\ni (s,y)\to z}v_{\lambda}(s,y)=\lambda,\quad z\in\partial ((0,T)\times B_{r}).\end{array}\right.$$
By Lemma \ref{lemma-supersl}, there is a solution $v^0(s,y)$ of the following problem:
$$\left\{\begin{array}{rl}&\displaystyle\frac{\partial v^0}{\partial s}+Lv\le\psi(v^0), \quad (s,y)\in(0,T)\times B_{r},\\
\\
&\displaystyle v^{0}(s,x)\ge z_{\psi},\quad (s,y)\in(0,T)\times B_{r},\\
\\
&\displaystyle\lim_{ (0,T)\times B_{r} \ni(s,y)\to z}v^0(s, y)={\infty},\quad z \in \mathcal{T}_{(0,T)\times B_{r}}.\end{array}\right.$$
where $\mathcal{T}_{(0,T)\times B_{r}}:=((0,T)\times \partial B_{r})\cup(\{T\}\times \overline{B_{r}})$
is a total set of $\partial((0,T)\times B_{r})$.
By Lemma \ref{CP}, for any $\lambda>0$, $$v_{\lambda}(s, y)\le v^0(s,y)\quad\mbox{ for all } (s,y)\in (0,T)\times B_{r}.$$
Letting $\lambda\uparrow\infty$ in \eqref{2.18}, we get that for any $t\in (0,T)$
$$-\log \p_{\delta_{x}}(\widetilde{X}^{B_{r}}_{t}=0)=\lim_{\lambda\to\infty}v_{\lambda}(T-t, x)\le v^0(T-t, x)<{\infty}.$$
Since $T>0$ is arbitrary, we get \eqref{2.17} for all $t>0$ and $r\in (0,\delta]$.\qed

\medskip

Note that condition \eqref{int-psi} is stronger than $\int^{{\infty}}\psi({\lambda})^{-1}d{\lambda}<{\infty}$, which is usually called
 Grey's condition.
Grey's condition
is a necessary and sufficient condition
 for a CB process to become extinct in a finite time with positive probability.
 Lemma \ref{lem2.8}
says that for a superdiffusion with spatially dependent branching mechanism, if the branching mechanism is locally dominated by a spatially independent branching mechanism which satisfies a condition stronger than the
 Grey's condition,
then this superdiffusion
become extinct in a finite time with positive probability.

 Recall that $\alpha (x)$ is the non-negative bounded Borel measurable function in \eqref{2.1} and $E_{+}=\{x\in E:\ \alpha(x)>0\}$.

\begin{proposition} \label{prop2.9}
 For a $(\xi,\psi_{\beta})$-superdiffusion, it holds that
 $E_{+} \subset E_{0}$.
\end{proposition}

\proof
Since all the first partials of $\psi_{\beta}$ are continuous,
$\alpha$ and $\beta$ in $\psi_{\beta}$ are continuous functions.
Thus for
any $x\in E_{+}$, there exists a neighborhood $D$ of $x$ such that $\sup_{y\in D}\beta(y)\le \beta^{*}<{\infty}$ and $\inf_{y\in D}\alpha(y)\ge \alpha^{*}>0$.
We have
\begin{equation}\nonumber
\psi_{\beta}(y,\lambda)\ge -\beta^{*}\lambda+\alpha^{*}\lambda^{2}=:\psi^{*}(\lambda)\quad\mbox{ for }y\in D\mbox{ and }\lambda\ge 0.
\end{equation}
It is easy to verify that $\psi^{*}$ satisfies \eqref{int-psi}.
Hence $x\in E_{0}$ by
Lemma \ref{lem2.8}.
\qed

In the rest of this subsection, we assume Condition 1' holds. We will establish the existence of $\mathbb{N}^{*}_{x}$ for all $x\in E_{+}$ under Condition 1'.

First we note that there exists a superdiffusion
$((X_{t})_{t\ge 0};\p^{*}_{\mu},\mu\in\mf)$
with $\p^*_{\mu}(X_{0}=\mu)=1$ such that for all $\mu\in\mf$, $f\in\mathcal{B}^{+}_{b}(E)$ and $t\ge 0$,
    \begin{equation}
    \p^*_{\mu}\left[e^{-\langle f,X_{t}\rangle}\right]=e^{-\langle u^{*}_{f}(t,\cdot),\mu\rangle},
    \nonumber
    \end{equation}
where $u^{*}_{f}(t,x)$ is the unique nonnegative locally bounded solution to the integral equation given by \eqref{1.8}.
The branching mechanism of this superdiffusion is $\psi^{*}_{\beta^{*}}(x,\lambda):=-\beta^{*}(x)\lambda+\psi^{*}_{0}(x,\lambda)$ where $\beta^{*}$ and $\psi^{*}_{0}$ are given in Proposition \ref{prop2.4}. Since $w$ is only locally bounded on $E$, $\beta^{*}$ is bounded from above
but may not be bounded from below. Hence the branching mechanism $\psi^{*}_{\beta^{*}}(x,\lambda)$ does not satisfy the usual assumptions in Section \ref{sec2.1}. Nevertheless, one may consider the $(\xi^{D},\psi^{*}_{\beta^{*}})$-superdiffusion in every bounded open domain $D$, where the underlying spatial motion is killed upon hitting
 the boundary of $D$. Then by using an increasing sequence of compactly embedded domains to approximate $E$, the process $(X,\p^{*}_{\mu})$ can be defined as a distributional limit of these $(\xi^{D},\psi^{*}_{\beta^{*}})$-superdiffusions. (see, for example, the argument before Theorem 2.4 in \cite{EKW}.)   We
remark  that this method of construction may fail in general when the spatial motion has discontinuous sample paths
since the process can potentially jump everywhere when it exits from the domain $D$.

   \medskip

Let $E^{*}_{0}$ be the set of points in $E$ where the Kuznetsov measure $\mathbb{N}^{*}_{x}$ corresponding to the $(\xi,\psi^{*}_{\beta^{*}})$-superdiffusion exists.
It follows that for every $x\in E^{*}_{0}$, every $t>0$ and $f\in\mathcal{B}^{+}_{b}(E)$,
\begin{equation}
\mathbb{N}^{*}_{x}\left(1-e^{-\langle f,X_{t}\rangle}\right)
=u^{*}_{f}(t,x).\nonumber
\end{equation}

According to \cite[Theorem I.1.1 and Theorem I.1.2]{D93},
for any bounded open set $D\Subset E$
and $t\ge 0$, there exists a random measure $\widetilde{X}^{*D}_{t}$
concentrated on the boundary of $[0,t)\times D$,
such that for every $\mu\in\mf$ and $\widetilde{f}\in\mathcal{B}^{+}_{b}([0,t]\times E)$,
 \begin{equation}
 \p^*_{\mu}\left[e^{-\langle\widetilde{f},\widetilde{X}^{*D}_{t}\rangle}\right]=e^{-\langle\widetilde{u}^{*D}_{\widetilde{f}}(t, \cdot),\mu\rangle},
 \nonumber
 \end{equation}
 where $\widetilde{u}^{*D}_{\widetilde{f}}(t,x)$ is the unique nonnegative locally bounded solution to the following integral equation: \begin{equation}
 \widetilde{u}^{*D}_{\widetilde{f}}(t, x)=\Pi_{x}\left[\widetilde{f}(t\wedge \tau_{D}, \xi_{t\wedge \tau_{D}})\right]-
 \Pi_{x}\left[\int_{0}^{t\wedge \tau_{D}}\psi^*_{\beta^*}\left(\xi_{s},\widetilde{u}^{*D}_{\widetilde{f}}(t-s, \xi_{s})\right)ds\right].\label{2.5'}
 \end{equation}
 Define $\widetilde{w}(t,x):=w(x)$ for $(t,x)\in [0,{\infty})\times E$.
 Using the local boundedness of $w$ on $E$, we can deduce that
 for any bounded open set $D\Subset E$,
 \begin{equation}\label{1.31}
 \widetilde{w}(t,x)=\Pi_{x}\left(\widetilde{w}(t\wedge\tau_{D},\xi_{t\wedge \tau_{D}})\right)-\Pi_{x}\left(\int_{0}^{t\wedge\tau_{D}}\psi_{\beta}(\xi_{s},\widetilde{w}(t-s,\xi_{s}))ds\right).
 \end{equation}

The next proposition shows that $\{\widetilde{X}^{*D}_{t};t\ge 0\}$ is a Doob's $h$-transformed process of $\{\widetilde{X}^{D}_{t};t\ge 0\}$ via the function $\widetilde{w}$.
\begin{prop}\label{prop2.10}
Let $\mw:=\{\mu\in\mf:\langle w,\mu\rangle<{\infty}\}$.
Suppose $D\Subset E$ is a bounded open set
and $\mu\in\mw$.
Then
for every $t\ge 0$ and $\widetilde{f}\in\mathcal{B}^{+}_{b}([0,t]\times E)$,
\begin{equation}\label{martingale-change}
e^{-\langle w,\mu\rangle} \p_{\mu}\left[e^{-\langle\widetilde{f}+\widetilde{w},\widetilde{X}^{D}_{t}\rangle}\right]= \p^*_{\mu}\left[e^{-\langle\widetilde{f},\widetilde{X}^{*D}_{t}\rangle}\right].
\end{equation}
\end{prop}

\proof Since the random measure $\widetilde{X}^{D}_{t}$ is supported on the boundary of $[0,t)\times D$ and $\widetilde{w}$ is locally bounded on $[0,{\infty})\times E$,
it follows by \eqref{2.4} and \eqref{2.5} that
$$e^{-\langle w,\mu\rangle} \p_{\mu}\left[\exp\left(-\langle\widetilde{f}+\widetilde{w},\widetilde{X}^{D}_{t}\rangle\right)\right]=\exp\left(-\langle \widetilde{u}^{D}_{\widetilde{w}+\widetilde{f}}(t,\cdot)-w(\cdot),\mu\rangle\right),$$
where $\widetilde{u}^{D}_{\widetilde{w}+\widetilde{f}}(t,x)$ is the unique nonnegative locally bounded solution to the following integral equation
\begin{equation}\label{1.35}
\widetilde{u}^{D}_{\widetilde{w}+\widetilde{f}}(t,x)=\Pi_{x}\left[(\widetilde{w}+\widetilde{f})(t\wedge\tau_{D},\xi_{t\wedge \tau_{D}})\right]-\Pi_{x}\left(\int_{0}^{t\wedge\tau_{D}}\psi_{\beta}(\xi_{s},\widetilde{u}^{D}_{\widetilde{w}+\widetilde{f}}(t-s,\xi_{s}))ds\right)
\end{equation}
Using \eqref{1.35} and \eqref{1.31}, it is straightforward to check that $\widetilde u^{*D}_{\widetilde f}(t,x):=\widetilde u^D_{\widetilde{f}+\widetilde w}(t,x)-w(x)$ is the unique nonnegative locally bounded solution to \eqref{2.5'}. Thus we get \eqref{martingale-change}.\qed

 \begin{lemma}\label{lem2.11}
 If $x\in E$ satisfies the condition of
   Lemma \ref{lem2.8},
 then $x\in E^{*}_{0}$.
\end{lemma}

\proof
Based on the argument in the beginning of this subsection, it suffices
to prove that there exists $\delta>0$ such that
\begin{equation}
\nonumber
\p^*_{\delta_x}(\widetilde{X}^{*B}_{t}=0)>0 \quad \mbox{ for all } t>0\mbox{ and }B=B(x,r)\mbox{ with }r\le\delta.\end{equation}
We choose $\delta$ sufficiently small such that
$B\subset D$, and $\psi_\beta(y, \lambda)\ge \psi(\lambda)$ for all $y\in B$ and $\lambda\ge 0$, where $D$ and $\psi$ satisfy the conditions of
 Lemma \ref{lem2.8}.
It follows from Proposition \ref{prop2.10} that
$$\begin{array}{rll}\displaystyle\p^*_{\delta_x}(\widetilde{X}^{*B}_{t}=0)&=&\displaystyle \lim_{\lambda\to{\infty}}\p_{\delta_x}(e^{-\langle\lambda, \widetilde{X}^{*B}_{t}\rangle})\\
&=&\displaystyle e^{-w(x)}\lim_{\lambda\to{\infty}}\p_{\delta_x}(e^{-\langle\widetilde w+\lambda, \widetilde{X}^{B}_{t}\rangle})\\
&=&\displaystyle e^{-w(x)}\lim_{\lambda\to\infty}e^{-\widetilde u^B_{\widetilde w+\lambda}(t,x)},\end{array}$$
where $\widetilde u^B_{\widetilde w+\lambda}(t,x)$ is the unique nonnegative solution of
\eqref{2.5}
with initial condition $\widetilde w+\lambda$.
Then $v(s, y):=\widetilde{u}^{B}_{\widetilde w+\lambda}(t-s, y)$ for $(s,y)\in[0,t]\times B,$ is a solution to the following boundary problem:
$$\left\{\begin{array}{rl}&\displaystyle\frac{\partial v}{\partial s}+Lv=\psi_\beta(y,v(s,y)), \quad (s,y)\in(0,t)\times B\\
&\displaystyle\lim_{(0,t)\times B\ni (s,y)\to z}v(t, y)=\widetilde{w}(z)+\lambda,\quad z\in\partial ((0,t)\times B).\end{array}\right.$$
Since $\psi_\beta(y,\lambda)\ge\psi(\lambda)$  for $y\in B$ and $\lambda\ge 0$, applying similar arguments as in the proof of
 Lemma \ref{lem2.8},
we get $\lim_{\lambda\to{\infty}}\widetilde u^B_{\widetilde w+\lambda}(t,x)<{\infty}$. Therefore $\p^*_{\delta_x}(\widetilde{X}^{*B}_{t}=0)>0$ for $t>0$. Hence $x\in E^{*}_{0}$.\qed

\begin{proposition}\label{prop2.12}
For a $(\xi,\psi_{\beta})$-superdiffusion that
satisfies Condition 1', it holds that $E_{+}\subset E^{*}_{0}$.
\end{proposition}

\proof This proposition follows from
Lemma \ref{lem2.11} (instead of Lemma \ref{lem2.8})
in the same way as the proof of Proposition \ref{prop2.9}. We omit the details here.\qed

\subsection{Verifications of examples}\label{verification}

\noindent\textbf{Example \ref{E:8.1}:}
It is easy to check that Assumptions 1-4 are
satisfied with $\lambda_{1}=-\beta$, $h(x)\equiv 1$ and
$w(x)\equiv
z_{\psi}$. Then the $h$-transformed semigroup $P^{h}_{t}$ is given by
$$P^{h}_{t}f(x):=\frac{e^{-\beta t}}{h(x)}\,\Pi_{x}\left[e^{\beta t}h(\xi_{t})f(\xi_{t})\right]=P_{t}f(x)\quad\mbox{ for }f\in\mathcal{B}^{+}_{b}(\R^{d}).$$
This implies that
the $h$-transformed process $\xi^{h}$ is still an inward OU process with generate $\mathcal{L}$. It is known that the transition density of $\xi^{h}$ with respect to $m$ is given by
$$p^{h}(t,x,y)=\left(1-e^{-2ct}\right)^{-d/2}\exp\left(-\frac{c}{e^{2ct}-1}\left(|x|^{2}+|y|^{2}-2e^{ct}x\cdot y\right)\right)\quad\mbox{ for }t>0\mbox{ and }x,y\in\R^{d}.$$
In view of this, one can easily check that $\xi^{h}$ has the Feller property, that is, $P^{h}_{t}$ maps $C_{0}(\R^{d})$ to $C_{0}(\R^{d})$ and $\lim_{t\to 0}\|P^{h}_{t}f-f\|_{\infty}=0$ for all $f\in C_{0}(\R^{d})$.
It remains to show that Assumption 5 holds for this example.
Let $Z$ be the skeleton process. We need to show that
for any $\mu\in \mc$, $\sigma>0$ and $f\in\mathcal{B}^{+}(E)$ with $\frac{fw}{h}$ bounded,
\begin{equation}
\lim_{n\to{\infty}}e^{\lambda_{1}n\sigma}\langle f,Z_{n\sigma}\rangle =\langle f,wh\rangle W^{h/w}_{\infty}(Z)\quad\pp_{\mu}
\mbox{-a.s.}
\label{7.20}
\end{equation}
For $t\ge 0$, $Z_{t}=\sum_{i=1}^{N_{0}}Z^{i,0}_{t}$, where $N_{0}=\langle 1,Z_{0}\rangle$ and $Z^{i,0}$ denotes the independent subtree of the skeleton initiated by the $i$th particle at time $0$. Recall that under $\pp_{\mu}$ for $\mu\in\mc$, $N_{0}$ is a Poisson random variable with mean $\langle w,\mu\rangle$. Moreover, given $Z_{0}$, $Z^{i,0}$ follows the same distribution as $(Z;\pp_{\cdot,\delta_{z_{i}(0)}})$. Thus we have
\begin{eqnarray}
&&\pp_{\mu}\left(\lim_{n\to{\infty}}e^{\lambda_{1}n\sigma}\langle f,Z_{n\sigma}\rangle=\langle f,wh\rangle W^{h/w}_{\infty}(X)\right)\nonumber\\
&\ge&\pp_{\mu}\left(\bigcap_{i=1}^{N_{0}}\{\lim_{n\to{\infty}}e^{\lambda_{1}n\sigma}\langle f,Z^{i,0}_{n\sigma}\rangle=\langle f,wh\rangle W^{h/w}_{\infty}(Z^{i,0})\}\right)\nonumber\\
&=&\pp_{\mu}\left[\prod_{i=1}^{N_{0}}\pp_{\cdot,\delta_{z_{i}(0)}}\left(\lim_{n\to{\infty}}e^{\lambda_{1}n\sigma}\langle f,Z_{n\sigma}\rangle=\langle f,wh\rangle W^{h/w}_{\infty}(Z)\right)\right].\nonumber
\end{eqnarray}
Hence to prove \eqref{7.20}, it suffices to prove that for every $x\in E $,
\begin{equation}
\pp_{\cdot,\delta_{x}}\left(\lim_{n\to{\infty}}e^{\lambda_{1}n\sigma}\langle f,Z_{n\sigma}\rangle=\langle f,wh\rangle W^{h/w}_{\infty}(Z)\right)=1.\label{7.21}
\end{equation}
Sufficient conditions for \eqref{7.21} are given in \cite{EKW} (see conditions (2.23)-(2.26), (2.34) and (2.35) there).
In this example, conditions (2.23)-(2.26) in \cite{EKW} hold for $\sigma_{1}=\sigma_{3}=p$, $\sigma_{2}=2$ and $\varphi_{1}(x)=\varphi_{2}(x)\equiv 1$, condition (2.34) holds for $a(t)=\sqrt{(\frac{-\lambda_{1}}{c}+\sigma)t}$ and condition (2.35) in \cite{EKW} holds for $K=1$.
Hence we prove \eqref{7.20}.

\medskip

\noindent\textbf{Example \ref{E:8.2}:}
Assumptions 1-4 are satisfied with $\lambda_{1}=cd-\beta$, $h(x)=\left(\frac{c}{\pi}\right)^{d/2}e^{-c|x|^{2}}$ and $w(x)\equiv z_{\psi}$. The transition density of $\xi$ with respect to $m$ is given by
$$p(t,x,y)=\left(e^{2ct}-1\right)^{-d/2}\exp\left(-\frac{c}{1-e^{-2ct}}\left(|x|^{2}+|y|^{2}-2e^{-ct}x\cdot y\right)\right)\quad\mbox{ for }t>0\mbox{ and }x,y\in\R^{d}.$$
The $h$-transformed process $\xi^{h}$ is symmetric with respect to the measure $\widetilde{m}(dx)=h^{2}m(dx)=\left(\frac{c}{\pi}\right)^{d/2}e^{-c|x|^{2}}dx$, and the transition density of $\xi^{h}$ with respect to $\widetilde{m}$ is given by
$$p^{h}(t,x,y)=\left(1-e^{-2ct}\right)^{-d/2}\exp\left(-\frac{c}{e^{2ct}-1}\left(|x|^{2}+|y|^{2}-2e^{ct}x\cdot y\right)\right)\quad\mbox{ for }t>0\mbox{ and }x,y\in\R^{d}.$$
This implies that $\xi^{h}$ is an inward OU process on $\R^{d}$ with infinitesimal generator $\mathcal{L}^{h}:=\frac{1}{2}\Delta-cx\cdot\nabla$. Thus $\lim_{t\to 0}\|P^{h}_{t}f-f\|_{\infty}=0$ for any $f\in C_{0}(\mathbb{R}^{d})$.
In view of Example 4.2 in \cite{EKW}, Assumption 5 is also satisfied.

\medskip

\noindent\textbf{Example \ref{E:8.3}:} Suppose $\widehat{\psi}_{\beta}(x,\lambda)  :=
-\beta(x)\lambda+\alpha(x)\lambda^{2}$.
Let $\mathbb{X}$ be a $(Y,\psi_{\beta})$-superdiffusion and $\widehat{\mathbb{X}}$ be a $(Y,\widehat{\psi}_{\beta})$-superdiffusion. For $x\in \R^{d}$, let $w_{\mathbb{X}}(x):=-\log\p_{\delta_{x}}\left(\exists t\ge 0: \langle 1,\mathbb{X}_{t}\rangle =0\right)$ and $w_{\widehat{\mathbb{X}}}(x):=-\log\p_{\delta_{x}}\left(\exists t\ge 0: \langle 1,\widehat{\mathbb{X}}_{t}\rangle =0\right)$. We know from \cite{EP} that the function $w_{\widehat{\mathbb{X}}}$ is continuous on $\R^{d}$ and solves the equation $\mathcal{L}u-\widehat{\psi}_{\beta}(x,u)=0$ on $\R^{d}$. Since $\psi_{\beta}\ge \widehat{\psi}_{\beta}$ pointwise, it follows by \cite[Lemma 4.5]{EKW} that $w_{\mathbb{X}}\le w_{\widehat{\mathbb{X}}}$ pointwise. Thus $w_{\mathbb{X}}$ is locally bounded on $\R^{d}$ and hence is continuous on $\R^{d}$ by \cite[Lemma 2.1]{EKW}. Recall that $X$ is a
$(\xi,\psi_{\beta})$-superdiffusion where $\xi=Y^{E}$ is the subprocess of $Y$ killing upon leaving $E$.
We mentioned in Section \ref{sec2.2} that one may think of $X_{t}$ describing the mass in $\mathbb{X}_{t}$ which historically avoids exiting $E$.

Since the function $\rho$ in the generator is bounded between two positive constants, we know from \cite{KS} that $\xi$ has a positive continuous transition density $p_{E}(t,x,y)$ with respect to $m$, and for each $T>0$, there exist positive $c_{i}$, $i=1,\cdots,4$ such that for every $(t,x,y)\in (0,T]\times E\times E$,
\begin{equation}
c_{1}f_{E}(t,x,y)
t^{-d/2}\exp\left(-\frac{c_{2}|x-y|^{2}}{t}\right)\le p_{E}(t,x,y)\le c_{3}f_{E}(t,x,y)
t^{-d/2}\exp\left(-\frac{c_{4}|x-y|^{2}}{t}\right),\label{8.17}
\end{equation}
where $f_{E}(t,x,y):=\left(1\wedge \frac{\delta_{E}(x)}{\sqrt{t}}\right)\left(1\wedge \frac{\delta_{E}(y)}{\sqrt{t}}\right)$ and $\delta_{E}(x)$ denotes the Euclidean distance between $x$ and the boundary of $E$.
Since $\beta$ is bounded,
it follows that the Feynman-Kac semigroup
$P^{\beta}_{t}$ admits an integral kernel $p_{E}^{\beta}(t,x,y)$ which is positive, symmetric and continuous in $(x,y)$ for each $t>0$, and
\begin{equation}
e^{-\|\beta\|_{\infty}t}p_{E}(t,x,y)\le p_{E}^{\beta}(t,x,y)\le e^{\|\beta\|_{\infty}t}p_{E}(t,x,y)\quad\mbox{ for all } x,y\in E\mbox{ and }t>0.\nonumber
\end{equation}
Thus $p^{\beta}_{E}(t,x,y)$ satisfies the same two-sided estimates \eqref{8.17} with possibly different constants $c_{i}>0$, $i=4,\cdots,8$. By this estimate, $\int_{E\times E}p_{E}^{\beta}(t,x,y)^{2}m(dx)m(dy)=\int_{E}p_{E}^{\beta}(2t,x,x)m(dx)<{\infty}$ for every $t\in (0,T]$. Thus $P^{\beta}_{t}$ is a Hilbert-Schimidt operator in $L^{2}(E,m)$ and hence is compact. The infinitesimal generator of $P^{\beta}_{t}$ is $\mathcal{L}^{(\beta)}:=\left(\mathcal{L}+\beta\right)|_{E}$ with zero Dirichlet boundary condition. We use $\sigma(\mathcal{L}^{\beta})$ to denote the spectrum set of $\mathcal{L}^{\beta}$. It then follows from Jentzch's theorem that
$\lambda_{1}:=\inf\{-\lambda:\ \lambda\in\sigma(\mathcal{L}^{(\beta)})\}$
is a simple eigenvalue and a corresponding eigenfunction $h$ can be chosen to be nonnegative with $\int_{E}h(x)^{2}m(dx)=1$.  We assume $\lambda_{1}<0$. Since $h(x)=e^{\lambda_{1}t}\int_{E}p_{E}^{\beta}(t,x,y)h(y)m(dy)$, $h$ is continuous and positive on $E$. Moreover, by the estimate \eqref{8.17}, we have
\begin{equation}
c_{9}\left(1\wedge \delta_{E}(x)\right)\le h(x)\le c_{10}\left(1\wedge \delta_{E}(x)\right) \quad\mbox{for all } x\in E.\label{8.18}
\end{equation}
Here $c_{9},\ c_{10}$ are positive constants independent of $x$.
Therefore, Assumptions 2-3 hold.

Let $w(x):=-\log\mathrm{P}_{\delta_{x}}\left(\exists t\ge 0:\ \langle 1,X_{t}\rangle =0\right)$.
Since $h$ is bounded on $E$, it is easy to see that
condition \eqref{llogl1} (or, equivalently, Assumption 4(i))
holds.
Thus by Theorem \ref{them2} and Remark \ref{rm for them2},
$\p_{\delta_{x}}\left(W^{h}_{\infty}(X)>0\right)>0$ for every $x\in E$.
It follows that
$$
w(x)\ge -\log\left(1-\p_{\delta_{x}}\left(\limsup_{t\to{\infty}}\langle 1,X_{t}\rangle >0\right)\right)\ge -\log\left(1-\p_{\delta_{x}}\left(W^{h}_{\infty}(X) >0\right)\right)>0\quad\mbox{ for }x\in E.$$
On the other hand, since $X_{t}$ describes only part of the mass in $\mathbb{X}_{t}$, we have $\mathrm{P}_{\delta_{x}}\left(\exists t\ge 0:\ \langle 1,X_{t}\rangle =0\right)\ge \mathrm{P}_{\delta_{x}}\left(\exists t\ge 0:\ \langle 1,\mathbb{X}_{t}\rangle =0\right)$ and hence $w(x)\le w_{\mathbb{X}}(x)$ on $E$. This together with the continuity of $ w_{\mathbb{X}}$ on $\R^{d}$ implies that $w$ is bounded on $E$ and hence is continuous on $E$ by \cite[Lemma 2.1]{EKW}.
This shows that Condition 1  holds,
and, consequently, Assumption 1 is satisfied.

We now show that Assumption 4(ii) holds for this example.
By condition \eqref{8.22} and the boundedness of $w$, one can easily prove that $\frac{\psi_{\beta}(x,w(x))}{w(x)}=-\beta(x)+\alpha(x)w(x)+\int_{(0,{\infty})}\frac{e^{-w(x)y}-1+w(x)y}{w(x)}\,\pi(x,dy)$ is a bounded function on $E$. Let $\{D_{n}:n\ge 1\}$ be a sequence of bounded domains with smooth boundaries such that $D_{n}\Subset D_{n+1}\Subset E$ for $n\ge 1$ and $\bigcup_{n=1}^{{\infty}}D_{n}=E$. We know from the argument in the beginning of Section \ref{sec5} that for $t>0$, $x\in E$ and $n$ sufficiently large so that $x\in D_{n}$,
\begin{eqnarray}
w(x)&=&\Pi_{x}\left[w(\xi_{t\wedge \tau_{D_{n}}})\exp\left(-\int_{0}^{t\wedge \tau_{D_{n}}}\frac{\psi_{\beta}(\xi_{s},w(\xi_{s}))}{w(\xi_{s})}\,ds\right)\right].\nonumber\end{eqnarray}
Consequently,
\begin{eqnarray}
w(x)&\ge&\exp\left(-t\sup_{y\in E}\frac{\psi_{\beta}(y,w(y))}{w(y)}\right)\Pi_{x}\left(w(\xi_{t\wedge \tau_{D_{n}}})\right)\nonumber\\
&=&e^{-c_{13}t}\left(\Pi_{x}\left[w(\xi_{t});t<\tau_{D_{n}}\right]+\Pi_{x}\left[w(\xi_{\tau_{D_{n}}});t\ge\tau_{D_{n}}\right]\right)\nonumber\\
&\ge&e^{-c_{13}t}\Pi_{x}\left[w(\xi_{t});t<\tau_{D_{n}}\right]=e^{-c_{13}t}\Pi_{x}\left[w(Y_{t});t<\tau_{D_{n}}\right].\nonumber
\end{eqnarray}
Thus by letting $n\to{\infty}$, we get
$w(x)\ge e^{-c_{13}t}\Pi_{x}\left[w(Y_{t});t<\tau_{E}\right]=e^{-c_{13}t}\Pi_{x}\left[w(Y^{E}_{t})\right]$.
Using this and the heat kernel estimates in \eqref{8.17}, we have for $x\in E$
\begin{eqnarray}
w(x)&\ge&e^{-c_{13}}\Pi_{x}\left[w(Y^{E}_{1})\right]\nonumber\\
&=&e^{-c_{13}}\int_{E}p_{E}(1,x,y)w(y)m(dy)\nonumber\\
&\ge&c_{1}e^{-c_{13}}\int_{E}(1\wedge \delta_{E}(x))(1\wedge \delta_{E}(y))e^{-c_{2}|x-y|^{2}}w(y)m(dy)\nonumber\\
&\ge&c_{1}e^{-c_{13}-c_{2}\mathrm{diam}(E)^{2}}(1\wedge \delta_{E}(x))\int_{E}(1\wedge \delta_{E}(y))w(y)m(dy)\nonumber\\
&\ge&c_{14}(1\wedge \delta_{E}(x)).\label{8.23}
\end{eqnarray}
Here $c_{14}>0$ is  constant and the last inequality comes from
the fact that $(1\wedge \delta_{E}(y))w(y)$ is positive everywhere on $E$.
By \eqref{8.23} and \eqref{8.18}, we have for $x\in E$,
\begin{eqnarray}
\int_{(0,{\infty})}r^{2}e^{-w(x)r}\pi(x,dr)&\le&\int_{(0,{\infty})}r^{2}e^{-c_{14}(1\wedge \delta_{E}(x))r}\pi(x,dr)\nonumber\\
&\le&\int_{(0,{\infty})}r^{2}e^{-c_{15}h(x)r}\pi(x,dr)\nonumber\\
&=&\frac{1}{h(x)}\int_{(0,{\infty})}r\log^{*}r\left(\frac{rh(x)}{\log^{*}(rh(x))}e^{-c_{15}h(x)r}\right)\frac{\log^{*}(rh(x))}{\log^{*}r}\,\pi(x,dr)\nonumber\\
&\le&\frac{1}{h(x)}\int_{(0,{\infty})}r\log^{*}r\left(\frac{rh(x)}{\log^{*}(rh(x))}e^{-c_{15}h(x)r}\right)
\frac{\log^{*}(r\|h\|_{\infty})}{\log^{*}r}\,\pi(x,dr),\nonumber
\end{eqnarray}
where $c_{15}$ is a positive constant.
It then follows from condition \eqref{8.22}, and the fact that functions $y\mapsto \frac{y}{\log^{*}y}e^{-c_{15}y}$ and $y\mapsto \frac{\log^{*}(y\|h\|_{\infty})}{\log^{*}y}$ are bounded on $(0,{\infty})$ that
$\int_{(0,{\infty})}r^{2}e^{-w(x)r}\pi(x,dr)\le c_{16}h(x)^{-1}$ on $E$. Immediately we get $\langle \left(\int_{(0,{\infty})}r^{2}e^{-w(\cdot)r}\pi(\cdot,dr)\right)^{2},1\wedge h^{4}\rangle\le c_{16}^{2}\int_{E}(h(x)^{-2}\wedge h(x)^{2})m(dx)<{\infty}$.
Hence Assumption 4(ii) holds.

We next show Assumption 6 holds.
 Let $P^{h}_{t}$ be the $h$-transformed semigroup from $P^{\beta}_{t}$ given by \eqref{e:1.8}.
 $P^{h}_{t}$ admits an integral kernel $p_{E}^{h}(t,x,y)$ with respect to the measure $\widetilde{m}(dx):=h(x)^{2}m(dx)$, which
is related to $p_{E}^{\beta}(t,x,y)$ by
$$p_{E}^{h}(t,x,y)=e^{\lambda_{1}t}\frac{p_{E}^{\beta}(t,x,y)}{h(x)h(y)}\quad\mbox{ for }x,y\in E \mbox{ and }t\ge 0.$$
Define $\widetilde{a}_{t}(x):=p_{E}^{h}(t,x,x)$ for $x\in E$ and $t>0$. Clearly by \eqref{8.17} and \eqref{8.18}, we have
$\sup_{x\in E}\widetilde{a}_{t}(x)<{\infty}$ for every $t>0$.
 Let $e_{\beta}(t):=\int_{0}^{t}\beta(\xi_{s})ds$. Suppose $f\in C_{0}(E)$.
 Then for
 any given ${\eps}>0$, there is $\delta>0$ so that $|f(x)-f(y)|<{\eps}$ whenever $|x-y|<\delta$. By the two-sided estimate of $p_{E}^{\beta}(t,x,y)$, for
 sufficiently small
 $t\in (0,1]$,
\begin{eqnarray}
\sup_{x\in E}\left|P^{h}_{t}f(x)-f(x)\right|&=&\sup_{x\in E}\frac{\left|e^{\lambda_{1}t}P^{\beta}_{t}(hf)(x)-h(x)f(x)\right|}{h(x)}\nonumber\\
&=&\sup_{x\in E}e^{\lambda_{1}t}\frac{\left|\Pi_{x}\left[e_{\beta}(t)h(\xi_{t})(f(\xi_{t})-f(\xi_{0}))\right]\right|}{h(x)}\nonumber\\
&\le&{\eps}+\sup_{x\in E}e^{\lambda_{1}t}\frac{\left|\Pi_{x}\left[e_{\beta}(t)h(\xi_{t})(f(\xi_{t})-f(\xi_{0}));|\xi_{t}-\xi_{0}|\ge\delta\right]\right|}{h(x)}\nonumber\\
&\le&{\eps}+e^{(\lambda_{1}+\|\beta\|_{\infty})t}\|h\|_{\infty}\|f\|_{\infty}\sup_{x\in E}\frac{\Pi_{x}\left(|\xi_{t}-\xi_{0}|\ge\delta\right)}{h(x)}\nonumber\\
&\le&{\eps}+c_{11}\sup_{x\in E}\frac{\delta_{E}(x)}{1\wedge \delta_{E}(x)}\int_{y\in E:\ |y-x|\ge \delta}t^{-(d+1)/2}\exp\left(-c_{4}|x-y|^{2}/t\right)dy\nonumber\\
&\le&{\eps}+c_{12}(1+\mathrm{diam}(E))\int_{|z|\ge\delta}t^{-(d+1)/2}e^{-c_{4}|z|^{2}/t}dz\nonumber\\
&=&{\eps}+c_{12}(1+\mathrm{diam}(E))t^{-1/2}\int_{\delta t^{-1/2}}^{{\infty}}r^{d-1}e^{-c_{4}r^{2}}dr.\nonumber
\end{eqnarray}
It follows that $\lim_{t\to 0}\|P^{h}_{t}f-f\|_{\infty}=0$ for all $f\in C_{0}(E)$.

It remains to prove that Assumption 5 holds. Fix $\phi\in \mathcal{B}^{+}_{b}(E)$, $\mu\in\mc$ and $\sigma>0$. For $m,n\in\mathbb{N}$, we have
\begin{equation}
e^{\lambda_{1}(m+n)\sigma}\langle \frac{h}{w}\phi,Z_{(m+n)\sigma}\rangle-\langle \phi,h^{2}\rangle W^{h/w}_{\infty}(Z)=I(m,n)+II(m,n)+III(n),\label{8.19}
\end{equation}
where $I(m,n):=e^{\lambda_{1}(m+n)\sigma}\langle \frac{h}{w}\phi,Z_{(m+n)\sigma}\rangle-\pp_{\mu}\left(\left.e^{\lambda_{1}(m+n)\sigma}\langle \frac{h}{w}\phi,Z_{(m+n)\sigma}\rangle\right|\mathcal{F}_{n\sigma}\right)$,
$$II(m,n):=\pp_{\mu}
\left(
\left.
e^{\lambda_{1}(m+n)\sigma}\langle \frac{h}{w}\phi,Z_{(m+n)\sigma}\rangle
\right|\mathcal{F}_{n\sigma}
\right)
-\langle \phi,h^{2}\rangle W^{h/w}_{n\sigma}(Z),
$$
and $III(n):=\langle \phi,h^{2}\rangle \left(W^{h/w}_{n\sigma}(Z)-W^{h/w}_{\infty}(Z)\right)$.
Note that by the Markov property of $Z$,
\begin{eqnarray}
\pp_{\mu}
\left(
\left.
e^{\lambda_{1}n\sigma}\langle \frac{h}{w}\phi,Z_{(m+n)\sigma}\rangle
\right|\mathcal{F}_{n\sigma}
\right)
&=&e^{\lambda_{1}(m+n)\sigma}\pp_{\cdot,Z_{n\sigma}}\left
(\langle \frac{h}{w}\phi,Z_{m\sigma}\rangle\right)\nonumber\\
&=&e^{\lambda_{1}(m+n)\sigma}\langle \frac{1}{w}P^{\beta}_{m\sigma}(h\phi),Z_{n\sigma}\rangle=e^{\lambda_{1}n\sigma}\langle\frac{h}{w}P^{h}_{m\sigma}\phi,Z_{n\sigma}\rangle.\nonumber
\end{eqnarray}
Thus $II(m,n)=e^{\lambda_{1}n\sigma}\langle\frac{h}{w}P^{h}_{m\sigma}g,Z_{n\sigma}\rangle$, where $g(x):=\phi(x)-\langle \phi,h^{2}\rangle$. Since $\langle g,h^{2}\rangle=0$, $\left|P^{h}_{m\sigma}g(x)\right|\le e^{-\lambda_{h}(m\sigma-\frac{1}{2})}\tilde{a}_{1}(x)^{1/2}\|g\|_{L^{2}(E,h^{2}m)}$ for $m\sigma>1/2$. Here $\lambda_{h}>0$ denotes the spectral gap in $\sigma(\mathcal{L}^{ (\beta)})$.  Thus we have for $m\sigma>1/2$
$$|II(m,n)|=\left|e^{\lambda_{1}n\sigma}\langle \frac{h}{w}P^{h}_{m\sigma}g,Z_{n\sigma}\rangle\right|\le c_{14}e^{-\lambda_{h}m\sigma}\sup_{x\in E}\tilde{a}_{1}(x)^{1/2}W^{h/w}_{n\sigma}(Z).$$
It is shown in the proof of the equivalence between Assumption 5 and Assumption 5' that for every $m\in\mathbb{N}$,
$\lim_{n\to{\infty}}I(m,n)=0$ $\pp_{\mu}$-a.s. Thus by letting $n\to {\infty}$ and then $m\to {\infty}$ in
\eqref{8.19},
we obtain that
$$\lim_{n\to{\infty}}e^{\lambda_{1}n\sigma}\langle \frac{h}{w}\phi,Z_{n\sigma}\rangle -\langle \phi,h^{2}\rangle W^{h/w}(Z)=0\quad\pp_{\mu}\mbox{-a.s.}$$

\medskip

\noindent\textbf{Example \ref{E:8.4}:}
Let $\lambda_{h}:=\lambda_{2}-\lambda_{1}$ be the spectral gap in $\sigma(\mathcal{L}^{(\beta ) })$,
where$
\lambda_{2}:=\inf \left\{\lambda\in\sigma(- \mathcal{L}^{(\beta ) }):  \lambda\not= \lambda_{1}
\right\}.$
 Since $P^{\beta}_{t}h(x)=e^{-\lambda_{1}t}h(x)$, it follows by property (ii) and H\"{o}lder inequality that $h\in L^{4}(E,m)$.
Let $g_{0}(\theta):=\log^{*}\theta/\theta$ for $\theta\in (0,+\infty)$.
Then
\begin{eqnarray}
\langle \int_{(0,+\infty)}r\log^{*}(rh(\cdot))\pi(\cdot,dr),h^{2}\rangle
&=&\langle \int_{(0,+\infty)}g_{0}(rh(\cdot))r^{2}\pi(\cdot,dr),h^{3}\rangle\nonumber\\
&\le&\|g_{0}\|_{\infty}\langle \int_{(0,+\infty)}r^{2}\pi(\cdot,dr),h^{3}\rangle\nonumber\\
&<&+\infty.\nonumber
\end{eqnarray}
The last inequality comes from the fact that $h\in L^{4}(E,m)$ and condition \eqref{8.1}.
Thus we have by Theorem \ref{them2} and Remark \ref{rm for them2} that $\p_{\delta_{x}}\left(W^{h}_{\infty}(X)>0\right)>0$ for each $x\in E$.
Hence
$$w(x)\ge -\log\left(1-\p_{\delta_{x}}\left(\limsup_{t\to+\infty}\langle 1,X_{t}\rangle >0\right)\right)\ge-\log\left(1-\p_{\delta_{x}}\left(W^{h}_{\infty}(X)>0\right)\right)>0.$$
In view of Remark \ref{rm1}, $w(x)$ is a bounded function on $E$ under our assumptions.
Using the boundedness of $w$ and \eqref{8.1}, it is easy to verify that
Assumptions 1-4 and Assumption 6
are satisfied by this example. Next we will show that
Assumption 5' is also satisfied,
that is, for all $\mu\in\mc$, $\phi\in\mathcal{B}^{+}_{b}(E)$, $\sigma>0$ and some $m\in \mathbb{N}$,
\begin{equation}
\lim_{n\to+\infty}e^{\lambda_{1}n\sigma}\langle \frac{h}{w}P^{h}_{m\sigma}\phi,Z_{n\sigma}\rangle=\langle \phi ,h^{2}\rangle W^{h/w}_{\infty}(Z)\quad\pp_{\mu}
\mbox{-a.s.}
\label{8.4}
\end{equation}
It is showed in the proof of Lemma \ref{lem4.2} that
$$e^{\lambda_{1}n\sigma}\langle \frac{h}{w}\,P^{h}_{m\sigma}\phi,Z_{n\sigma}\rangle =e^{\lambda_{1}(m+n)\sigma}
\sum_{i=1}^{N_{n\sigma}}\pp_{\mu}\left(\langle \phi h,I^{i,n\sigma}_{m\sigma}\rangle\,|\,\mathcal{F}_{n\sigma}\right)
+e^{\lambda_{1}n\sigma}\langle \frac{h}{w}\,\theta^{*}_{\phi h}(m\sigma, \cdot),Z_{n\sigma}\rangle.$$
Thus by Lemma \ref{lem1.1}, Lemma \ref{lem4.1} and Theorem \ref{them1}, under Assumptions 1-3 and
Assumption 4(i),
\eqref{8.4} is equivalent to that
\begin{equation}
\lim_{n\to+\infty}e^{\lambda_{1}(m+n)\sigma}\sum_{i=1}^{N_{n\sigma}}\langle\phi h,I^{i,n\sigma}_{m\sigma}\rangle=\langle \phi,h^{2}\rangle W^{h}_{\infty}(X)\quad\pp_{\mu}
\mbox{-a.s.}
\label{8.6}
\end{equation}
Since $e^{\lambda_{1}(m+n)\sigma}\langle \phi h, X_{(m+n)\sigma}\rangle \ge e^{\lambda_{1}(m+n)\sigma}\sum_{i=1}^{N_{n\sigma}}\langle\phi h,I^{i,n\sigma}_{m\sigma}\rangle$ and $\phi\in\mathcal{B}^{+}_{b}(E)$, \eqref{8.6} (or equivalently \eqref{8.4}) is equivalent to
\begin{equation}
\lim_{n\to+\infty}e^{\lambda_{1}n\sigma}\langle \phi h,X_{n\sigma}\rangle=\langle \phi,h^{2}\rangle W^{h}_{\infty}(X)\quad\pp_{\mu}
\mbox{-a.s.}
\label{8.7}
\end{equation}
Fix $\phi\in\mathcal{B}^{+}_{b}(E)$ and $\mu\in\mc$. Let $g(x):=\phi(x)-\langle \phi,h^{2}\rangle$ for $x\in E$.
We have by \eqref{var} that for $t>0$,
\begin{eqnarray}
\pp_{\mu}\left[\left(e^{\lambda_{1}t}\langle gh,X_{t}\rangle\right)^{2}\right]
&=&\langle P^{h}_{t}g,h\mu\rangle^{2}+\int_{0}^{t}e^{\lambda_{1}s}\langle P^{h}_{s}\left[\left(2\alpha+\int_{(0,+\infty)}y^{2}\pi(\cdot,dy)\right)h\left(P^{h}_{t-s}g\right)^{2}\right],h\mu\rangle ds\nonumber\\
&=:&I(t)+II(t).\nonumber
\end{eqnarray}
Since $g\in L^{2}(E,\widetilde{m})$ and $\langle g,h^{2}\rangle=0$, we have by \eqref{p3} for $r>t_{0}/2$ and $x\in E$,
\begin{equation}\label{8.9}
|P^{h}_{\lambda}g(x)|\le e^{-\lambda_{h}(r-\frac{t_{0}}{2})}\widetilde{a}_{t_{0}}(x)^{1/2}\|g\|_{L^{2}(E,\widetilde{m})}=:c_{1}e^{-\lambda_{h}r}\widetilde{a}_{t_{0}}(x)^{1/2}.
\end{equation}
Thus $I(t)\le c_{2}e^{-2\lambda_{h}t}$ for all $t>t_{0}/2$ and some $c_{2}>0$. On the other hand, for $r\in (0,t_{0}/2]$ and $x\in E$,
\begin{equation}\label{8.10}
|P^{h}_{\lambda}g(x)|\le \|g\|_{\infty}\le \|g\|_{\infty}e^{\lambda_{h}(\frac{t_{0}}{2}-r)}=:c_{3}e^{-\lambda_{h}r}.
\end{equation}
Let ${\eps}\in (0,-\lambda_{1}\wedge 2\lambda_{h})$, we have by \eqref{8.9} and \eqref{8.10} that
\begin{eqnarray}
II(t)&\le&\|2\alpha+\int_{(0,+\infty)}y^{2}\pi(\cdot,dy)\|_{\infty}\int_{0}^{t}e^{\lambda_{1}s}\langle P^{h}_{s}\left(h\left(P^{h}_{t-s}g\right)^{2}\right),h\mu \rangle ds\nonumber\\
&\le&c_{4}\int_{0}^{t}e^{\lambda_{1}s-2\lambda_{h}(t-s)}\langle P^{h}_{s}\left(h\left(1\vee \widetilde{a}_{t_{0}}\right)\right),h\mu \rangle ds\nonumber\\
&\le&c_{4}e^{-{\eps} t}\int_{0}^{+\infty}e^{(\lambda_{1}+{\eps})s}\langle P^{h}_{s}\left(h\left(1\vee \widetilde{a}_{t_{0}}\right)\right),h\mu \rangle ds.\label{8.11}
\end{eqnarray}
Since $h\in L^{4}(E,m)$, it follows by property (ii) that $h\left(1\vee \widetilde{a}_{t_{0}}\right)\in L^{2}(E,\widetilde{m})$. Thus by \eqref{p3}, $\langle P^{h}_{s}\left(h\left(1\vee \widetilde{a}_{t_{0}}\right)\right),h\mu\rangle$ is bounded from above for $s$ sufficiently large. Therefore the integral in the right hand side of \eqref{8.11} is finite. It implies that $\sum_{n=1}^{+\infty}\pp_{\mu}\left[\left(e^{\lambda_{1}n\sigma}\langle gh,X_{n\sigma}\rangle\right)^{2}\right]<+\infty$, and hence we obtain \eqref{8.7} by Borel-Cantelli Lemma.

\medskip

\noindent\textbf{Example \ref{E:8.5}:}
By a similar argument as that for \cite[IV.5]{CMS}, we have
\begin{equation}
h(x)\ge c_{0}|x|^{-d-\alpha}
\mbox{ for }|x|\ge 1 ,
\label{7.6}
\end{equation}
for some positive constant $c_{0}$.
Let $\lambda_{2}$ be the second bottom of the spectrum of $\sigma(\mathcal{E}^{(\beta ) })$, that is $\lambda_{2}:=\inf\{\mathcal{E}^{(\beta ) }(u,u):\ u\in \mathcal{F}\mbox{ with }\int_{\R^{d}}u(x)h(x)dx=0\mbox{ and }\int_{\R^{d}}u(x)^{2}dx=1\}.$ Then the spectral gap $\lambda_{h}:=\lambda_{2}-\lambda_{1}>0$.
Recall that the $h$-transformed semigroup $P^{h}_{t}$ is an $\widetilde{m}$-symmetric semigroup with $\widetilde{m}(dx)=h(x)^{2}dx$.  $P^{h}_{t}$ admits an integral kernel $p^{h}(t,x,y)$, given by \eqref{density-trans}, with respect to the measure $\widetilde{m}(dx)$. Let $p^{\beta}(t,x,y)$ be the integral kernel of $P^{\beta}_{t}$ with respect to the Lebesgue measure, which satisfies the estimates in \eqref{density-FK}.
For any $f\in C_{c}(\R^{d})$, $t>0$ and $x\in \R^{d}$, we have by \eqref{7.1} and \eqref{density-FK} that
\begin{eqnarray}
\left|P^{h}_{t}f(x)-f(x)\right|
&\le&\int_{\R^{d}}p^{h}(t,x,y)\left|f(y)-f(x)\right|\widetilde{m}(dy)\nonumber\\
&=&\frac{e^{\lambda_{1}t}}{h(x)}\int_{\R^{d}}p^{\beta}(t,x,y)h(y)\left|f(y)-f(x)\right|dy\nonumber\\
&\le&
c_{1}
\frac{e^{(\lambda_{1}+\|\beta\|_{\infty})t}}{h(x)}\int_{\R^{d}}\left(t^{-d/\alpha}\wedge\frac{t}{|x-y|^{d+\alpha}}\right)
h(y)|f(y)-f(x)|dy\label{9.13}
\end{eqnarray}
Let $f_{t}(x):=\frac{1}{h(x)}\int_{\R^{d}}\left(t^{-d/\alpha}\wedge\frac{t}{|x-y|^{d+\alpha}}\right)h(y)|f(y)-f(x)|dy$.
Suppose $\mbox{supp}(f)\subset B(0,R)$ for some $R>1$.
For $x\in B(0,2R)$ and $z\in\R^{d}$, $h(z)/h(x)\le\|h\|_{\infty}/\inf_{y\in B(0,2R)}h(y)<{\infty}$, and
we have
\begin{eqnarray}
f_{t}(x)&\le&c_{2}\int_{\R^{d}}\left(t^{-d/\alpha}\wedge\frac{t}{|x-y|^{d+\alpha}}\right)|f(y)-f(x)|dy\nonumber\\
&=&c_{2}\int_{\R^{d}}\left(t^{-d/\alpha}\wedge\frac{t}{|z|^{d+\alpha}}\right)
|f(x+z)-f(x)|dz.
\label{7.3}
\end{eqnarray}
Since $f\in C_{c}(\R^{d})$ is uniformly continuous on $\R^{d}$, for any ${\eps}>0$, there exists $\delta>0$, such that $|f(z_{1})-f(z_{2})|\le {\eps}$ whenever $|z_{1}-z_{2}|<\delta$. Thus by \eqref{7.3},
\begin{eqnarray}
f_{t}(x)&\le&c_{2}{\eps}\int_{|z|<\delta}t^{-d/\alpha}\wedge\frac{t}{|z|^{d+\alpha}}\,dz+2c_{2}\|f\|_{\infty}
\int_{|z|\ge \delta}t^{-d/\alpha}\wedge\frac{t}{|z|^{d+\alpha}}\,dz\nonumber\\
&=&c_{2}{\eps}\int_{|u|<\delta t^{-1/\alpha}}1\wedge\frac{1}{|u|^{d+\alpha}}\,du+2c_{2}\|f\|_{\infty}
\int_{|u|\ge \delta t^{-1/\alpha}}1\wedge\frac{1}{|u|^{d+\alpha}}\,du\label{7.4}
\end{eqnarray}
Since $\int_{\R^{d}}1\wedge\frac{1}{|u|^{d+\alpha}}\,du<{\infty}$, by letting $t\to 0$ and ${\eps}\to 0$ in \eqref{7.4}, we get $\sup_{x\in B(0,2R)}f_{t}(x)\to 0$.
On the other hand, for $x\not\in B(0,2R)$, by the fact $\mbox{supp}(f)\subset B(0,R)$, the boundedness of $h$ and \eqref{7.6}, we have
\begin{eqnarray}
f_{t}(x)&=&\frac{1}{h(x)}\int_{|y|\le R}\left(t^{-d/\alpha}\wedge\frac{t}{|x-y|^{d+\alpha}}\right)h(y)|f(y)|dy\nonumber\\
&\le&c_{3}\int_{|y|\le R}\left(t^{-d/\alpha}\wedge\frac{t}{|x-y|^{d+\alpha}}\right)|x|^{d+\alpha}dy\nonumber\\
&\le&c_{3}\int_{|y|\le R}\frac{t|x|^{d+\alpha}}{(|x|-|y|)^{d+\alpha}}\,dy\nonumber\\
&\le&c_{3}\int_{|y|\le R}\frac{t|x|^{d+\alpha}}{(|x|-R)^{d+\alpha}}\,dy\nonumber\\
&\le&2^{d+\alpha}c_{3}t.\nonumber
\end{eqnarray}
Thus $\sup_{x\not\in B(0,2R)}f_{t}(x)\to 0$ as $t\to 0$. Hence by \eqref{9.13} we conclude that
\begin{equation}
\lim_{t\to 0}\|P^{h}_{t}f-f\|_{\infty}=0\label{7.7}
\end{equation}
for all $f\in C_{c}(\R^{d})$. Since $C_{c}(\R^{d})$ is dense in $(C_{0}(\R^{d}),\|\cdot\|_{\infty})$, \eqref{7.7} is true for all $f\in C_{0}(\R^{d})$.

Recall that $w(x)=-\log \p_{\delta_{x}}\left(\exists t\ge 0,\ \langle 1,X_{t}\rangle =0\right)$.
The argument in Remark \ref{rm1} shows that
$w(x)$ is bounded on $E$ and so Condition 1 is satisfied.
Let $W^{h}_{t}(X):=e^{\lambda_{1}t}\langle h,X_{t}\rangle$. Then $W^{h}_{t}(X)$ is a nonnegative $\p_{\mu}$-martingale with respect to $\mathcal{F}_{t}:=\sigma\{X_{s}:s\in [0,t]\}$ for all $\mu\in \mf$. Let $W^{h}_{\infty}(X):=\lim_{t\to{\infty}}W^{h}_{t}(X)$. Note that by \eqref{var},
\begin{eqnarray}
\p_{\delta_{x}}\left(W^{h}_{t}(X)^{2}\right)&=&h(x)^{2}+e^{2\lambda_{1}t}\int_{0}^{t}P^{\beta}_{s}
\left[\widehat{\alpha}\left(P^{\beta}_{t-s}h\right)^{2}\right](x)ds\nonumber\\
&=&h(x)^{2}+\int_{0}^{t}e^{2\lambda_{1}s}P^{\beta}_{s}
\left[\widehat{\alpha}h^{2}\right](x)ds\nonumber\\
&\le&h(x)^{2}+\|\widehat{\alpha}h\|_{\infty}\int_{0}^{t}e^{2\lambda_{1}s}P^{\beta}_{s}h(x)ds\nonumber\\
&=&h(x)^{2}+\|\widehat{\alpha}h\|_{\infty}h(x)\int_{0}^{t}e^{\lambda_{1}s}ds\nonumber\\
&\le&h(x)^{2}+c_{4}h(x).\label{8.14}
\end{eqnarray}
Here $\widehat{\alpha}(x):=2\alpha(x)+\int_{(0,{\infty})}y^{2}\pi(x,dy)$ for $x\in E$ and $c_{4}$ is a positive constant independent of $t$.
Thus
for any $\mu\in {\cal M}^h_F(E)$,
$\{W^{h}_{t}(X):\ t\ge 0\}$ is
an
$L^{2}$-bounded nonnegative martingale and hence
$W^{h}_{t}(X)$ converges to $W^{h}_{\infty}(X)$ $\p_{\mu}$-a.s. and in $L^{2}(P_{\mu})$.
In particular $\p_{\mu}\left(W^{h}_{\infty}(X)>0\right)>0$.
Thus $\p_{\delta_{x}}\left(\exists t>0,\ \langle 1,X_{t}\rangle=0\right)\le 1-\p_{\delta_{x}}\left(W^{h}_{\infty}(X)>0\right)<1$, and so
$w(x)>0$ for all $x\in\R^{d}$.
Furthermore,
By Cauchy-Schwartz inequality and \eqref{8.14},
we have
\begin{equation}
\p_{\delta_{x}}\left(W^{h}_{\infty}(X)>0\right)
\ge \frac{ \left( \p_{\delta_{x}}\left[W^{h}_{\infty}(X)\right] \right)^{2}}
{\p_{\delta_{x}}\left[ W^{h}_{\infty}(X)^{2}\right]}=
\frac{h(x)^{2}}{\lim_{t\to{\infty}}\p_{\delta_{x}}\left[ W^{h}_{t}(X)^{2}\right]}
\ge \frac{h(x)}{h(x)+c_{4}} .\nonumber
\end{equation}
Thus we get $$w(x)\ge -\log \left(1-\p_{\delta_{x}}\left(W^{h}_{\infty}(X)>0\right)\right)\ge \log(h(x)+c_{4})-\log c_{4}.$$
Using this, \eqref{7.6} and the fact that $h$ is bounded away from 0 and ${\infty}$ on compact sets, we conclude that the function $h(x)/w(x)$ is bounded from above on $\R^{d}$.
We see from the above arguments that Assumptions 1-4 and Assumption 6 are satisfied by this example.

Finally we will show that
Assumption 5 is also satisfied.
Let $Z$ be the skeleton process. We need to show that for all $\mu\in \mc$, $\sigma>0$ and $f\in\mathcal{B}^{+}(E)$ with $\frac{fw}{h}$ bounded,
\begin{equation}
\lim_{n\to{\infty}}e^{\lambda_{1}n\sigma}\langle f,Z_{n\sigma}\rangle =\langle f,wh\rangle W^{h/w}_{\infty}(Z)\quad\pp_{\mu}\mbox{-a.s.}\label{8.20}
\end{equation}
Fix $f\in\mathcal{B}^{+}(E )$ with $fw/h$ bounded. Let $g(x):=f(x)-\langle f,wh\rangle \frac{h}{w}(x)=\left(\frac{fw}{h}(x)-\langle f,wh\rangle\right)\frac{h}{w}(x)$ for $x\in E $.
Through the same argument as in Example \ref{E:8.1},
to prove \eqref{8.20}, it suffices to prove that for every $x\in E$,
\begin{equation}
\pp_{\cdot,\delta_{x}}\left(\lim_{n\to{\infty}}e^{\lambda_{1}n\sigma}\langle g,Z_{n\sigma}\rangle=0\right)=1.\label{8.21}
\end{equation}
We observe that $g$, $\frac{gw}{h}$ are bounded functions on $\R^{d}$ by the boundedness of $\frac{h}{w}$ and $\frac{fw}{h}$. For every $x\in \R^{d}$,
$$\pp_{\cdot,\delta_{x}}\left[\left(e^{\lambda_{1}t}\langle g,Z_{t}\rangle\right)^{2}\right]=\frac{h(x)}{w(x)}\left(I(t,x)+II(t,x)\right),$$
where $I(t,x):=e^{\lambda_{1}t}P^{h}_{t}\left(\frac{w}{h}g^{2}\right)(x)$,
$II(t,x):=\int_{0}^{t}e^{\lambda_{1}s}P^{h}_{s}\left[(2\alpha+b)\frac{h}{w}P^{h}_{t-s}\left(\frac{w}{h}g\right)^{2}\right](x)ds$.
For $I$, we have
\begin{equation}
I(t,x)\le\|\frac{w}{h}g\|_{\infty}\|g\|_{\infty}e^{\lambda_{1}t}.\label{7.12}
\end{equation}
For $II$, we have
\begin{equation}
II(t,x)\le \|2\alpha+b\|_{\infty}\|\frac{h}{w}\|_{\infty}\|\frac{w}{h}g\|_{\infty}\int_{0}^{t}e^{\lambda_{1}s}P^{h}_{s}\left[
\left|P^{h}_{t-s}\left(\frac{w}{h}g\right)\right|\right](x)ds.\label{7.13}
\end{equation}
Since $\langle \frac{w}{h}g,h^{2}\rangle=0$,
by \eqref{p1}
and H\"{o}lder inequality for any $t>s>0$
\begin{eqnarray}
P^{h}_{s}\left[\left|P^{h}_{t-s}\left(\frac{w}{h}g\right)\right|\right](x)
&=&\int_{\R^{d}}p^{h}(s,x,y)\left|P^{h}_{t-s}\left(\frac{w}{h}g\right)\right|\widetilde{m}(dy)\nonumber\\
&\le&\left(\int_{\R^{d}}p^{h}(s,x,y)^{2}\widetilde{m}(dy)\right)^{1/2}\|P^{h}_{t-s}\left(\frac{w}{h}g\right)
\|_{L^{2}(\R^{d},\widetilde{m})}\nonumber\\
&\le&e^{-\lambda_{h}(t-s)}p^{h}(2s,x,x)^{1/2}\|\frac{w}{h}g\|_{L^{2}(\R^{d},\widetilde{m})}\label{7.14}
\end{eqnarray}
Let ${\eps}\in (0,-\lambda_{1}\wedge \lambda_{h})$. Note that by \eqref{p3}, for each $x$, $p^{h}(s,x,x)$ is bounded from above when $s$ is suffciently large. This together with \eqref{7.13} and \eqref{7.14} yield that for every $x\in\R^{d}$,
\begin{equation}
II(t,x) \le c_{5}\int_{0}^{t}e^{\lambda_{1}s-\lambda_{h}(t-s)}ds
 \le c_{5}e^{-{\eps} t}\int_{0}^{t}e^{(\lambda_{1}+{\eps})s}ds
 \le c_{6}e^{-{\eps} t}.\label{7.15}
\end{equation}
Thus by \eqref{7.12} and \eqref{7.15}, we get $\sum_{n=1}^{{\infty}}\pp_{\cdot,\delta_{x}}\left[\left(e^{\lambda_{1}n\sigma}\langle g,Z_{n\sigma}\rangle\right)^{2}\right]<{\infty}$. Consequently, \eqref{8.21} follows by Borel-Cantelli lemma.

\bigskip

\noindent{\bf Acknowledgement.} We thank Zenghu Li for helpful discussions on Kuznetsov
measure for superprocesses.

\medskip

{\bf Zhen-Qing Chen}

Department of Mathematics, University of Washington, Seattle,
WA 98195, USA

E-mail: zqchen@uw.edu

\medskip
{\bf Yan-Xia Ren}

LMAM School of Mathematical Sciences \& Center for
Statistical Science, Peking
University,

Beijing, 100871, P.R. China.

E-mail: yxren@math.pku.edu.cn

\medskip
{\bf Ting Yang}

School of Mathematics and Statistics, Beijing Institute of Technology, Beijing, 100081, P.R.China;

Beijing Key Laboratory on MCAACI, Beijing, 100081, P.R. China.

Email: yangt@bit.edu.cn

\end{document}